\documentclass[10pt,final]{article}
\usepackage[left=2.54cm,top=2.54cm,right=2.54cm,bottom=2.54cm]{geometry}
\usepackage{fancyhdr}
\usepackage{afterpage}
\usepackage{setspace}
\usepackage{bibspacing}
\newfont{\toto}{msbm10 at 12 pt}
\newfont{\ithd}{cmr9}

\usepackage{amsmath,amsthm,amsfonts,amssymb}
\usepackage[pdftex]{graphicx}
\usepackage[T1]{fontenc}

\graphicspath{{}}

\usepackage{setspace}
\usepackage{titlesec}
\usepackage{wallpaper}
\usepackage{emptypage}
\usepackage[titletoc]{appendix}

\usepackage[utf8]{inputenc}
\usepackage[ngerman,english]{babel}
\usepackage[T1]{fontenc}
\usepackage{marvosym} 

\usepackage{graphicx}
\usepackage{caption}
\usepackage{subcaption}
\usepackage{xcolor}
\usepackage{tabularx}

\setkeys{Gin}{draft=false}

\usepackage[final]{animate}

\usepackage{mathtools}
\usepackage{amssymb}
\usepackage{amsfonts}
\usepackage{amsbsy} 
\usepackage{amsopn} 
\usepackage{amstext} 
\usepackage{amsxtra} 
\usepackage{bbm} 
\usepackage{dsfont} 

\usepackage{cancel} 
\usepackage{amscd} 
\usepackage[all]{xy} 
\usepackage{esvect} 
\usepackage{esint} 
\usepackage{empheq} 

\usepackage{eso-pic,picture}
\usepackage{lscape}
\usepackage{blindtext}
\usepackage{listings}
\usepackage{url}
\usepackage{textcomp} 
\usepackage{exscale} 
\usepackage{xspace} 
\usepackage{scalerel} 
\usepackage{stackengine} 
\usepackage{hhline}
\usepackage{lineno}
\usepackage{fancyhdr}
\usepackage{longtable}

\usepackage{hyphenat} 
\usepackage{testhyphens} 

\input ushyphex
\PassOptionsToPackage{svgnames,dvipsnames}{xcolor}
\usepackage{tcolorbox}
\tcbuselibrary{skins,breakable}

\mathcode`@="8000
{\catcode`\@=\active\gdef@{\mkern1mu}}

\usepackage{lmodern}
\usepackage{microtype}

\usepackage[numbers]{natbib}

\newcounter{thmctr}
\newtheorem{theorem}[thmctr]{Theorem}
\newtheorem{lemma}[thmctr]{Lemma}
\newtheorem{corollary}[thmctr]{Corollary}
\newtheorem{remark}[thmctr]{Remark}
\newtheorem{proof}{Proof}

\addto\captionsenglish{}

\usepackage[notref,notcite]{showkeys} 
\usepackage[mode=multiuser]{fixme} 
\fxusetheme{color}

\RequirePackage[l2tabu]{nag}

\usepackage{varioref} 
\usepackage{hyperref} 
\usepackage[capitalise]{cleveref} 

\let\ccref\cref
\renewcommand{\cref}[1]{\mbox{\ccref{#1}}}
\crefname{appendix}{}{}
\crefname{listing}{Algorithm}{Algorithms}
\crefname{equation}{}{}
\crefname{table}{Tab.~}{Tabs.~}
\crefname{lemma}{Lemma~}{Lemmata~}

\renewcommand{\l}{\ensuremath{\left(}} 
\renewcommand{\r}{\ensuremath{\right)}} 
\newcommand{\lb}{\ensuremath{\left[}}
\newcommand{\rb}{\ensuremath{\right]}}
\newcommand{\lc}{\ensuremath{\left\lbrace}}
\newcommand{\rc}{\ensuremath{\right\rbrace}}

\newcommand{\transp}{\ensuremath{^{\textsf T}}}

\newcommand\restr[2]{{
\left.\kern-\nulldelimiterspace 
#1 
\vphantom{|} 
\right|_{#2} 
}}

\newcommand{\dx}{\ensuremath{\,\mathrm{d}\vec x}}

\newcommand{\ds}{\ensuremath{\,\mathrm{d} s}}

\renewcommand{\d}{\,\mathrm{d}} 

\renewcommand{\div}{\ensuremath{\nabla \cdot}}

\newcommand{\R}{\ensuremath{\mathbb{R}}}

\renewcommand{\phi}{\ensuremath{\varphi}}

\renewcommand{\epsilon}{\ensuremath{\varepsilon}}

\renewcommand{\vec}[1]{\ensuremath{\boldsymbol{#1}}}

\newcommand{\ie}{, \mbox{i.\,e.,}~}
\newcommand{\eg}{\mbox{e.\,g.,}~}
\newcommand{\etal}{\mbox{et al.\ }}

\newcommand{\wrt}{\mbox{w.\,r.\,t.}~}
\newcommand{\cf}{{cf.\ }~}

\definecolor{myred}{RGB}{182,59,87}
\definecolor{myblue}{RGB}{116,135,197}
\definecolor{mygreen}{RGB}{52,152,126}
\definecolor{myorange}{RGB}{227,149,87}

\definecolor{green}{rgb}{0.0,0.6,0.0}

\definecolor{brightgray}{RGB}{198,198,198}
\definecolor{darkgray}{RGB}{70,66,66}
\definecolor{redbrown}{RGB}{164,93,72}

\definecolor{MoSTBrightBlue}{RGB}{0,129,200}
\definecolor{MoSTDarkBlue}{RGB}{0,78,146}
\definecolor{MoSTVeryDarkBlue}{RGB}{0,58,107}

\allowdisplaybreaks 

\newcommand*\reallywidehat[1]{%
\savestack{\tmpbox}{\stretchto{%
\scaleto{%
\scalerel*[\widthof{\ensuremath{#1}}]{\kern-.6pt\bigwedge\kern-.6pt}%
{\rule[-\textheight/2]{1ex}{\textheight}}
}{\textheight}%
}{0.5ex}}%
\stackon[1pt]{#1}{\tmpbox}%
}

\makeatletter
\newcommand*{\coloneqq}{\ensuremath{\mathrel{\rlap{%
\raisebox{0.38ex}{$\cdot$}}\raisebox{-0.38ex}{$\cdot$}}=}}
\makeatother

\makeatletter
\newcommand*{\eqqcolon}{\ensuremath{=\mathrel{\rlap{%
\raisebox{0.38ex}{$\cdot$}}\raisebox{-0.38ex}{$\cdot$}}}}
\makeatother

\definecolor{keywordcolor}{rgb}{0,0.25,0.5}
\definecolor{commentcolor}{rgb}{0.2,0.5,0.2}
\definecolor{stringcolor}{rgb}{0.5,0.5,0.2}
\newcommand{\link}[2]{{\color{blue}\href{#1}{\texttt{#2}}}}

\newcommand{\doi}[1]{\textrm{\textsc{doi:}} \link{http://dx.doi.org/#1}{\nolinkurl{#1}}}

\let\myurl\url
\renewcommand{\url}[1]{\textcolor{blue}{\myurl{#1}}}

\addto{\captionsenglish}{}

\graphicspath{{../img/hyperbolic/}{../img/afc/}{../img/swe/}
{../img/theory/}{../img/discontinuous/}}

\newcommand{\f}{\ensuremath{\mathbf f}}

\newcommand{\n}{\ensuremath{\vec n}}
\renewcommand{\vv}{\ensuremath{\vec v}}
\newcommand{\x}{\ensuremath{\vec x}}

\newcommand{\id}{\ensuremath{\mathcal I}}

\newcommand{\mij}{\ensuremath{m_{ij}}}
\newcommand{\cij}{\ensuremath{\mathbf c_{ij}}}
\newcommand{\dij}{\ensuremath{d_{ij}}}
\newcommand{\fij}{\ensuremath{f_{ij}}}
\newcommand{\gij}{\ensuremath{g_{ij}}}
\newcommand{\uij}{\ensuremath{\bar u_{ij}}}
\newcommand{\lij}{\ensuremath{\alpha_{ij}}}
\newcommand{\hij}{\ensuremath{\bar h_{ij}}}
\newcommand{\hvij}{\ensuremath{\overline{(h\vv)}_{ij}}}
\newcommand{\afcsum}{\ensuremath{\sum_{j\in \mathcal N_i\setminus\{i\}}}}
\newcommand{\afcijtext}{\smash{$i\in \{1,\ldots,N\}$}, \smash{$j \in \mathcal N_i \setminus \{i\}$}}
\newcommand{\afcijmath}{\ensuremath{i\in \{1,\ldots,N\},~j \in \mathcal N_i \setminus \{i\}}}

\newcommand{\fhv}{\ensuremath{\f^{h\vv}}}
\newcommand{\spacingmin}{\phantom{a}\hspace{-1.8mm}}
\newcommand{\spacingmax}{\phantom{i}\hspace{-1.8mm}}

\title{\textbf{Bound-preserving and entropy-stable algebraic flux correction schemes for the shallow water equations with topography}}

\author{H. Hajduk$^{*}$ and D. Kuzmin$^{*}$\\
Corresponding author: hennes.hajduk@math.tu-dortmund.de\\\\
$^{*}$ TU Dortmund University, Institute of Applied Mathematics (LS III)\\ Vogelpothsweg 87, D-44227 Dortmund, Germany
}
\date{}

\begin{document}

\maketitle
\afterpage{\fancyhead{}}

\centerline{
\begin{minipage}[t]{150mm}
\textbf{Abstract:} 
A well-designed numerical method for the shallow water equations (SWE) should ensure well-balancedness, nonnegativity of water heights, and entropy stability.
For a continuous finite element discretization of a nonlinear hyperbolic system without source terms, positivity preservation and entropy stability can be enforced using the framework of algebraic flux correction (AFC).
In this work, we develop a well-balanced AFC scheme for the SWE system including a topography source term.
Our method preserves the lake at rest equilibrium up to machine precision.
The low-order version represents a generalization of existing finite volume approaches to the finite element setting.
The high-order extension is equipped with a property-preserving flux limiter.
Nonnegativity of water heights is guaranteed under a standard CFL condition.
Moreover, the flux-corrected space discretization satisfies a semi-discrete entropy inequality.
New algorithms are proposed for realistic simulation of wetting and drying processes.
Numerical examples for well-known benchmarks are presented to evaluate the performance of the scheme.
\end{minipage}
}
\vskip0.5cm

\section{Introduction}

In this paper, we develop a new property-preserving method for the shallow water equations (SWE) using algebraic flux correction (AFC) tools to enforce relevant inequality constraints \cite{kuzmin2012}.
The proposed methodology is based on the monolithic convex limiting (MCL) strategy developed in \cite{kuzmin2020} for homogeneous (systems of) conservation laws.
We use continuous piecewise (multi-)linear finite elements as baseline discretization but extensions to higher-order elements and discontinuous Galerkin approaches are expected to work similarly \cite{hajduk2021,kuzmin2020a,rueda-ramirez-arxiv}.
Our scheme is provably well-balanced \wrt \textit{lake at rest} equilibria, guarantees semi-discrete entropy stability and preserves positivity of water heights under a mild time step restriction.
Moreover, preservation of local bounds is enforced for the water height and velocity components.
Last but not least, we explore new approaches to realistic simulation of wetting and drying processes.
We are not aware of any other unstructured grid method that provides all of the above properties at once.
However, similar and alternative approaches to enforcing individual constraints can be found in the literature on numerical methods for the shallow water equations.
We begin this paper with a brief introduction to the state of the art.

Many well-balanced methods use the \textit{hydrostatic reconstruction} technique developed by Audusse \etal \cite{audusse2004}.
Originally proposed in the context of finite volume methods, it yields approximations that preserve lake at rest scenarios, ensure nonnegativity of water heights under standard CFL conditions, and satisfy a semi-discrete entropy inequality.
Hydrostatic reconstructions achieve these properties by properly balancing flux and source terms.
However, even the original low order hydrostatic reconstruction method does not satisfy fully discrete entropy inequalities \cite[Sec.~2.2]{audusse2004}.
This issue is addressed by Berthon \etal \cite{berthon2019}, who increase the amount of artificial viscosity to construct a method that is entropy stable in the fully discrete sense.
However, the final numerical example of their study indicates that their scheme still violates the fully discrete entropy inequality in the presence of nonflat bathymetry and dry states.

Noelle \etal \cite{noelle2007} transform to \textit{equilibrium variables} and use equilibrium-limited reconstructions in their high order finite volume methods.
The main focus of their work is on exact preservation of moving water equilibria in addition to lakes at rest.
They show that their method captures such states exactly if all stationary shocks are located at cell interfaces and Roe's numerical flux is employed.
Additionally, they prove a Lax--Wendroff-type theorem \cite[Thms.~3.14 and~3.17, respectively]{noelle2007}.

Another type of well-balanced finite volume discretizations of the SWE is the family of \textit{central-upwind schemes}.
The ones presented by Kurganov and Petrova \cite{kurganov2007a} are well balanced for the lake at rest and positivity preserving for the water height.
As in \cite{audusse2004}, these properties are achieved by performing compatible reconstructions for the conserved unknowns and the bathymetry.
To ensure positivity preservation for the water heights, the algorithm is enhanced with a generalized minmod limiter.
Additionally, a modification for numerical treatment of wetting and drying is proposed in \cite{kurganov2007a}.
In \cref{sec:swe-num}, we test this approach and some new alternatives in the context of our flux-limited finite element schemes.
In contrast to the central-upwind methods presented in \cite{kurganov2007a}, the applicability of AFC tools is not restricted to Cartesian grids.

Fjordholm \etal \cite{fjordholm2011} present well-balanced and entropy-conservative/-stable finite volume schemes for the SWE with topography.
Contrary to \cite{audusse2004}, their approach does not rely on reconstructions of the bathymetry.
Instead, a transformation to equilibrium variables is used to ensure well-balancedness for moving water equilibria.
Moreover, Fjordholm \etal generalize Tadmor's entropy stability condition to the case of nonflat topography and use it to design numerical fluxes.
It is admitted in \cite{fjordholm2011} that oscillations around discontinuities may produce negative water heights.
This shortcoming could be cured by employing a positivity-preserving limiter.

In \cite{ricchiuto2009}, Ricchiuto and Bollermann design residual distribution schemes for the shallow water equations.
As in our case, linear continuous finite elements are used in the baseline discretization.
The method preserves lake at rest scenarios and guarantees nonnegativity of the water height under CFL-like constraints.
In the SWE context, the use of residual distribution schemes introduces some complications, as mentioned in the conclusions of \cite{ricchiuto2009}.

Wintermeyer \etal \cite{wintermeyer2017} discretize the SWE using high order discontinuous Galerkin spectral element methods.
A proof of entropy conservation/stability is provided for suitable choices of numerical fluxes.
Well-balancedness \wrt lake at rest scenarios is achieved despite difficulties caused by the presence of metric terms in the case of curvilinear elements.
Although the entropy-stable DG scheme developed in \cite{wintermeyer2017} is not bound preserving, it seems to be well suited for algebraic flux correction.

Azerad \etal \cite{azerad2017} present a property-preserving finite element method that is well balanced \wrt the lake at rest.
Their scheme incorporates hydrostatic reconstructions into a nodal continuous Galerkin formulation.
Detailed analysis is performed in \cite{azerad2017} for a low order version, which is a generalization of the algebraic Lax--Friedrichs method to the case of a nonflat bottom.
Second order of accuracy is recovered by adjusting the numerical viscosity.
Extensions \cite{guermond2018a} of the schemes developed in \cite{azerad2017} are based on flux-corrected transport (FCT)-type limiting and incorporate a regularized friction term into the model.
The AFC methodology that we propose in the present paper differs from the one developed in \cite{azerad2017} in the formulation of the low order method and in the limiting strategy.
Our algorithm provides all desired properties without using hydrostatic reconstructions.

This paper is based on \cite[Ch.~4]{hajduk2022} of the first author's Ph.D. thesis and is organized as follows.
Following this introduction, we summarize the general-purpose MCL methodology for conservation laws, aspects regarding the shallow water equations and our primary objectives.
Next, we carefully extend the MCL framework to the nonconservative SWE system, whilst enforcing positivity of water heights and entropy-stability in the semi-discrete sense.
In addition, we discuss how to achieve exact preservation of lake at rest configurations.
Subsequently, we discuss existing and new approaches for simulating wetting and drying scenarios.
Numerical results are presented and conclusions are drawn in the last sections.

\section{Preliminaries}

\subsection{Monolithic convex limiting for conservation laws}\label{sec:conslaw}

Algebraic flux correction is a general framework for enforcing inequality constraints in numerical methods for conservation laws \cite{kuzmin2012}.
In this section, we review recent advances in the field of AFC for hyperbolic problems of the form
\begin{equation}\label{conslaw}
\frac{\partial u}{\partial t}+\nabla\cdot\mathbf{f}(u)=0\qquad
\mbox{in}\ \Omega\times\mathbb{R}_+.
\end{equation}
Using a continuous ($\mathbb{P}_1$ or $\mathbb{Q}_1$) finite element space \smash{$\mathrm V_h=\text{span}\{\varphi_i\}_{i=1}^N$} for discretization in space, we obtain
\begin{equation}\label{galerkin}
\sum_{j\in\mathcal N_i}m_{ij}\frac{\mathrm{d}u_j}{\mathrm{d} t}+
\sum_{j\in\mathcal N_i\backslash\{i\}}(\mathbf f_j
-\mathbf f_i)\cdot\mathbf{c}_{ij}=0,
\end{equation}
where
\begin{align*}
m_{ij}=\int_\Omega \varphi_i\,\varphi_j\dx, \qquad\text{and}\qquad
\mathbf{c}_{ij}=\int_\Omega \varphi_i\,\nabla\varphi_j\dx
\end{align*}
are entries of the consistent mass matrix, and the discrete gradient/divergence operator, respectively.
We denote by $\mathcal N_i$ i the stencil of node $\vec x_i$ and use the shorthand notation $\mathbf f_i = \mathbf f(u_i)$ for nodal fluxes.
Note that many other space discretizations can be written in this generic form.

Since the Galerkin space discretization \eqref{galerkin} may violate maximum principles and entropy conditions, we replace it
by the modified semi-discrete scheme \cite{kuzmin2020,kuzmin2012}
\begin{equation}\label{afc}
m_{i}\frac{\mathrm{d}u_i}{\mathrm{d} t}=
\sum_{j\in\mathcal N_i\backslash\{i\}}[d_{ij}(u_j-u_i)-
(\mathbf f_j-\mathbf f_i)\cdot\mathbf{c}_{ij}+f_{ij}^*],
\end{equation}
where $m_i=\sum_{j\in\mathcal N_i}m_{ij}>0$ and $f_{ij}^*$ is a suitably constrained approximation to
\begin{align}\label{eq:fij}
f_{ij}=m_{ij}(\dot u_i-\dot u_j)+d_{ij}(u_i-u_j).
\end{align}
The artificial viscosity coefficients \smash{$d_{ij} = \max\{\lambda_{ij}^{\max}|\mathbf{c}_{ij}|,\,\lambda_{ji}^{\max}|\mathbf{c}_{ji}|\}$} are defined using the maximum wave speed $\lambda_{ij}^{\max}$ of a one-dimensional Riemann problem with the initial states $u_i$ and $u_j$ \cite{guermond2016}.
The nodal time derivatives \smash{$\dot u_j=\frac{\mathrm{d}u_j}{\mathrm{d} t}$} are defined by the solution of system \eqref{galerkin}. 

By construction, our AFC scheme \eqref{afc} reduces to \eqref{galerkin} for $f_{ij}^*=f_{ij}$.
The choice $f_{ij}^*=0$ corresponds to an algebraic Lax--Friedrichs (ALF) method.
Guermond and Popov \cite{guermond2016} proved that
this ``first-order'' approximation to \eqref{galerkin} is
invariant domain preserving (IDP) in the sense that the nodal values $u_i$ stay in a convex set $\mathcal A$ if all initial values belong to this set.
Moreover, the validity of a fully discrete entropy inequality can be shown for any convex entropy pair.
The crux of the proofs presented in \cite{guermond2016} is representation of the ALF scheme in terms of the \textit{bar states}
\begin{equation}
\label{barstate}
\bar u_{ij}=\frac{u_j+u_i}{2}
-\frac{(\mathbf f_j-\mathbf f_i) \cdot\mathbf{c}_{ij}}{2d_{ij}}
\end{equation}
such that $\bar u_{ij}\in\mathcal A$ if $\mathcal A$ is an invariant set of \eqref{conslaw} and $u_i,u_j\in\mathcal A$.
A ``second-order'' IDP scheme for general hyperbolic problems was designed in \cite{guermond2018} using \textit{convex limiting}
based on a localized flux-corrected transport
algorithm. 

The \textit{monolithic} convex limiting strategy that we favor in the present paper differs from FCT-like predictor-corrector approaches in that the fluxes $f_{ij}^*$ are used to correct the right-hand side of \eqref{afc} rather than a low-order solution obtained with the ALF method.
The semi-discrete IDP limiting criterion is given by \cite{kuzmin2020}
\begin{equation}
\bar u_{ij}\in\mathcal A_i\quad \Rightarrow \quad
\bar u_{ij}^*:=\bar u_{ij}+\frac{f_{ij}^*}{2d_{ij}}
\in\mathcal A_i,
\end{equation}
where $\mathcal A_i\subseteq\mathcal A$ is a convex set of
admissible values and $\bar u_{ij}^*$ is a flux-corrected
counterpart of the ALF bar state $\bar u_{ij}$ defined by
\eqref{barstate}.
Additionally, we constrain the fluxes $f_{ij}^*$ to satisfy an entropy stability condition which implies the validity of a semi-discrete entropy inequality; see \cite{kuzmin2020e,kuzmin2020f} for details.
Integration in time is performed using a strong stability
preserving (SSP) Runge--Kutta method \cite{gottlieb2001,gottlieb2011,shu1988}.

\subsection{The shallow water equations}

In the previous section, we discussed monolithic convex limiting strategies and related concepts that were developed in \cite{kuzmin2020,kuzmin2020e} for conservation laws.
We now extend these techniques to a system of balance laws.
In particular, we consider the shallow water equations with a nonconservative topography term.
This important nonlinear system of partial differential equations reads
\begin{align}\label{eq-swe-topo}
\frac{\partial}{\partial t} \begin{bmatrix}
h \\ h\vv
\end{bmatrix} + \div \begin{bmatrix}
h \vv\transp \\ h \vv \otimes \vv + \frac g 2 h^2 \id
\end{bmatrix}
+ \begin{bmatrix}
0\\
gh\nabla b
\end{bmatrix} = 0 
\end{align}
Here $h$ is the total water height, $\vv\in \R^d$, $d\in \{1,2\}$ is the velocity vector, $g$ is the gravitational constant, $\id \in \R^{d\times d}$ is the identity matrix, and $b$ is the bottom topography or \textit{bathymetry}.
Additional theoretical and numerical challenges arise when it comes to solving \eqref{eq-swe-topo} instead of the system of conservation laws corresponding to the case $b\equiv{}$const.
We refer to \cite[Ch.~3]{bouchut2004} and references cited therein for a review of the theory of balance laws and some aspects of the shallow water equations with topography.

\subsection{Objectives}\label{sec:obj}

The goal of this paper is to generalize the bound-preserving and entropy-stable MCL schemes to the inhomogeneous SWE \eqref{eq-swe-topo}.
One key requirement that we deemed essential in the development of our algorithms is that they represent generalizations of the corresponding schemes from \cref{sec:conslaw} for the flat bottom case.
Another desirable property of numerical methods for balance laws is their \textit{well-balancedness}.
For some hyperbolic systems, there exist certain steady states\ie solutions $u(\x)$ that are independent of $t$ because the flux and source terms are in equilibrium.
A numerical method for solving such a system is called well balanced if it captures the simplest steady states exactly in the discrete setting.
In the derivation of the SWE one uses the chain rule to rewrite the term $gh\nabla (h+b)$ as \smash{$\frac g2\nabla h^2 + gh\nabla b$}.
This decomposition suggests that system \eqref{eq-swe-topo} admits the so-called \textit{lake at rest} steady state solution
\begin{align}\label{eq-lakeatrest}
\vv = \vec 0, \qquad h\nabla(h+b) = 0.
\end{align}
This configuration corresponds to a still body of water that is unperturbed by external forces, such as in- and outflows through domain boundaries.
Note that the second identity in \eqref{eq-lakeatrest} does not imply that the free surface elevation $H=h+b$ has to be a global constant, as is the case for a classical lake at rest.
In fact, \eqref{eq-lakeatrest} allows variations in $H$, as long as the water height $h$ is zero at the same physical location.
This case corresponds to a so-called \textit{dry state} that occurs whenever the topography $b$ exceeds the water depth $h$.
For an island that rises from a body of water, every point on the surface of the island represents a dry state.

Besides lake at rest configurations, other types of equilibria exists for the SWE\@.
In the absence of friction and/or Coriolis forces, \eqref{eq-swe-topo} admits so-called \textit{moving water equilibria} steady states.
In 1D, such configurations occur if the discharge $hv$ as well as the expression \smash{$\frac12 v^2 + g(h+b)$} remain constant \cite[Ch.~3]{bouchut2004}, \cite{kurganov2007a}.
While lake at rest scenarios can be preserved with simple numerical treatments (see, \eg \cite{audusse2004,berthon2019,kurganov2007a,
fjordholm2011,azerad2017}), moving water equilibria require advanced well-balancing techniques (see for instance \cite{noelle2007}).
The incorporation of such approaches into our flux correction schemes is a topic of its own and will not be attempted in this work.
Instead, we focus on well-balancing \wrt lake at rest configurations.
Nevertheless, some numerical examples of moving water equilibria are solved numerically in this paper.

Another important aspect of numerical methods for the SWE and related models is the need to deal with \textit{wetting and drying} scenarios (see, \eg \cite{ricchiuto2009,barros2015,vater2015}) in which simulations may crash if no special measures are implemented.
In many examples of practical interest, there exist islands rising from the body of water but the interface between these islands and the water surface is moving.
The dry land masses are then modeled by allowing the bottom topography to exceed the values of the free surface elevation at the same location.
Even in the case that the resolution is sufficient to capture the interface, it can be quite difficult to accurately resolve the moving shoreline with numerical methods.
In \cref{sec:wad}, we present two new ways of dealing with this issue and compare our results with some existing approaches.

\section{Limiting for the shallow water equations with topography}\label{sec:swe-afc}

Let us now extend the schemes discussed in \cref{sec:conslaw} to the inhomogeneous hyperbolic system \eqref{eq-swe-topo} step by step.
We first choose a target discretization suitable for flux correction procedures.
Next, we derive a low order method for which all desirable properties (conservation, numerical and physical admissibility, entropy stability) are guaranteed.
Then we recover the target scheme by including raw antidiffusive fluxes into the algorithm.
Finally, these fluxes are limited in a way that ensures preservation of both local and global bounds, as well as semi-discrete entropy stability.

The most important considerations regarding the low order method and flux-corrected schemes were already discussed in the context of general systems of conservation laws in \cref{sec:conslaw}.
For brevity, we focus solely on aspects that need to be modified and refer to \cite{kuzmin2020} and \cite[Sec.~3.3]{hajduk2022} for the rest.
In particular, the treatment of boundary terms does not present any additional difficulties because these terms are the same for discretizations of \eqref{eq-swe-topo} with constant and spatially variable bathymetry $b$.
Therefore, we omit all boundary terms in the following presentation but remark that they generally need to be incorporated into the algorithm.
This task is achieved by choosing appropriate external Riemann data based on the physics of the problem and incorporating boundary contributions weakly via a numerical flux.

Conceptually, the only difference compared to the case of a flat topography is the presence of the nonconservative term $gh\nabla b$ in the momentum equation.
The consistent Galerkin discretization of this term produces the nodal contribution
\begin{align}\label{eq-source-cons}
g\sum_{k\in \mathcal N_i} h_k \int_\Omega \varphi_i \,\varphi_k \nabla b \dx
\end{align}
to the $i$th component of the momentum equation.
In this formula, $h_k$ denotes the value of the discretized water height at node $\x_k$ and $\varphi_k$ are the corresponding piecewise (multi-)linear continuous Lagrange basis functions of the finite element space $\mathrm V_h$.
It is common \cite{kurganov2007a} to approximate $b$ by its piecewise (multi)-linear continuous interpolant $b_h\in \mathrm V_h$, defined by
\begin{align*}
b_h(\x) = \sum_{j=1}^N b_j \varphi_j(\x), \qquad b_j \coloneqq b(\x_j).
\end{align*}
If the bathymetry $b$ is discontinuous in node $\x_j$, one may set $b_j$ equal to any of the one-sided limits in cells containing $\x_j$ or to an average of these limits.
Alternatively, a projection can be used to obtain $b_h$ from $b$.

In \cref{sec:conslaw}, we approximated the inviscid flux $\f(u_h)$ using a group finite element formulation \cite{fletcher1983,barrenechea2017b} to derive a quadrature rule for the corresponding integral.
Discretization \eqref{eq-source-cons} of the source term must be approximated similarly for our method to be well balanced.
With this goal in mind, we replace \eqref{eq-source-cons} by the quadrature-based version
\begin{align}\label{eq-source-quad}
g \afcsum \frac{h_i+h_j}2 (b_j - b_i)\,\cij,
\end{align}
which is similar to what is done, for instance, in \cite[Eq.~(3.8)]{audusse2004}, \cite[Eq.~(2.6)]{kurganov2007a}, and \cite[Eq.~(2.7)]{fjordholm2011}.
Note that if $h_k\equiv{}$const for $k\in\mathcal N_i$, then \eqref{eq-source-quad} equals \eqref{eq-source-cons} with $b$ replaced by $b_h$.
In this sense, \eqref{eq-source-quad} is similar to the quadrature rule based on the group finite element approximation of $\f(u_h)$.
The approximation \eqref{eq-source-quad} to \eqref{eq-source-cons} is second order accurate if $h_k$ is constant for $k\in\mathcal N_i$ and first order accurate otherwise.
\begin{remark}
In principle, it is possible to compensate the quadrature error due to the source term approximation \eqref{eq-source-quad} in the process of flux correction.
If this is desired, one needs to decompose the difference between \eqref{eq-source-cons} and \eqref{eq-source-quad} into edge contributions, add the corresponding correction terms \smash{$\vec{f}_{ij}^b$} to the raw antidiffusive fluxes \smash{$\vec{f}_{ij}^{h\vv}=-\vec{f}_{ji}^{h\vv}$} of the momentum equation, and modify the limiting formula because \smash{$\vec{f}_{ij}^b\ne -\vec{f}_{ji}^b$} in general.
At an early stage of developing our method, we performed a preliminary study that showed the feasibility of this approach.
In the final version, we disregard first order quadrature errors caused by using \eqref{eq-source-quad} instead of \eqref{eq-source-cons} because in real-life applications such errors are likely to be negligible compared to measurement errors in the bathymetry data.
\end{remark}
Inserting \eqref{eq-source-quad} into the semi-discrete momentum equation and approximating other terms as in \cref{sec:conslaw}, we obtain the quadrature-based target scheme
\begin{subequations}\label{eq-swe-target}
\begin{align}\label{eq-h-target}
\sum_{j=1}^N \mij \frac{\d h_j}{\d t} ={}&
-\afcsum ((h\vv)_j - (h\vv)_i)\transp\cij,\\
\sum_{j=1}^N \mij \frac{\d (h\vv)_j}{\d t} ={}&
-\afcsum \big[ \fhv_j - \fhv_i + g \frac{h_i+h_j}{2} (b_j - b_i)\,\id\big]\, \cij, \label{eq-mom-target}
\end{align}
\end{subequations}
where
\begin{align*}
\fhv_i =
\frac 1 {h_i} (h\vv)_i\otimes (h\vv)_i
+ \frac g 2 h_i^2\,\id, \qquad i\in \{1,\ldots,N\}.
\end{align*}
It is worth checking at this stage whether \eqref{eq-swe-target} is well balanced for lake at rest configurations \eqref{eq-lakeatrest} and clarify the meaning of well-balancedness in this context.
If $(h\vv)_i = \vec 0$ for all $i\in \{1,\ldots,N\}$, then
\eqref{eq-mom-target} reduces to
\begin{align*}
\sum_{j=1}^N \mij \frac{\d (h\vv)_j}{\d t} ={}& -\afcsum
\frac g2 [h_j^2-h_i^2 + (h_i+h_j)(b_j-b_i)]\,\cij \\
={}& -\afcsum \frac g2 (h_j+h_i)[h_j - h_i + b_j -b_i]\,\cij.
\end{align*}
Assuming for now that $h_i\ge 0$ for all $i\in \{1,\ldots,N\}$ (a condition that we enforce later on), we see that the right hand side of this expression is zero if and only if
\begin{align}\label{eq-wb-cond}
h_i=h_j=0\qquad \text{or} \qquad H_i = H_j \qquad \forall \afcijmath,
\end{align}
where $H_i\coloneqq h_i + b_i$, $i\in\{1,\ldots,N\}$ are the coefficients of the discrete free surface elevation $H_h \in \mathrm V_h$.
If \eqref{eq-wb-cond} holds, the discharge is unperturbed, and thus the right hand side of the continuity equation \eqref{eq-h-target} is zero.

\subsection{Low order method}
As in the case of a system of conservation laws, we need to modify \eqref{eq-swe-target} to obtain a property-preserving semi-discretization.
To this end, we perform row sum mass lumping and include Rusanov (local Lax--Friedrichs) artificial dissipation.
However, the presence of the nonconservative term makes matters more involved.
Special care needs to be taken, for instance, to ensure semi-discrete entropy stability, positivity preservation for water heights, and well-balancedness.
To construct a low order method that meets all of our requirements, let us begin with the straightforward generalization
\begin{subequations}\label{eq-naive}
\begin{align}\label{eq-naive-h}
m_i \frac{\d h_i}{\d t} ={}&
\afcsum \big[ \dij (h_j - h_i) - ((h\vv)_j - (h\vv)_i)\transp\cij\big],\\
m_i \frac{\d (h\vv)_j}{\d t} ={}&
\afcsum \big[ \dij ((h\vv)_j - (h\vv)_i) - (\f^{h\vv}_j - \f^{h\vv}_i)\,\cij \notag\\
&\qquad\quad{}- \frac g2 (h_i+h_j) (b_j - b_i)\,\cij\big]\label{eq-naive-mom}
\end{align}
\end{subequations}
of the algebraic Lax--Friedrichs method applied to the SWE with flat bathymetry.
As before, the artificial viscosity coefficients $\dij$ are defined by \cite{leveque2002}
\begin{align}\label{eq:dij}
\dij = \max\{\lambda_{ij}^{\max}|\mathbf{c}_{ij}|,
\lambda_{ji}^{\max}|\mathbf{c}_{ji}|\}, \qquad
\lambda_{ij}^{\max} = \max \lc\left|\vv_i \cdot \frac{\mathbf{c}_{ij}}{|\mathbf{c}_{ij}|}\right| + \sqrt{gh_i},~\left|\vv_j \cdot \frac{\mathbf{c}_{ij}}{|\mathbf{c}_{ij}|}\right| + \sqrt{gh_j}\rc.
\end{align}
We will modify \eqref{eq-naive} step by step until we are able to prove the desired properties.

A first observation regarding \eqref{eq-naive} is that this scheme does not preserve the lake at rest \eqref{eq-lakeatrest} if the given velocity is zero, the free surface elevation is constant ($h_i+b_i=H_i=H_j=h_j+b_j$ for all pairs of nodes) but $h_i \neq h_j$ for some $j\ne i$.
In this scenario, the flux $\dij (h_j - h_i)$ of the semi-discrete continuity equation \eqref{eq-naive-h} will disturb the equilibrium and produce nonphysical waves.
This issue can be resolved by replacing \eqref{eq-naive-h} with
\begin{align}\label{eq-swe-conti-disc}
m_i \frac{\d h_i}{\d t} ={}&
\afcsum \lb \dij (H_j - H_i) - ((h\vv)_j - (h\vv)_i)\transp\cij\rb.
\end{align}
This discretization preserves the lake at rest in the case $H\equiv{}$const.
However, the theory that was used to prove positivity preservation the low order method does not carry over to systems of balance laws.
For the SWE with flat bottom, nonnegativity of water heights follows from the fact that the low order bar states are averaged exact solutions of the Riemann problem \cite{guermond2016}.
If source terms are included, the so-defined intermediate states may fail to stay in the admissible set \smash{$\mathcal A^{\max}=\{(h,h\vv\transp)\transp \in \R^{d+1}: h\ge 0\}$} of the homogeneous SWE\@.
Thus, we need to enforce the nonnegativity constraint for $h_i$ by modifying the discretization of the continuity equation.
To this end, we notice that $H_j-H_i = h_j-h_i + (b_j-b_i)$ and introduce a bathymetry limiter \smash{$\lij^b \in [0,1]$} that transforms \eqref{eq-swe-conti-disc} into
\begin{align}
m_i \frac{\d h_i}{\d t} ={}&
\afcsum \lb \dij (h_j - h_i + \lij^b (b_j -b_i)) - ((h\vv)_j - (h\vv)_i)\transp\cij\rb \label{eq-swe-conti-lim1} \\ ={}&
\afcsum 2\dij \big[ \hij - h_i + \frac{\lij^b}{2}(b_j - b_i)\big] = \afcsum 2\dij ( \hij^b - h_i ).\notag
\end{align}
Here $\hij$ is the first component of the usual low order bar state \eqref{barstate} and
\begin{align}\label{eq-hij}
\hij^b \coloneqq \hij + \frac{\lij^b}{2}(b_j - b_i).
\end{align}
The correction factor $\lij^b$ is used to ensure that \smash{$\hij^b \ge 0$}.
This condition holds for $\lij^b=0$ by definition of $\hij^b$.
However, the largest admissible value of $\lij^b\in[0,1]$ should be employed for consistency reasons.
To maintain the conservation property of the semi-discrete continuity equation, we impose the usual symmetry condition $ \lij^b=\alpha_{ji}^b$.
Note that for $b_j-b_i\ge 0$, the use of $\lij^b=1$ in \eqref{eq-hij} cannot produce negative \smash{$\hij^b$} provided that \smash{$\hij\ge 0$}.
In this case, however, the limiter may need to act to enforce the condition \smash{$\bar h_{ji}^b \ge 0$}.
These considerations lead to the definition
\begin{align}\label{eq-alphab}
\lij^b = \begin{cases}
\min \lc 1, \frac{2\bar h_{ji}}{b_j-b_i} \rc & \text{if } b_i -b_j < 0, \\
1 & \text{if } b_i - b_j = 0,\\
\min \lc 1, \frac{2\hij}{b_i-b_j} \rc & \text{if } b_i -b_j > 0.
\end{cases}
\end{align}
This approach to enforcing the nonnegativity of water heights in the low order method is equivalent to the correction procedure proposed by Audusse \etal \cite[Sec.~2.2]{audusse2015}.
The authors of \cite{audusse2015} also impose the conservation and positivity requirements, which yields a 1D finite volume version of our water height limiter based on \eqref{eq-alphab}.

It is worth checking how the limiter \eqref{eq-alphab} behaves for lake at rest configurations.
\begin{lemma}[Well-balancedness of the positivity-preserving limiter]\label{lem:lar}
Let $u_h\in \mathrm V_h^{d+1}$ be a lake at rest solution\ie assume that \eqref{eq-wb-cond} holds in addition to $h_i \ge 0$ and $(h\vv)_h=\vec 0$.
Then for any $i\in \{1,\ldots,N\}$ and $j\in\mathcal N_i\setminus\{i\}$ we have
\begin{enumerate}
\item[i)]
$\lij^b(b_j-b_i)=0$ if $h_i=h_j=0$ and
\item[ii)]
$\lij^b=1$ if $H_i=H_j$.
\end{enumerate}
In either case, the application of the bathymetry limiter \eqref{eq-alphab} does not perturb the lake at rest state because the right hand side of \eqref{eq-swe-conti-lim1} is zero for the given data.
\end{lemma}
\noindent
\textbf{Proof:} See \cite[Pf. of Lem. 4.2]{hajduk2022}. \hfill $\square$
\medskip
%
%

At this stage, one may be tempted to use the spatial semi-discretization consisting of \eqref{eq-swe-conti-lim1} and \eqref{eq-naive-mom} as a low order method for algebraic flux correction.
While this version is already usable, it does not yet ensure semi-discrete entropy stability for general bathymetry.
For this reason, we modify the momentum equation \eqref{eq-naive-mom} as follows
\begin{align}\notag
m_i \frac{\d (h\vv)_i}{\d t} ={}& \afcsum \big[ \dij \big( (h\vv)_j - (h\vv)_i + \frac{\vv_i+\vv_j}{2} \lij^b (b_j - b_i) \big) \\
& \quad\quad\quad{}- (\f_j^{h\vv} - \f_i^{h\vv})\,\cij - g \frac{h_i+h_j}{2} \lij^b (b_j - b_i) \, \cij \big].
\label{eq-swe-mom}
\end{align}
The term \smash{$\frac12(\vv_i+\vv_j) \lij^b\,d_{ij}(b_j - b_i)$} is included for entropy stabilization purposes.
For consistency reasons, we apply the correction factor $\lij^b$ to all bathymetry fluxes.

\begin{remark}
Even though the bathymetry plays the role of a parameter in the SWE model, the correction factor $\lij^b$ adjusts the source term contribution to the momentum equation.
The consistency error introduced in this way is acceptable because $\lij^b\ne 1$ is used only for (neighbors of) dry states.
A similar concept is employed in the popular hydrostatic reconstruction approach \cite{audusse2004,berthon2019,azerad2017} in which topography values are locally adjusted to guarantee nonnegativity of water heights and well-balancedness.
\end{remark}
Let us now discuss how to generalize Tadmor's entropy stability condition \cite{tadmor1987,tadmor2003} to our setting and verify it for the low order method.
An entropy pair of the SWE with nonflat topography is given by \cite{wintermeyer2017}
\begin{align}\label{eq-entropypair}
\eta(u,b) = \frac 1 2\l gh^2 + h|\vv|^2\r + ghb, \qquad
\vec q(u,b) = \big( g(h+b) + \frac 1 2|\vv|^2 \big) h\vv.
\end{align}
As of now, our proofs of entropy stability are, in fact, limited to this entropy pair.
The entropy variable and potential corresponding to \eqref{eq-entropypair} read
\begin{align}\label{eq-entr-swe}
v(u,b) = \begin{bmatrix}
g(h+b) - \frac 1 2 |\vv|^2\\ \vv
\end{bmatrix} = v(u,0) + \begin{bmatrix}
gb\\\vec 0
\end{bmatrix}, \qquad \vec\psi(u,b) = \vec \psi(u) = \frac g2h^2 \vv.
\end{align}
A generalized version of Tadmor's entropy stability condition was derived by Fjordholm \etal \cite[Sec.~2.1]{fjordholm2011} in the context of finite volume methods for structured grids.
Adapting this generalization to our continuous FEM setting, we arrive at
\begin{align}\notag
\frac{\dij}2 P_{ij} \le{}& (\vec \psi_j - \vec \psi_i)\cdot\cij
+ \frac{(v(u_i,0) - v(u_j,0))\transp}{2} (\f_j + \f_i)\,\cij \\&
+ g \bigg[ \frac{h_i+h_j}{2} \frac{\vv_i+\vv_j}{2} - \frac{(h\vv)_i+(h\vv)_j}{2} \bigg]\transp \cij \lij^b(b_j - b_i) \eqqcolon Q_{ij},
\label{eq-mod-tadmor}
\end{align}
where
\begin{align*}
P_{ij} \coloneqq{}&
\begin{bmatrix}
g[h_i - h_j + \lij^b(b_i-b_j)] -\frac{|\vv_i|^2-|\vv_j|^2}2 \\ \vec \vv_i - \vv_j
\end{bmatrix}\transp
\begin{bmatrix}
h_j - h_i + \lij^b(b_j - b_i) \\
(h\vv)_j - (h\vv)_i + \frac{\vv_i + \vv_j}{2} \lij^b (b_j -b_i)
\end{bmatrix}.
\end{align*}
Inequality \eqref{eq-mod-tadmor} imposes the upper bound \smash{$(\vec \psi_j - \vec \psi_i)\cdot\cij$} on the rates of entropy production/dissipation due to low order fluxes and source terms.
It turns out that the parameters $d_{ij}$ of our low order method \eqref{eq-swe-conti-lim1},\,\eqref{eq-swe-mom} can be chosen sufficiently large (overestimating the maximum speed if necessary) to ensure the validity of \eqref{eq-mod-tadmor}.

\begin{lemma}[Entropy stability of the low order method]\label{lem:dij-suf}
There exist coefficients $d_{ij}\ge 0$ such that condition \eqref{eq-mod-tadmor} holds for the numerical fluxes of the semi-discrete low order method defined by \eqref{eq-swe-conti-lim1} and \eqref{eq-swe-mom}.
\end{lemma}
\noindent
\textbf{Proof:} See \cite[Pf. of Lem. 4.4]{hajduk2022}. \hfill $\square$
\medskip

\noindent
To ensure entropy stability of the low order method \eqref{eq-swe-conti-lim1},\,\eqref{eq-swe-mom} in practice, we verify whether the coefficients $\dij$ defined by \eqref{eq:dij} are large enough to satisfy \smash{$\frac{\dij}2 P_{ij} \le Q_{ij}$}.
If this is not the case, we set $\dij=d_{ji}= 2\min\{0,Q_{ij},Q_{ji}\}/P_{ij}$.
In practical applications, such adjustments seem to be necessary only in the vicinity of dry states.
For such configurations, our approach of increasing the artificial viscosity does not significantly reduce the time step if an appropriate wetting and drying treatment is adopted.
\begin{remark}\label{rem:tadmor}
As an alternative to adjusting the diffusion coefficients $d_{ij}$, condition \eqref{eq-mod-tadmor} can be satisfied by further reducing the value of $\lij^b$.
Indeed, \eqref{eq-mod-tadmor} reduces to Tadmor's usual entropy stability condition \cite{tadmor1987,tadmor2003} for $\lij^b = 0$.
A formula to compute such $\lij^b$ can be derived similarly to the IDP pressure fix for the Euler equations discussed in \cite[Sec.~5.1]{kuzmin2020}.
\end{remark}
Let us now generalize the setting of \cref{sec:conslaw} to derive local semi-discrete entropy inequalities.
First, we rewrite the low order method \eqref{eq-swe-conti-lim1},\,\eqref{eq-swe-mom} as follows
\begin{align*}
m_i \frac{\d u_i}{\d t} ={}& \afcsum [\gij -(\f_j + \f_i)\,\cij + s_{ij}] + 2\f_i \afcsum \cij,
\end{align*}
where 
\begin{align*}
\gij = \dij\begin{bmatrix}
h_j - h_i + \lij^b(b_j - b_i) \\
(h\vv)_j - (h\vv)_i + \frac{\vv_i + \vv_j}{2} \lij^b (b_j -b_i)
\end{bmatrix}, \qquad
s_{ij} = \begin{bmatrix}
0\\ - g \frac{h_i+h_j}{2} \lij^b (b_j - b_i)\,\cij
\end{bmatrix}.
\end{align*}
Thus, $\gij = - g_{ji}$ for all \afcijtext~and $s_{ij}=s_{ji}$ if $\cij=-\mathbf c_{ji}$.
\begin{theorem}[Local semi-discrete entropy inequality]\label{thm:swe}
Consider the low order method \eqref{eq-swe-conti-lim1},\,\eqref{eq-swe-mom} satisfying \smash{$\frac{\dij}2 P_{ij} \le \min\{ Q_{ij}, Q_{ji}\}$} for all \afcijtext.
Define
\begin{align*}
G_{ij} \coloneqq{}& \frac{(v_i+v_j)\transp}{2} \gij + \frac{(v_i - v_j)\transp}{2} [(\f_i - \f_j)\,\cij + s_{ij}],\\
W_{ij} \coloneqq{}& \frac g 2 (1-\lij^b)(b_i-b_j)[\dij(h_j-h_i + \lij^b(b_j-b_i)) - ((h\vv)_i + (h\vv)_j)\transp\cij].
\end{align*}
Then for all $i \in\{1,\ldots,N\}$ the semi-discrete entropy inequalities
\begin{align}\label{eq-swe-entropy-ineq}
m_i \frac{\d \eta_i}{\d t} \le \afcsum [ G_{ij} + W_{ij} - (\vec q_j - \vec q_i)\cdot\cij ]
\end{align}
hold \wrt the entropy pair $(\eta,\vec q)$ defined by \eqref{eq-entropypair}.
\end{theorem}
\noindent
\textbf{Proof:} See \cite[Pf. of Thm 4.6]{hajduk2022}. \hfill $\square$
\medskip

\noindent
The consistency errors $W_{ij}$ can be attributed to the occurrence of dry or almost dry areas, which require the use of $\lij^b<1$.
For such states, even the validity of the continuous entropy inequality is questionable because the momentum equation of the SWE model does not describe the underlying physics correctly.
In particular, the absence of friction terms becomes an issue.
This argument justifies the presence of $W_{ij}$ in \eqref{eq-swe-entropy-ineq}.
\begin{corollary}[Global semi-discrete entropy inequality]\label{cor:entropy-swe}
Let the assumptions of \mbox{Theorem~\ref{thm:swe}} be fulfilled and assume that the spatial semi-discretization \eqref{eq-swe-conti-lim1}, \eqref{eq-swe-mom} satisfies
\begin{align*}
\mathbf c_{ij}=-\mathbf c_{ji} \qquad \text{or} \qquad u_i(t) = u_j(t) \qquad i \in \{1,\ldots,N\},~j\in \mathcal N_i\setminus \{i\},
\end{align*}
and, furthermore, the external Riemann data coincides with the semi-discrete approximation $u_h(t)$ in all boundary nodes.
(This assumption is true, for instance, at supercritical outlets.)
If, in addition, $\lij^b=1$ for all \afcijtext, then the following semi-discrete entropy inequality holds
\begin{align*}
\frac{\d}{\d t} \int_\Omega \big( \sum_{i=1}^N \eta_i\,\phi_i\big) \dx + \int_{\partial\Omega} \big(\sum_{i=1}^N \vec q_i\, \phi_i\big) \cdot \n \ds\le 0.
\end{align*}
\end{corollary}
\noindent
\textbf{Proof:} See \cite[Pf. of Cor. 4.7]{hajduk2022}. \hfill $\square$
\medskip

\noindent
We conclude the discussion of the low order method by formulating the bar state form of the momentum equation \eqref{eq-swe-mom}, which reads
\begin{align*}
m_i \frac{\d (h\vv)_i}{\d t} ={}& \afcsum 2\dij \big( \hvij^b - (h\vv)_i \big).
\end{align*}
Similarly to the bar state \eqref{eq-hij} of the water height, the bar state of the discharge
\begin{align*}
\hvij^b ={}& \hvij + \frac{\vv_i + \vv_j}{4}\lij^b(b_j-b_i) - \frac{g\frac{h_i+h_j}2 \lij^b(b_j-b_i)\,\cij}{2\dij}\\ ={}&
\frac{(h\vv)_i+(h\vv)_j + \frac{\vv_i+\vv_j}{2}\lij^b(b_j-b_i)}{2} - \frac{\l\f_j^{h\vv}-\f_i^{h\vv} + g\frac{h_i+h_j}2 \lij^b(b_j-b_i)\,\id \r \cij}{2\dij}
\end{align*}
consists of symmetric and skew-symmetric terms.
\begin{remark}
A similar concept based on intermediate states is employed in the work of Audusse \etal \cite{audusse2015}.
We already mentioned the fact that their limiter for the water height \cite[Sec.~2.2]{audusse2015} is equivalent to \eqref{eq-alphab}.
A minor difference between our low order method and the well-balanced scheme of Audusse \etal is that their finite volume method employs intermediate states based on the HLL Riemann solver \cite{harten1983a} instead of local Lax--Friedrichs-type bar states.
More importantly, our discretization of the momentum equation includes an entropy-stabilizing term, which is missing in \cite{audusse2015}.
The absence of this term might be the reason why no conclusive evidence regarding the validity of discrete entropy inequalities could be provided in \cite[Sec.~2.4]{audusse2015}.
\end{remark}

\subsection{Monolithic convex limiting}

Having derived the low order method \eqref{eq-swe-conti-lim1},\,\eqref{eq-swe-mom}, we now discuss the MCL methodology for the SWE with topography.
As in the case of conservation laws, we first need to define the raw antidiffusive fluxes $\fij=-f_{ji} \in \R^{d+1}$, \smash{$\fij\coloneqq (\fij^h, (\vec f_{ij}^{h\vv})\transp)\transp$} with which the target scheme \eqref{eq-swe-target} can be recovered from the low order method.
In contrast to \eqref{eq:fij}, we have to include additional terms due to modifications that make our low order method property preserving for nonflat topography.
A straightforward computation shows that if $\lij^b=1$ for all \afcijtext, then \eqref{eq-swe-target} can be recovered via
\begin{align*}
&\phantom{\frac{(h\vv)}{t}} m_i \frac{\d h_i}{\d t} =
\afcsum \big[ 2\dij (\hij^b - h_i) + \fij^h \big],\\
&\phantom{\frac{h}{t}} m_i \frac{\d (h\vv)_i}{\d t} = \afcsum \big[ 2\dij \big(\hvij^b - (h\vv)_i\big) + \vec f_{ij}^{h\vv}\big],
\end{align*}
where
\begin{align*}
\fij^h ={}& \mij \big(\dot h_i - \dot h_j\big) + \dij \lb h_i - h_j + \lij^b(b_i - b_j) \rb,\\
\vec f_{ij}^{h\vv} ={}& \mij \l\dot{(h\vv)}_i - \dot{(h\vv)}_j\r + \dij \big[ (h\vv)_i - (h\vv)_j + \frac{\vv_i+\vv_j}{2}\lij^b(b_i - b_j) \big].
\qquad\quad
\end{align*}
In our fully discrete scheme, we once more define the dotted quantities as low order time derivatives, which are computed from \eqref{eq-swe-conti-lim1} and \eqref{eq-swe-mom}, respectively.
For steady state problems, these quantities are set to zero because mass lumping errors need to be compensated only for time-dependent problems.

We are now in a position to present the generalized sequential limiting technique \cite{kuzmin2020,dobrev2018} with which we obtain flux-corrected counterparts \smash{$f_{ij}^{h,*}$ and $\vec f_{ij}^{h\vv,*}$ of $\fij^h$ and $\vec f_{ij}^{h\vv}$.}
In the first step of the sequential MCL algorithm, we limit the water height using
\begin{align*}
h_i^{\min} \coloneqq \min_{j\in \mathcal N_i\setminus\{i\}} \bar h_{ij}^b, \qquad
h_i^{\max} \coloneqq \max_{j\in \mathcal N_i\setminus\{i\}} \bar h_{ij}^b
\end{align*}
as local bounds of numerical admissibility conditions, which imply global positivity preservation for the water height if the bathymetry correction factor $\alpha_{ij}^b$ defined by \eqref{eq-alphab} is applied.
The limiting formula for the raw antidiffusive fluxes $f_{ij}^h$ becomes (\cf \cite[Eq.~(46)]{kuzmin2020})
\begin{align}\label{eq-mcl-H}
f_{ij}^{h,*} ={}& \begin{cases}
\spacingmin\min\lc
f_{ij}^h,\; 2d_{ij}\spacingmin\min \lc h_i^{\max} - \bar h_{ij}^b, \bar h_{ji}^b - h_j^{\min} \rc \rc & \text{if }f_{ij}^h \geq 0, \\
\spacingmax\max\lc
f_{ij}^h,\; 2d_{ij}\spacingmax\max \lc h_i^{\min} - \bar h_{ij}^b, \bar h_{ji}^b - h_j^{\max} \rc \rc & \text{if }f_{ij}^h\le 0.
\end{cases}
\end{align}
The corresponding flux-corrected bar states in the continuity equation can be written as
\begin{align*}
\hij^{b,*} = \hij^b + \frac{f_{ij}^{h,*}}{2\dij} = \hij + \frac{\lij^b(b_j-b_i)}2 + \frac{f_{ij}^{h,*}}{2\dij} = \bar h_{ij}^* + \frac{\lij^b(b_j-b_i)}2.
\end{align*}
Next, we need to limit $\vec f_{ij}^{h\vv}$ in a way that ensures the validity of numerical admissibility conditions for individual velocity (rather that discharge) components.
To construct local bounds for this step, we first define the velocity bar states as (\cf \cite[Eq.~(80)]{kuzmin2020})
\begin{align*}
\bar\vv_{ij} \coloneqq \frac{\hvij^b + \overline{(h\vv)}_{ji}^b}{\hij^b + \bar h_{ji}^b} = \frac{2\dij\big( \hvij + \overline{(h\vv)}_{ji} \big) - g\frac{h_i+h_j}2\lij^b (b_j-b_i)(\cij-\mathbf c_{ji})}{2\dij(\hij + \bar h_{ji})} = \bar\vv_{ji},
\end{align*}
which represents a generalization of the velocity bar states of the homogeneous SWE.
Note that the components $\frac12\lij^b(b_j-b_i)$ and $\frac14(\vv_i+\vv_j)\lij^b(b_j-b_i)$ of the bar states $\bar h_{ij}^b$, $\hvij^b$ and the corresponding antisymmetric components of the bar states $\bar h_{ji}^b$, $\overline{(h\vv)}_{ji}^b$ cancel out upon summation in the numerator and denominator of the second ratio.
The symmetric source terms add up in the numerator.

Let $\vv_i^{\min},\,\vv_i^{\max}\in\R^d$ be vectors containing local bounds to be imposed on individual components of the nodal velocity (alternative limiting strategies for vector fields, in particular ones that guarantee rotational invariance, are discussed in \cite{hajduk2019}).
Inequalities involving vectors should be understood componentwise.
As in the case of conservation laws, we limit the bar states of the momentum equation as follows (\cf \cite[Sec.~5.1]{kuzmin2020})
\begin{align}\label{eq-swe-v-constr}
\bar h_{ij}^*\vv_i^{\min} \le \hvij^{b,*} \coloneqq \hvij^b + \frac{\vec f_{ij}^{h\vv,*}}{2\dij} = \bar h_{ij}^*\bar\vv_{ij} + \frac{\vec g_{ij}^{h\vv,*}}{2\dij} \le \bar h_{ij}^* \vv_i^{\max},
\end{align}
where \smash{$\vec g_{ij}^{h\vv,*}$} is a limited counterpart of the flux
\begin{align*}
\vec g_{ij}^{h\vv} = \vec f_{ij}^{h\vv} + 2\dij \big( \hvij^b - \bar h_{ij}^* \bar\vv_{ij} \big).
\end{align*}
It is easy to verify that \smash{$\vec g_{ij}^{h\vv}+\vec g_{ji}^{h\vv}=\vec 0$} by construction.
This property must be preserved by the flux limiter.
From \eqref{eq-swe-v-constr} we derive the flux constraints for
\begin{align*}
\vec g_{ij}^{h\vv,*} = \begin{cases}
\spacingmin\min\lc
\vec g_{ij}^{h\vv},\; 2d_{ij}\spacingmin\min \lc \bar h_{ij}^* \l \vv_i^{\max} - \bar \vv_{ij}\r,\; \bar h_{ji}^* \l \bar \vv_{ij} - \vv_j^{\min}\r \rc \rc & \text{if } \vec g_{ij}^{h\vv} \geq 0, \\\spacingmax\max\lc
\vec g_{ij}^{h\vv},\; 2d_{ij}\spacingmax\max \lc \bar h_{ij}^* \l \vv_i^{\min} - \bar \vv_{ij}\r,\; \bar h_{ji}^*\l \bar \vv_{ij} - \vv_j^{\max}\r \rc \rc & \text{if }\vec g_{ij}^{h\vv}\le 0.
\end{cases}
\end{align*}
Note that this formula applies a scalar limiter of the form \eqref{eq-mcl-H} to each component of the flux vector $\vec g_{ij}^{h\vv}$.
Finally, we obtain the flux-corrected momentum bar states via
\begin{align*}
\vec f_{ij}^{h\vv,*} = \vec g_{ij}^{h\vv,*} - 2\dij\big( \hvij^b - \bar h_{ij}^* \bar\vv_{ij} \big), \qquad
\hvij^{b,*} = \hvij^b + \frac{\vec f_{ij}^{h\vv,*}}{2\dij}.
\end{align*}
What remains is to choose feasible bounds $\vv_i^{\min},\,\vv_i^{\max}$.
As for the SWE without topography source term, we should include the generalized velocity bar states $\bar{\vv}_{ij}$ in their definition.
To prove that a subsequent entropy fix based on \eqref{eq-mod-tadmor} cannot cause violations of \eqref{eq-swe-v-constr}, we need to extend the bounds by including the states \smash{$\overline{(h\vv)}_{ij}^b/\bar h_{ij}$}.
For reasons discussed in \cite[Lem.~3.16 \& Sec.~4.3.2]{hajduk2022}, we define the velocity bounds as follows
\begin{align*}
\vv_i^{\min} \coloneqq \min_{j\in \mathcal N_i\setminus \{i\}} \min \lc \bar \vv_{ij},~\frac{\overline{(h\vv)}_{ij}^b}{\bar h_{ij}} \rc, \qquad
\vv_i^{\max} \coloneqq \max_{j\in \mathcal N_i\setminus \{i\}} \max \lc \bar{\vv}_{ij},~\frac{\overline{(h\vv)}_{ij}^b}{\bar h_{ij}} \rc.
\end{align*}
Note that, contrary to the sequential approach for the SWE with flat bottom, the states \smash{$\overline{(h\vv)}_{ij}^b/\bar h_{ij}$} are neither symmetric nor antisymmetric in general.
Therefore, they need to be computed for all bar states, even the ones corresponding to pairs of interior nodes.

\subsection{Semi-discrete entropy fix}

At the current design stage, the semi-discrete bound-preserving scheme reads
\begin{align*}
m_i \frac{\d h_i}{\d t} ={}&
\afcsum \big[ 2\dij ( \hij^b - h_i ) + \fij^{h,*} \big] = 
\afcsum 2\dij (\hij^{b,*} - h_i),\\
m_i \frac{\d (h\vv)_i}{\d t} ={}& \afcsum \big[ 2\dij \big( \hvij^b - (h\vv)_i\big) + \vec f_{ij}^{h\vv,*}\big] = 
\afcsum 2\dij \big(\hvij^{b,*} - (h\vv)_i\big).
\end{align*}
To enforce a semi-discrete entropy inequality, we employ limiting coefficients $\beta_{ij}=\beta_{ji} \in [0,1]$ and entropy limited fluxes \smash{$\fij^{**}=\beta_{ij}\fij^*=\beta_{ij}(\fij^{h,*},(\vec f_{ij}^{h\vv,*})\transp)\transp$}.
Our approach represents a straightforward generalization of the entropy limiter used for conservation laws.
We adjust the Rusanov coefficients $\dij$ if necessary to guarantee that the low order method corresponding to $f_{ij}^*=0$ satisfies the entropy stability condition \smash{$\frac\dij2 P_{ij} \le Q_{ij}$}\ie \eqref{eq-mod-tadmor}.
The flux-corrected scheme with $\fij^*$ replaced by $\fij^{**}$ is entropy stable if
\begin{align}\label{eq-tadmor-swe}
\frac{\dij}2 P_{ij} + \frac{\beta_{ij}}{2} R_{ij} \le{}& Q_{ij},
\end{align}
where
\begin{align*}
R_{ij} \coloneqq \begin{bmatrix}
g[h_i - h_j + \lij^b(b_i - b_j) - \frac12 (|\vv_i|^2 - |\vv_j|^2)
\\ \vv_i - \vv_j
\end{bmatrix}\transp \fij^* = R_{ji}.
\end{align*}
Thus, we enforce \eqref{eq-tadmor-swe} by setting
\begin{align}\label{eq-entropy-lim-swe}
\beta_{ij} = \begin{cases} \displaystyle
\frac{2\min \{ Q_{ij}, Q_{ji} \} - \dij P_{ij}}{R_{ij}} & \text{if } R_{ij} > 2\min \{ Q_{ij}, Q_{ji}\} - \dij P_{ij}, \\
1 & \text{otherwise.}
\end{cases}
\end{align}
Since the Rusanov coefficients $\dij$ are chosen large enough for $2\min \{Q_{ij}, Q_{ji}\}\ge \dij P_{ij}$ to hold, \eqref{eq-entropy-lim-swe} produces $\beta_{ij}=\beta_{ji}\in[0,1]$.

The generalization of our monolithic limiting strategies for the SWE with topography is now complete.
Written in terms of the flux-corrected bar states
\begin{align*}
\bar h_{ij}^{b,**} = \bar h_{ij}^b + \frac{\beta_{ij}f_{ij}^{h,*}}{2\dij}, \qquad
\overline{(h\vv)}_{ij}^{b,**} = \overline{(h\vv)}_{ij}^b + \frac{\beta_{ij}\vec f_{ij}^{h\vv,*}}{2\dij},
\end{align*}
the resulting semi-discrete method reads
\begin{align*}
m_i \frac{\d h_i}{\d t} ={}&
\afcsum 2\dij \big( \hij^{b,**} - h_i\big),\\
m_i \frac{\d (h\vv)_i}{\d t} ={}&
\afcsum 2\dij \big(\hvij^{b,**} - (h\vv)_i\big).
\end{align*}
By construction, this finite element method is provably well balanced, bound preserving, and entropy stable.
We summarize its properties in the following theorem.

\begin{theorem}[Properties of flux correction schemes for \eqref{eq-swe-topo}]
The low order method and the flux-corrected schemes presented in this section
\begin{itemize}
\item[i)]
reduce to the corresponding algorithms discussed in \cref{sec:conslaw} if applied to the shallow water equations with flat topography,
\item[ii)]
produce nonnegative water heights under the CFL-like condition
\begin{align}\label{eq-cfl}
1 - \frac{\Delta t}{m_i}\afcsum{2\dij} \ge 0,
\end{align}
\item[iii)]
are well balanced for the lake at rest in the sense of Lemma~\ref{lem:lar}, and
\item[iv)]
satisfy the semi-discrete entropy inequalities \eqref{eq-swe-entropy-ineq} \wrt the entropy pair \eqref{eq-entropypair} if the correction factors $\beta_{ij}$ are either zero (low order method) or calculated using \eqref{eq-entropy-lim-swe}.
In the flux-corrected version of the scheme, the numerical fluxes $G_{ij}$ and consistency errors $W_{ij}$ appearing in \eqref{eq-swe-entropy-ineq} are replaced with
\begin{align*}
G_{ij}^* \coloneqq{}& G_{ij} + \frac{\beta_{ij}}{2} \big( \big[
g(h_i+b_i + h_j + b_j) - \frac12(|\vv_i|^2+|\vv_j|^2) \big] \fij^{h,*} + (\vv_i + \vv_j)\transp \vec f_{ij}^{h\vv,*} \big),\\
W_{ij}^* \coloneqq{}& W_{ij} + \frac g2(1-\lij^b)(b_i - b_j)\beta_{ij}\fij^{h,*},
\end{align*}
respectively.
Moreover, the statement of Corollary~\ref{cor:entropy-swe} remains valid.
\end{itemize}
\end{theorem}
\textbf{Proof:} See \cite[[Pf. of Thm. 4.9]{hajduk2022}. \hfill $\square$
\medskip

\section{Wetting and drying algorithms}\label{sec:wad}

Before moving on to numerical examples, we discuss some wetting and drying algorithms proposed in the literature and present our own approach.
Ricchiuto and Bollermann \cite[Sec.~4.3]{ricchiuto2009} set the velocity to zero if the water height is smaller than a prescribed tolerance for which they use the square of the normalized mesh size \smash{$(h/|\Omega|)^2$}.
Admittedly, their wetting and drying approach is more involved.
In particular, it incorporates information on the topography slope in wet-dry transition regions.
We have not tested this part of their algorithm but ran experiments with the version that sets the velocity to zero in dry regions.

A velocity fix that does guarantee continuous dependence on data is given by
\begin{align}\label{eq-agp-fix}
\tilde{\vv} =\frac{2h(h\vv)}{h^2 + \max\{h,\epsilon\}^2},
\end{align}
where $\epsilon \ll 1$.
Azerad \etal \cite{azerad2017} use \smash{$\epsilon=10^{-16}\max_{\x\in \overline \Omega }h_0(\x)$} in this formula.
A~problem with this approach is that if the water height approaches zero but the discharge does not, the use of \eqref{eq-agp-fix} produces velocities with magnitudes that tend to infinity, resulting in unrealistic CFL conditions and a blowup of kinetic energy.

Kurganov and Petrova suggest a similar fix \cite[Eqs.~(2.17),\,(2.21)]{kurganov2007a} in which the velocity is computed via
\begin{align}\label{eq-kp-fix}
\tilde{\vv} = \frac{\sqrt 2 h (h\vv)}{\sqrt{h^4 + \max\{h,\epsilon\}^4}}
\end{align}
and the parameter $\epsilon$ is set equal to the (normalized) mesh size.
In our experience, this choice introduces significant approximation errors because the mesh size is usually much larger than the thickness of a water layer that can be considered as dry.
On the other hand, this fix seems to be quite robust in practice.
Importantly, Kurganov and Petrova \cite[Eq.~(2.21)]{kurganov2007a} emphasize the need for adjusting the discharge by setting $(h\vv) = h\tilde{\vv}$ after calculating the velocity via \eqref{eq-kp-fix}.

Many more algorithms for wetting and drying processes exist besides the ones already mentioned.
Most of them work in a fashion similar to the approaches discussed above.
There are also schemes that can not directly be applied in the context of continuous finite elements.
For instance, Vater \etal \cite{vater2015} employ slope limiters to handle wetting and drying scenarios.

Let us now discuss a new nodal velocity correction based on the entropy of the shallow water system.
Here we restrict ourselves to the case of a flat topography because it is currently unclear to us whether an extension to the general case is feasible for our MCL schemes.
The underlying idea is based on the observation that unbounded velocities, which may occur in dry or nearly dry areas, result in blow ups of the kinetic energy and, therefore, of the entropy.
On the other hand, entropy analysis of the bar state form for the SWE with flat bathymetry provides an upper bound for the nodal entropy.
Violations of this bound (and the resulting lack of discrete entropy stability in practice) are caused not by the discretization but by the numerically unstable calculation of nodal velocities for the next step.
Thus, for entropy stability reasons, the magnitudes of nodal velocities should stay bounded in the vicinity of dry states.
In physics, this property is enforced by viscous friction, which is missing in our model.

Recall once more that for $b\equiv 0$, an entropy for the shallow water equations is the sum of potential and kinetic energies
(\cf \eqref{eq-entropypair}).
By convexity, the entropy of the state \smash{$\tilde u_i = (\tilde h_i, \widetilde{(h\vv)}_i)$} produced by a forward Euler update satisfies the estimate
\begin{align}
\eta(\tilde u_i) ={}& \eta \big( \big(1 - \frac{\Delta t}{m_i} \afcsum 2\dij \big) u_i + \frac{\Delta t}{m_i}\afcsum 2\dij \uij^{**} \big)
\notag \\ \le{}&
\big(1 - \frac{\Delta t}{m_i} \afcsum 2\dij \big) \eta(u_i) +
\frac{\Delta t}{m_i}\afcsum 2\dij \eta(\uij^{**}) \eqqcolon \eta_i^{\max}
\label{eq-entropy-bound}
\end{align}
under the CFL condition \eqref{eq-cfl}.
The value $\eta_i^{\max}$ can now be used to prohibit the occurrence of unbounded velocities that would lead to a violation of \eqref{eq-entropy-bound}.
Invoking the definition \eqref{eq-entropypair} of $\eta$, we enforce \eqref{eq-entropy-bound} by adjusting the nodal velocities as follows
\begin{align*}
\tilde \vv_i = \begin{cases}
\frac{\widetilde{(h\vv)}_i}{\tilde h_i} & \text{if } |\widetilde{(h\vv)}_i| \le h_i Q_i, \\
\frac{Q_i}{|\widetilde{(h\vv)}_i|} \widetilde{(h\vv)}_i & \text{if } |\widetilde{(h\vv)}_i| > h_i Q_i,
\end{cases}\qquad\text{where}\qquad Q_i \coloneqq \sqrt{\frac{2\eta_i^{\max}}{\tilde h_i} - g\tilde h_i}.
\end{align*}
We then follow Kurganov and Petrova \cite[Eq.~(2.21)]{kurganov2007a} and overwrite the nodal discharge by \smash{$\tilde h_i \tilde \vv_i$}.
The approach presented here for a forward Euler update directly carries over to other SSP RK methods, which are convex combinations of forward Euler steps.
Unfortunately, the entropy-based approach interferes with the well-balancedness property for the lake at rest unless the topography is flat.

To keep our scheme well balanced, we developed a wetting and drying algorithm that is based on the theory of laminar boundary layers (see for instance \cite{schlichting2017}).
As suggested by the above discussion of our entropy-based approach, a particular challenge for realistic treatment of wetting and drying processes is to obtain a physically correct model for the velocities in wet-dry transition regions.
According to the boundary layer theory, viscous friction effects should not be neglected in these areas.
For the SWE in particular, a bottom friction term should be incorporated into the system.
A derivation of the viscous SWE including bottom friction can be found in \cite{gerbeau2000}.
In essence, a nonconservative term $\sigma \vv$ is added on the left hand side of the momentum equation.
Here $\sigma>0$ is the bottom friction coefficient, which may generally depend on the solution and parameters of the SWE\@.
Particular models for $\sigma$ are discussed, for instance, in \cite[Sec.~2.7]{vreugdenhil1994} and \cite[Sec.~9.8]{cushman-roisin2011}.
Physical intuition tells us that wet-dry transitions occur in a boundary layer of thickness $0<\delta\ll 1$.
According to \cite[Ch.~2]{schlichting2017}, one may assume that inertial and viscous forces are in equilibrium and contributions of the material derivative can be neglected in the boundary layer, which in our case implies
\begin{align*}
g h \nabla H + \sigma \vv = 0.
\end{align*}
For nodes belonging to wet-dry zones\ie for $h_i \le \delta$, we use this identity to compute a nodal boundary layer velocity $\vv_i^{\text{BL}}$ via the lumped $\mathrm L^2$ projection
\begin{align}\label{eq-bl-fix}
m_i \vv^{\text{BL}}_i = - \frac g \sigma\,
h_i \sum_{j \in \mathcal N_i}^N H_j\, \cij.
\end{align}
Then the nodal velocity \smash{$\tilde \vv_i = (h\vv)_i/h_i$} is adjusted as follows
\begin{align}\label{eq-bl-fix2}
\tilde \vv_i = \frac{(h\vv)_i}{\max\{h_i,\delta\}} + \max\lc 0, \frac{\delta - h_i}{\delta} \rc \vv_i^{\text{BL}}.
\end{align}
Finally, we overwrite the discharge by \smash{$h_i \tilde \vv_i$} as in the energy-based version and in the algorithm proposed by Kurganov and Petrova \cite[Eq.~(2.21)]{kurganov2007a}.

Note that formula \eqref{eq-bl-fix2} ensures a continuous transition between \smash{$\tilde \vv_i$} for $h_i\in[0,\delta]$ and \smash{$\vv_i^{\text{SWE}} \coloneqq (h\vv)_i/h_i$} for $h_i\ge\delta$.
Similarly to wall function models for turbulent flows, it uses the inviscid SWE model for $h_i>\delta$ but adapts (the solution of) the momentum equation in the boundary layer, where viscous friction effects are dominant and some assumptions behind the derivation of the shallow water equations are invalid.

In all of the numerical experiments below, we set the bottom friction parameter and the boundary layer thickness to $\sigma=10$ and $\delta=10^{-3}$, respectively.
These values are chosen according to \cite{gerbeau2000} and the boundary layer theory \cite[Ch.~2]{schlichting2017}.
Note that our boundary layer has a thickness of 1 millimeter for the SWE without nondimensionalization.
In our opinion this constitutes a reasonable value for which a nodal state can be considered as almost dry and friction should come into play.

Barros \etal \cite[Sec.~3]{barros2015} propose an approach that models the impact of bottom friction on wetting and drying in a different way.
The underlying idea is to treat the sea bed as a porous medium.
Based on precomputed water depths, the authors of \cite{barros2015} distinguish between wet, dry, and transitional regions.
The water height and the bottom friction coefficient are adjusted in the latter two regimes.
Since we consider the SWE without a bottom friction term (our own fix based on boundary layer theory only adjusts the velocity and discharge), this approach cannot be directly pursued here.

\section{Numerical examples}\label{sec:swe-num}

Let us now apply the generalized flux correction schemes for the shallow water equations with topography to various one-dimensional benchmarks.
By default, we employ a uniform mesh consisting of 128 elements and adaptive SSP2 RK time stepping (Heun's method) with CFL parameter $\nu=0.5$, which is used to compute the time step $\Delta t$ in the first SSP RK stage via
\begin{align*}
\Delta t = \min_{i\in \{1,\ldots,N\}} \frac{\nu\,m_i}{\afcsum{2\dij}}.
\end{align*}
Below we consider classical steady state examples and various dam break problems, before testing our algorithms for an idealized parabolic lake.
The strategies under investigation include the algebraic local Lax--Friedrichs scheme\ie the low order method (LOW), the bound-preserving monolithic convex limiting approach without entropy fix (MCL), as well as the MCL algorithm enhanced by the semi-discrete entropy limiter (MCL-SDE).

\subsection{Steady problems}

We investigate the well-balancedness of our schemes by applying them to an exact lake at rest configuration as well as moving water equilibria.
By default, we employ the raw antidiffusive fluxes $\fij=\dij(u_j-u_i)$ for \afcijtext~in this section.
This choice is suitable for steady and weakly time-dependent problems.

\subsubsection{Lake at rest}

In our first test, we set $\Omega=(0,1)$, $g=1$, and $b(x) = \max(0, 0.25 - 5(x-0.5)^2)$.
Our initial conditions read $v_0 \equiv 0$ and $h_0(x) = \max\{0.2\mathrm H(0.5-x)+0.1\mathrm H(x-0.5), b(x)\}$, where $\mathrm H$ is the Heaviside function.
This test problem represents a lake at rest configuration, and is essentially the same as in \cite[Sec.~4.2]{liang2009} but many similar benchmarks exist in the literature.
In our case, there are two bodies of water with different depths that are separated from each other by a land mass.
Boundary conditions are realized as reflecting walls.
However, this choice does not affect the numerical results.

\begin{figure}[ht!]
\centering
\begin{subfigure}[b]{0.32\textwidth}
\caption{Free surface elevation}
\includegraphics[width=\textwidth]{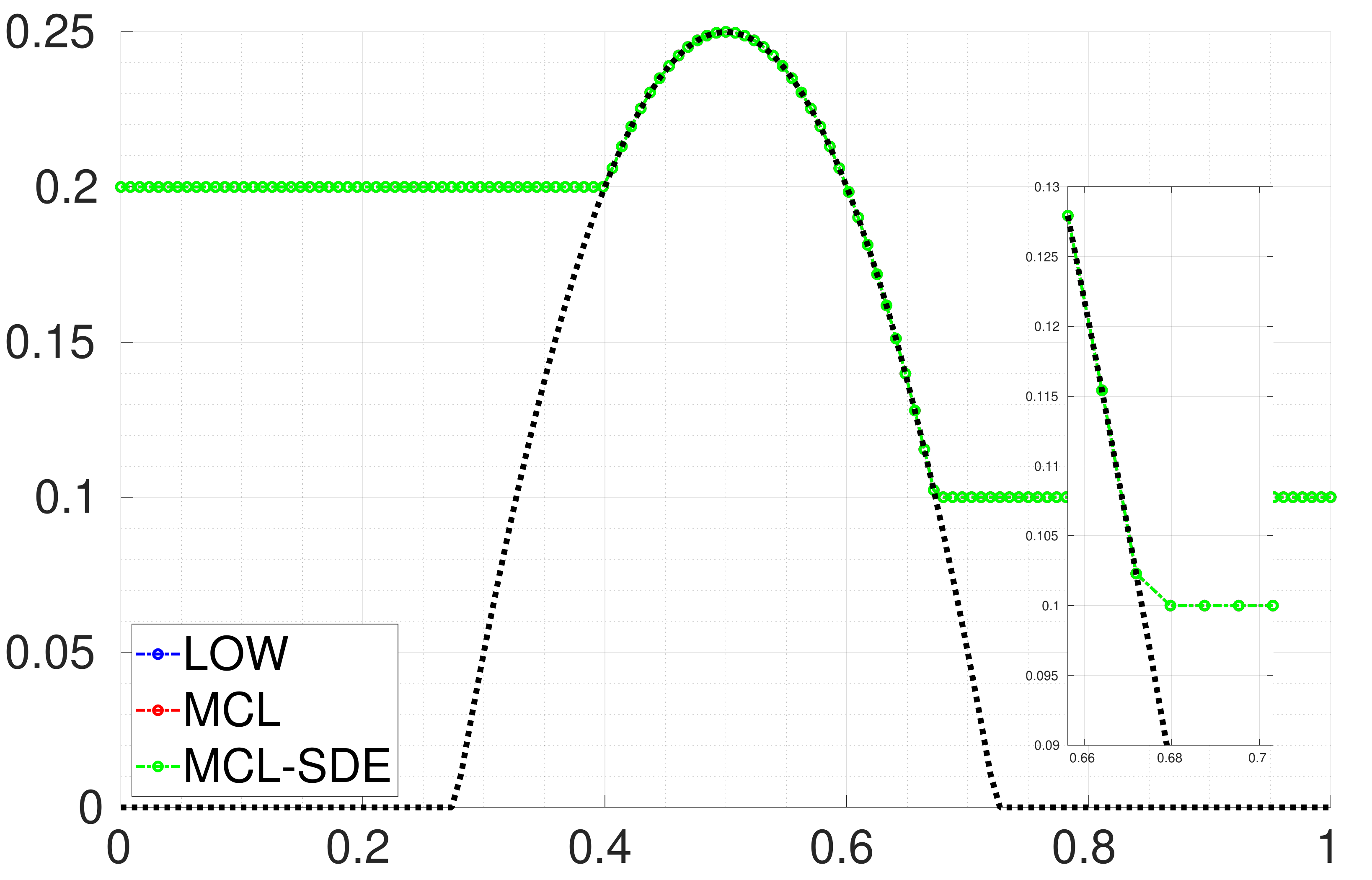}\label{fig:lar-1}
\end{subfigure}
\begin{subfigure}[b]{0.32\textwidth}
\caption{Discharge}
\includegraphics[width=\textwidth]{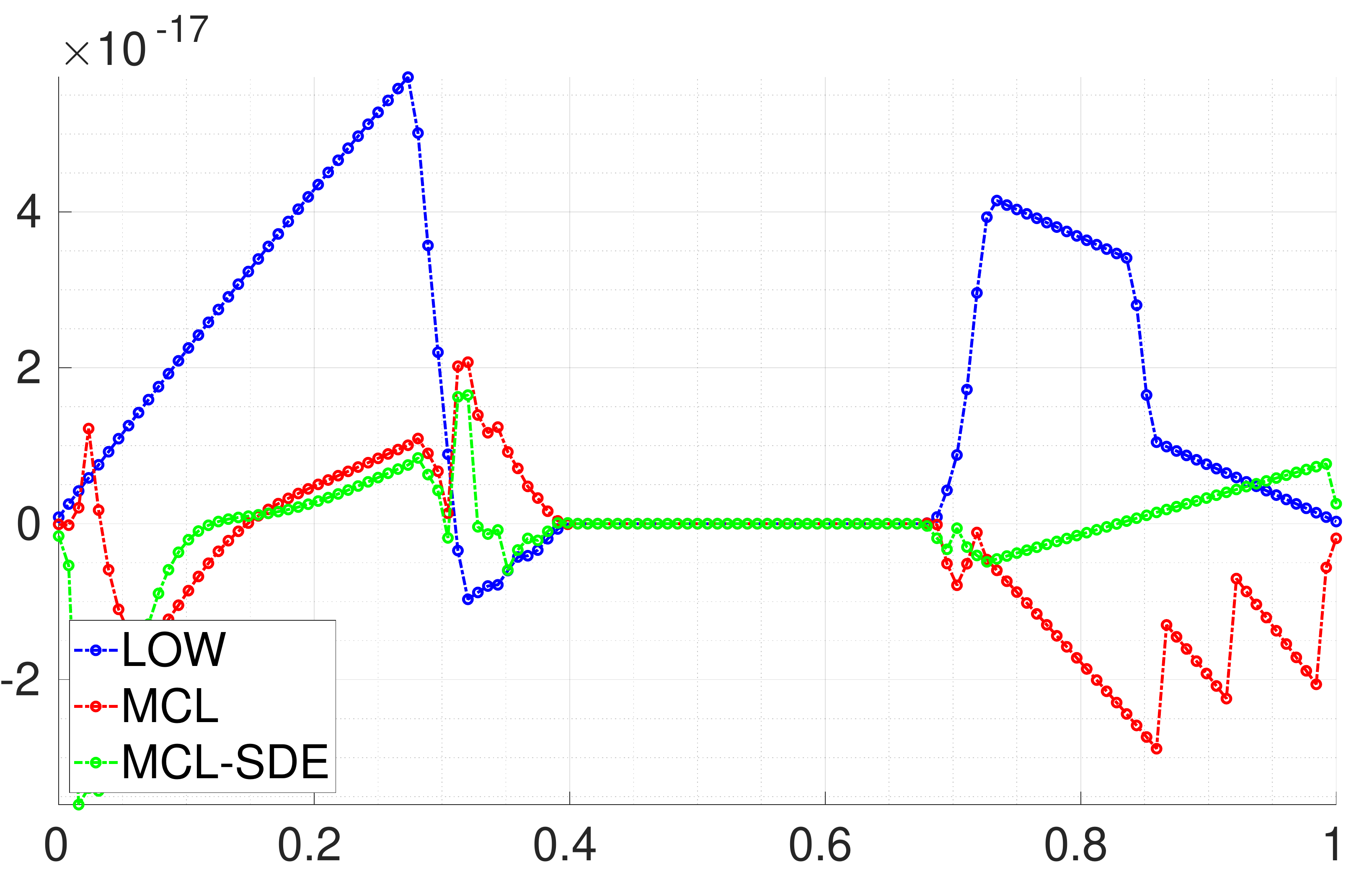}
\end{subfigure}
\begin{subfigure}[b]{0.32\textwidth}
\caption{Velocity}
\includegraphics[width=\textwidth]{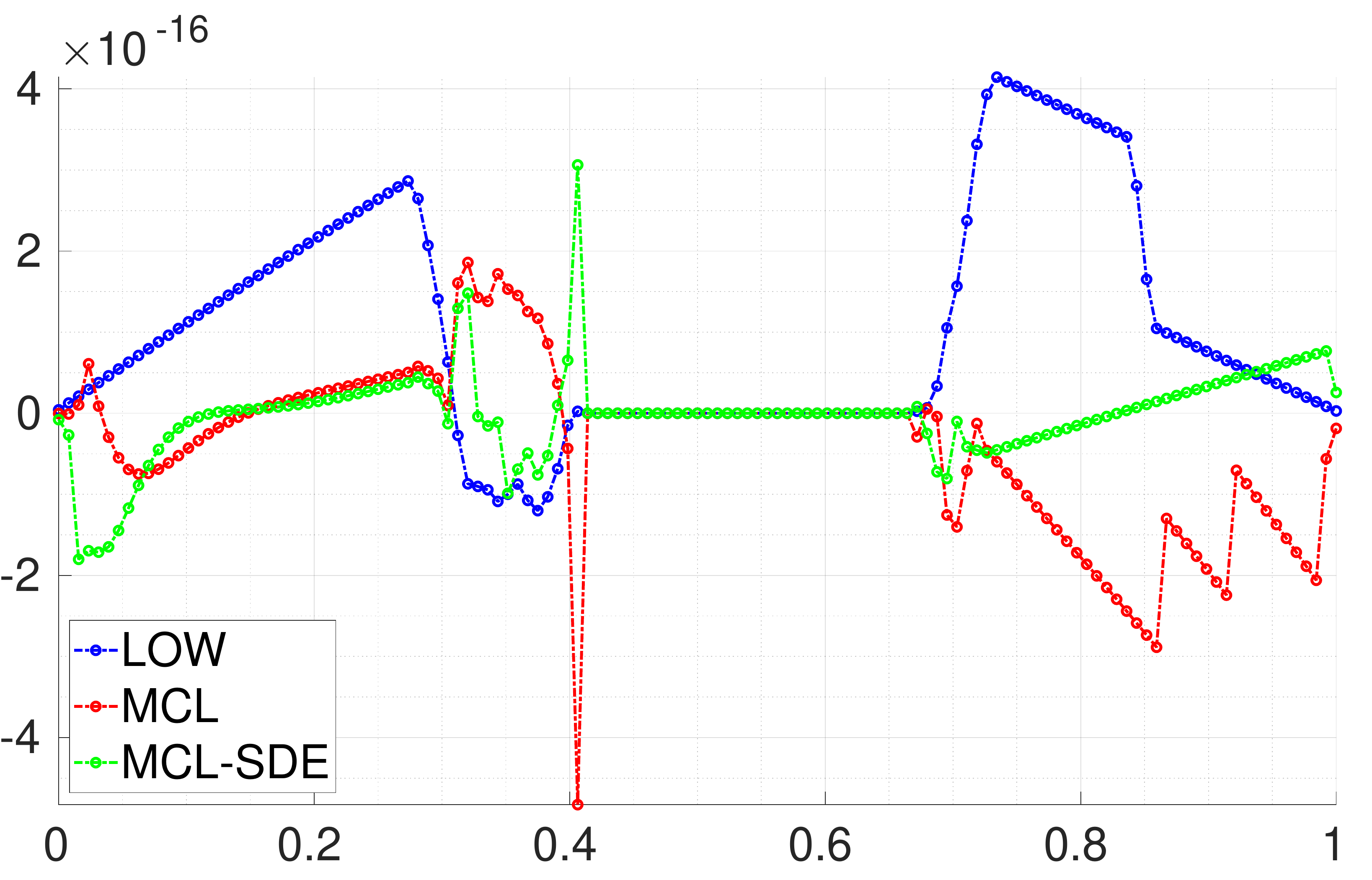}
\end{subfigure}
\caption{A lake at rest for the shallow water equations.
Approximations at $T=100$ obtained with adaptive SSP2 RK time stepping and $\nu=0.5$ on a uniform mesh consisting of 128 elements.}\label{fig:lar}
\end{figure}

We solve this problem numerically using the boundary layer-based wetting and drying approach \eqref{eq-bl-fix}--\eqref{eq-bl-fix2} and display the results in \cref{fig:lar}.
All methods clearly preserve the lake at rest scenario, which is why the free surface elevation profiles in \cref{fig:lar-1} are perfectly on top of each other.
Note that the oscillations observable in discharge and velocity are of the order of machine precision and do not amplify in the course of the simulation.
The stability of the approach becomes evident by realizing that a total of 22898 time steps was performed with each scheme to reach the very large end time $T=100$.

In this test, the employed combination of mesh and initial conditions does not capture the shorelines exactly.
Thus, the well-balancedness condition \eqref{eq-wb-cond} does not apply here.
As a consequence, elements containing a wet-dry transition can only resolve the free surface in these cells by introducing an artificial slope in the discrete water heights.
One can see this artifact in the zoomed region of \cref{fig:lar-1}.
It was our intention to show that, in practice, our methods remain well balanced for such practical examples, even if \eqref{eq-wb-cond} does not hold.
We also ran a similar experiment where \eqref{eq-wb-cond} is satisfied and well-balancedness is guaranteed by Lemma~\ref{lem:lar}.
For such problems our schemes preserve the lake at rest configuration up to machine precision without introducing any nonphysical slopes in the water height.
The discharge and velocity profiles obtained in this fashion are similar to those in \cref{fig:lar}.

\subsubsection{Moving water equilibria}

Next, we study three classical steady benchmarks \cite[Sec.~5.3]{vazquez-cendon1999}, \cite[Sec.~3.1]{delestre2013} as well a supercritical modification thereof.
In all cases, the spatial domain is $\Omega=(0,25)$ and the bathymetry is set to $b(x) = \max\{0, 0.2-0.05(x-10)^2\}$.
At first, we assume that no nondimensionalization has been performed, thus we use $g=9.81$.

In the first example, we employ $(h_0,(hv)_0) \equiv (2,0)$ as initial condition and prescribe $(hv)_{\text{in}}=4.42$ at the subcritical inlet on the left and $h_{\text{in}}=2$ at the subcritical outlet on the right.
In fact the flow is subcritical everywhere and the treatment of boundaries is in accordance with the physics of the problem.
The exact solution for this setup can be computed as discussed in \cite[Sec.~3.1]{delestre2013}.
As a result of the bump in the bathymetry, there appears a corresponding one in the free surface elevation.

Next, we consider a transcritical flow example without a shock, which is obtained with the initial and left boundary data $(h_0,(hv)_0) \equiv (0.66,0)$ and $(hv)_{\text{in}}=1.53$, respectively.
In this example, the type of the right boundary changes in time and is determined numerically, by computing the eigenvalues of the flux Jacobian for the internal state at the boundary.
The external boundary state is then set according to the underlying physics based on the available boundary data $h_{\text{in}}=0.66$ and $(hv)_{\text{in}} = 1.53$.
In this example, the flow becomes supercritical (this behavior is referred to as \textit{torrential} flow) at the bathymetry bump and to the right of it but remains subcritical at the left domain boundary.

\begin{figure}[ht!]
\centering
\begin{subfigure}[b]{0.32\textwidth}
\caption{Subcritical flow}\label{fig:subcritical}
\includegraphics[width=\textwidth]{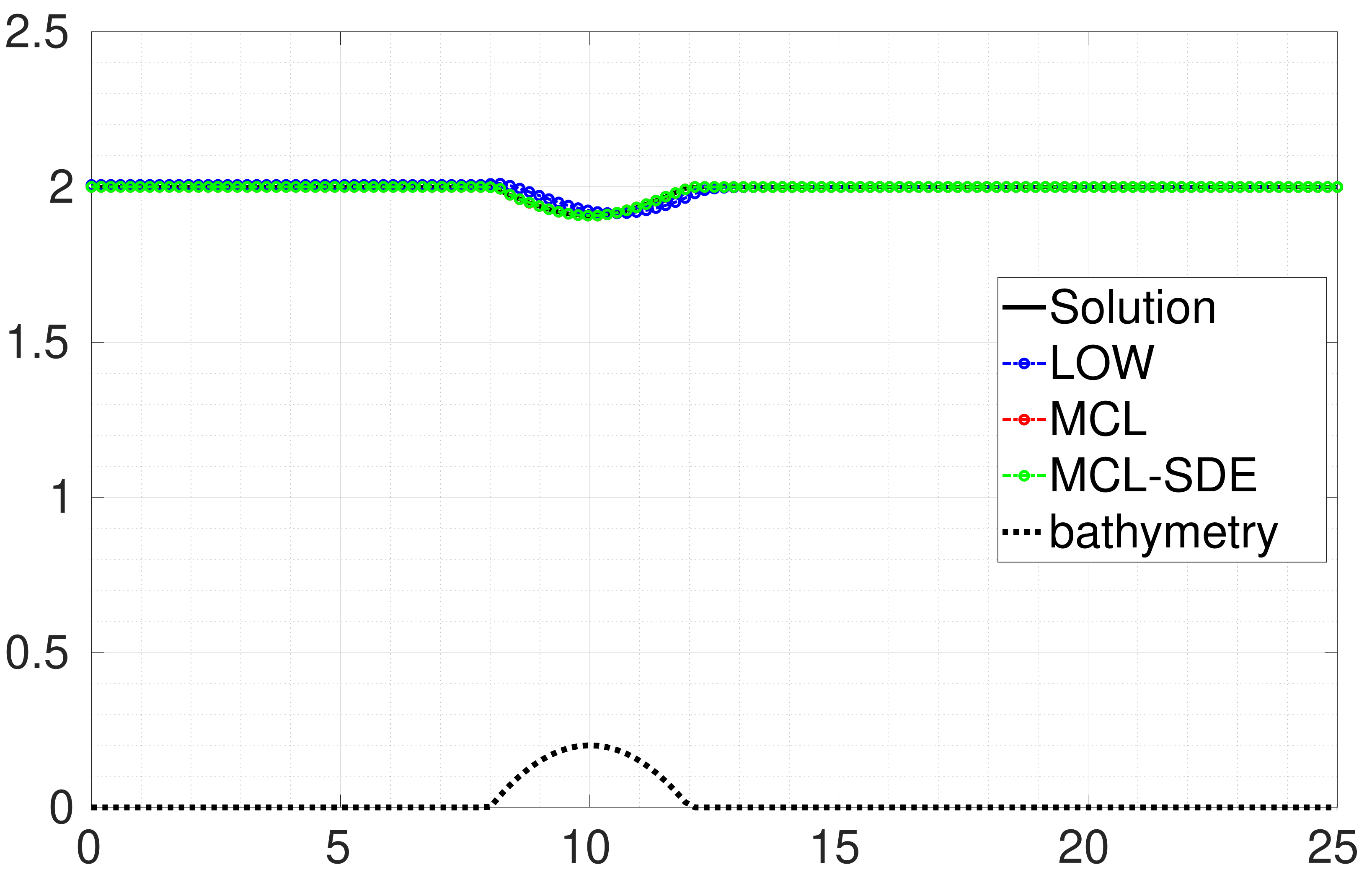}
\end{subfigure}
\begin{subfigure}[b]{0.32\textwidth}
\caption{Transcritical flow w/o shock}
\includegraphics[width=\textwidth]{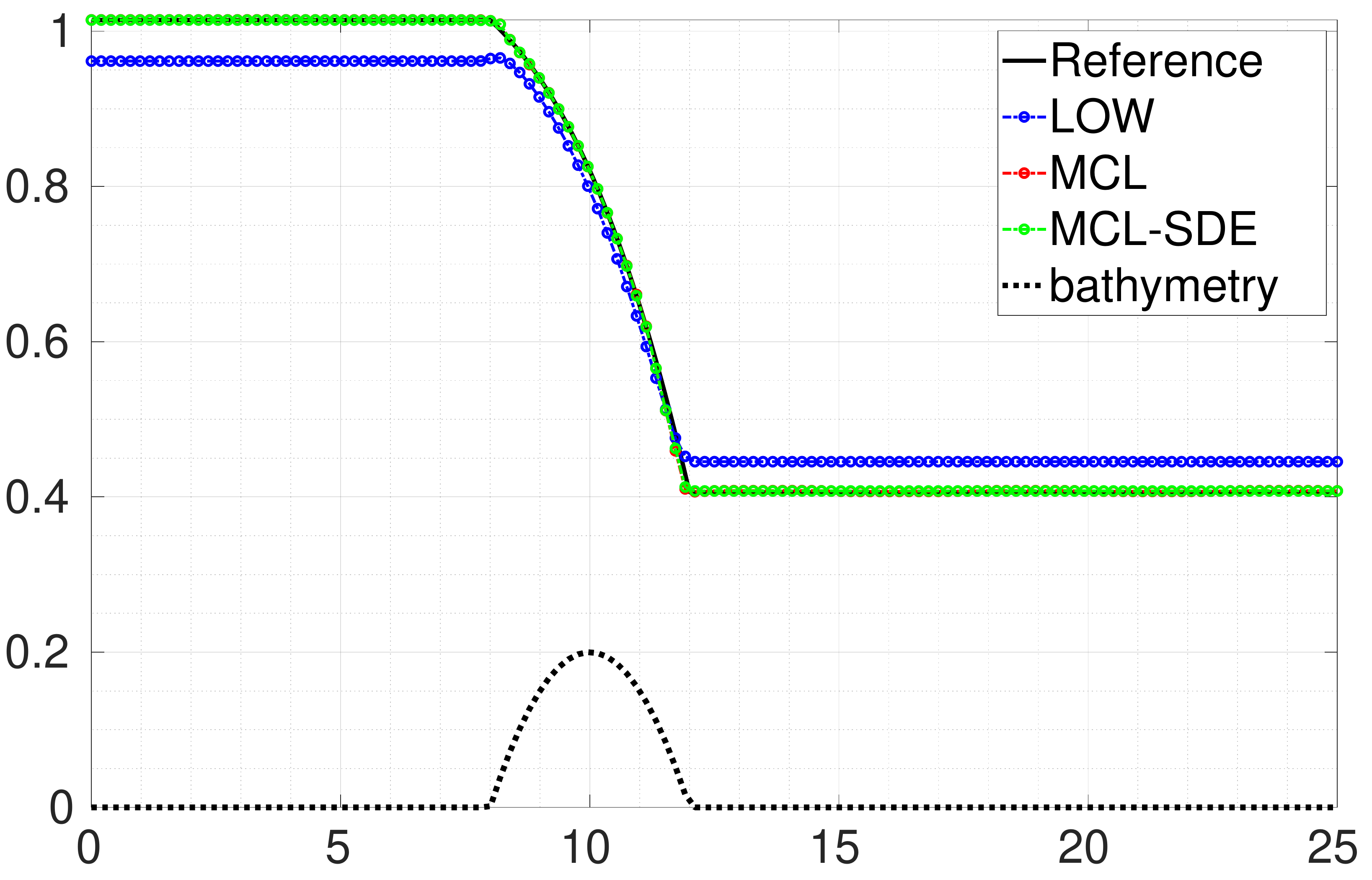}
\end{subfigure}
\begin{subfigure}[b]{0.32\textwidth}
\caption{Transcritical flow with shock}\label{fig:steady-c}
\includegraphics[width=\textwidth]{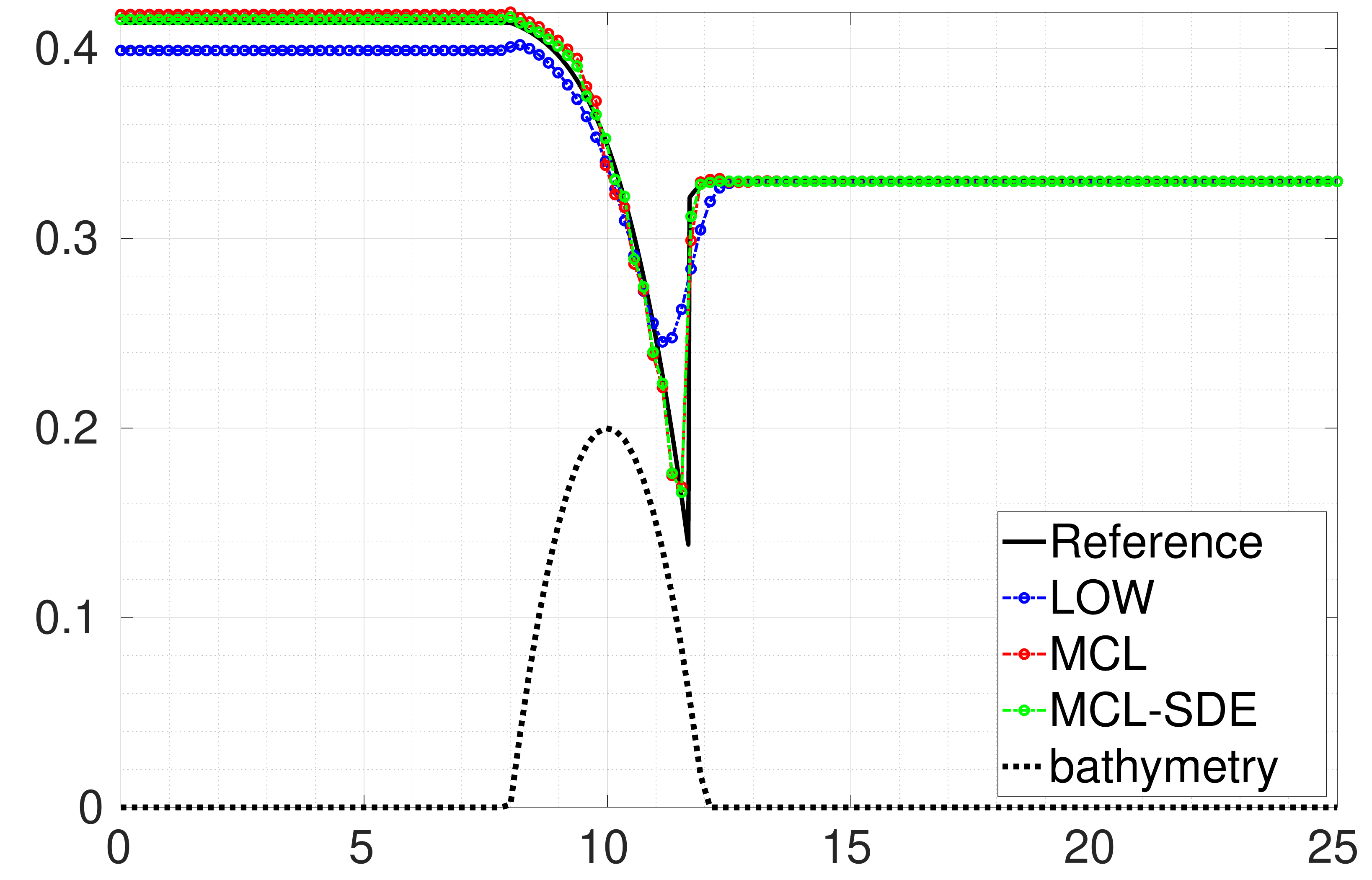}
\end{subfigure}
\caption{Moving water equilibria for the shallow water equations \cite{vazquez-cendon1999}.
Approximations to the free surface elevation at $T=400$ (a), $T=200$ (b), and $T=800$ (c) obtained with adaptive SSP2 RK time stepping and $\nu=0.5$ on a uniform mesh consisting of 128 elements.}\label{fig:steady}
\end{figure}

Another transcritical example is obtained by setting initial and boundary data as $(h_0,(hv)_0) \equiv (0.33,0)$, $(hv)_{\text{in}}=0.18$ on the left and $h_{\text{in}}=0.33$ on the right, respectively.
Again, the region around the bathymetry bump becomes torrential, and this time, a steady shock forms.
In this example however, the flow is subcritical, not only on the left of the area with elevated topography but also in the post-shock region.
Respective reference solutions obtained with the SWASHES software \cite{delestre2013} on uniform meshes of 1\,000 cells are displayed in \cref{fig:steady} for both transcritical examples.

The obtained free surface elevations for the three test problems are displayed in \cref{fig:steady}.
With the employed resolution, one can clearly see that the LOW profiles do not quite attain the exact values in regions where the exact solutions are constant.
As expected, LOW also smears the shock in \cref{fig:steady-c} significantly.
On the other hand, the agreement of the flux-corrected approximations with the respective exact or reference solutions is satisfactory.

All three of the above examples are classical steady benchmarks. 
Thus, we check, whether the approximations converge to steady states.
Unfortunately, this is only the case for the low order method, not for the flux-limited schemes.
A variety of reasons for this lack of convergence can be imagined.
In these examples it can be due to the fact, that our schemes are not exactly well balanced for moving water equilibria.

In our final steady example, we modify the above configurations by assuming the system to be in nondimensional form.
Thus, we set $g=1$.
As initial condition we use $(h_0,(hv)_0) \equiv (1,2.1)$, which corresponds to supercritical flow.
Thus, supercritical in- and outlet boundary conditions are prescribed at $x=0$ and $x=25$, respectively.
Again, the bathymetry bump produces a corresponding feature in the free surface elevation.
Contrary to the subcritical case displayed in \cref{fig:subcritical}, the bump is pointing upward in this example.
The exact solution to this problem can be derived as in \cite[Sec.~3.1]{delestre2013}.

\begin{figure}[ht!]
\centering
\begin{subfigure}[b]{0.32\textwidth}
\caption{Free surface elevation}\label{fig:supercritical-1}
\includegraphics[width=\textwidth]{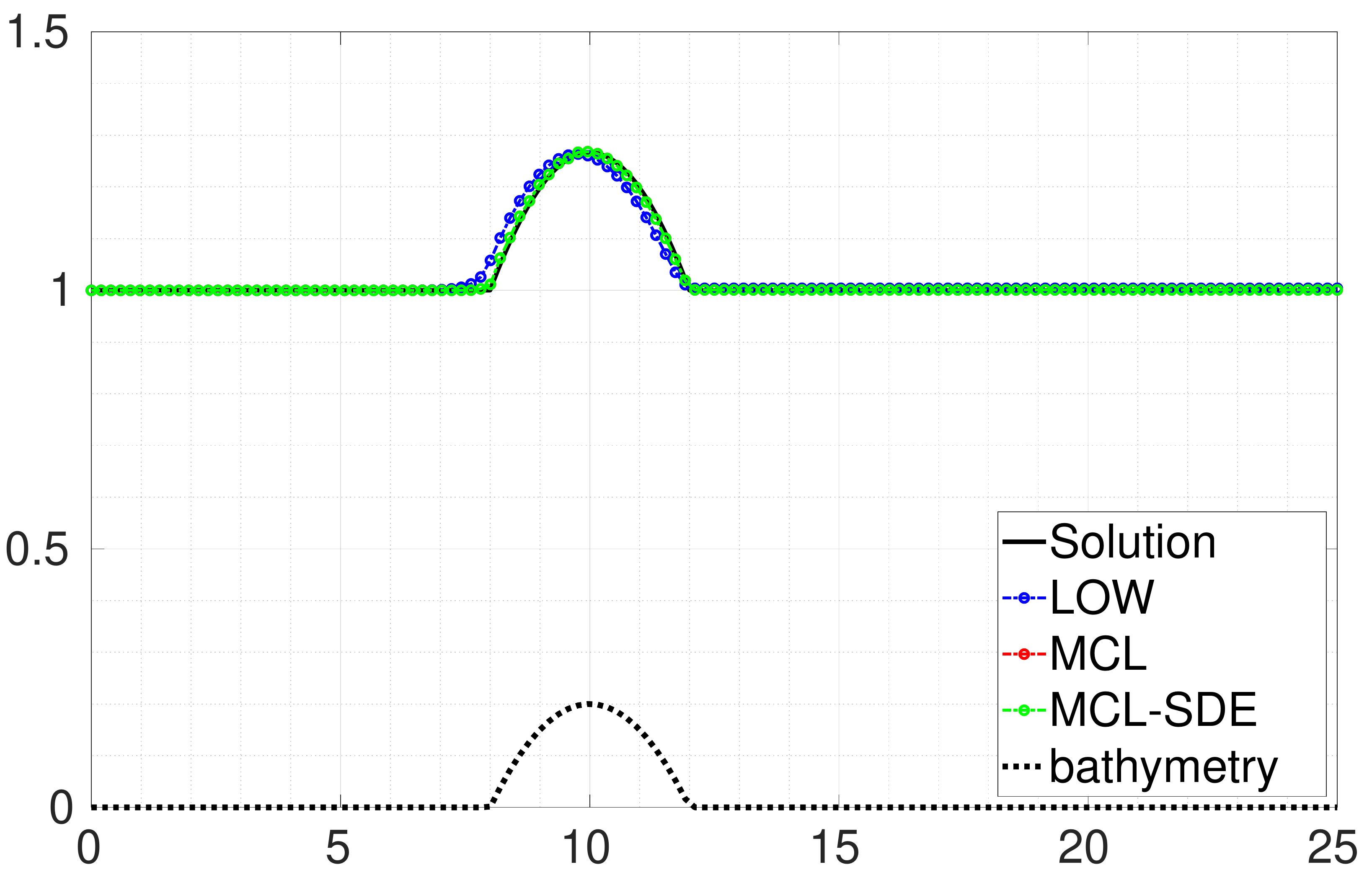}
\end{subfigure}
\begin{subfigure}[b]{0.32\textwidth}
\caption{Discharge}
\includegraphics[width=\textwidth]{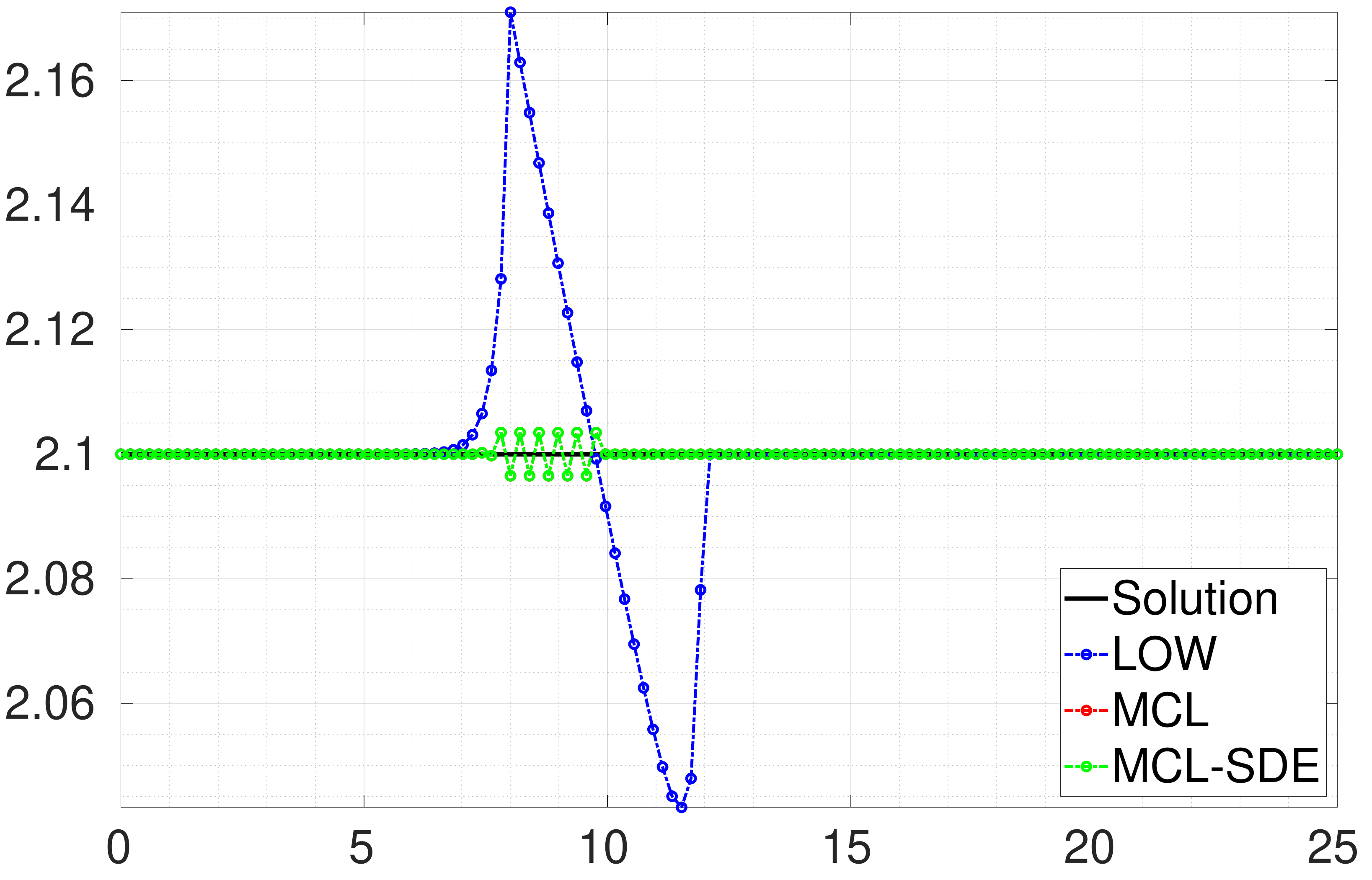}
\end{subfigure}
\begin{subfigure}[b]{0.32\textwidth}
\caption{Velocity}
\includegraphics[width=\textwidth]{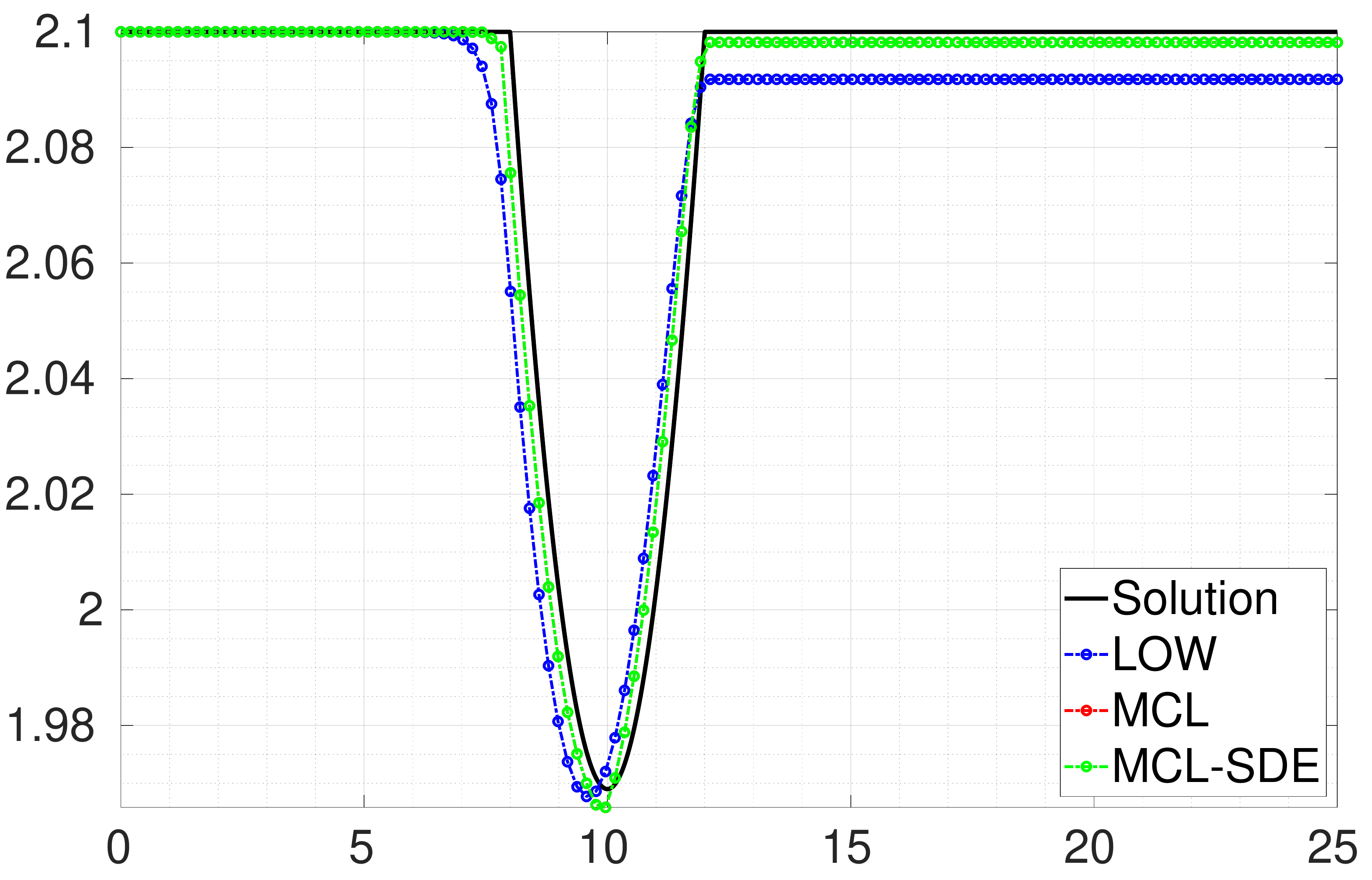}
\end{subfigure}
\caption{Supercritical moving water equilibrium for the shallow water equations.
Approximations at steady state obtained with adaptive SSP2 RK time stepping and $\nu=0.5$ on a uniform mesh consisting of 128 elements.}\label{fig:supercritical}
\end{figure}

In this example, all three schemes under investigation do converge to the steady states displayed in \cref{fig:supercritical}.
These profiles were obtained on a uniform mesh consisting of 128 elements.
To rule out that these are isolated instances, we increased the spatial and temporal resolutions by factors of two and four.
The steady state residual in each run eventually drops below the threshold of $10^{-12}$, although our schemes are not exactly well balanced for moving water equilibria.
It is quite remarkable, that the displayed oscillatory discharge profiles represent discrete steady states.
Nevertheless, a combination of AFC with strategies that guarantee well-balancedness \wrt to moving water equilibria (for instance, as in \cite{noelle2007}) is an important topic for future research.

If we include the term \smash{$\mij (\dot u_i^{\mathrm L} - \dot u_j^{\mathrm L})$} in the antidiffusive fluxes $\fij$, the MCL-SDE method still converges to steady state for all three resolutions under consideration, while the MCL method without entropy fix does not.
Therefore, it may be a good idea to employ an entropy fix in practical computations.

The fact that the low order discharge in \cref{fig:supercritical} looks quite different from its flux-corrected counterparts is due to the term \smash{$\frac12(\vv_i+\vv_j)\lij^b(b_j-b_i)$} present in the low order method.
Its influence on the numerical approximations is reduced in the process of flux limiting.
A preliminary version of our low order method without this term did not produce the slight phase errors visible in the low order approximations in \cref{fig:supercritical-1} and \cref{fig:subcritical}.
It is somewhat interesting that these deviations from the respective exact solutions are opposite in the subcritical and supercritical cases.

\subsection{Dam breaks}\label{sec:dam}

Having studied some problems with steady state solutions in the previous section, we now perform experiments for the generalized Riemann problem
\begin{align*}
\frac{\partial}{\partial t}\begin{bmatrix}
h\\hv
\end{bmatrix} + \frac{\partial}{\partial x}
\begin{bmatrix}
hv \\ hv^2 + \frac g 2h^2
\end{bmatrix} + \begin{bmatrix}
0\\gh\frac{\partial b}{\partial x}
\end{bmatrix} ={}& 0 \qquad \text{in } \Omega \times (0,T),\\
h_0(x) = \begin{cases}
h_L & \text{if } x < x_0,\\
h_R & \text{if } x > x_0,
\end{cases}
\qquad v_0 \equiv{}& 0,
\end{align*}
where $\Omega \subset \R$.
We use values $h_L > h_R$ in the three below tests.
This setup corresponds to an idealized dam located at $x_0\in \Omega$ that is removed at time $t=0$.
As a result a water wave propagates into the region $x> x_0$, while a rarefaction wave travels in the opposite direction.

\subsubsection{Wet dam break over flat topography}\label{sec:wet}

First, we set $g=1$ and consider an example with flat bottom topography\ie $b\equiv 0$.
Thus, we may apply the standard limiting techniques for conservation laws, instead of their generalized versions for the SWE with a topography source term.
If both $h_L$ and $h_R$ are positive, the generalized Riemann problem is referred to a wet dam break.
Such tests represent relatively mild test cases, which are similar to Sod's shock tube problem \cite{sod1978} for the Euler equations.
A difference is that there is one fewer unknown in the system, and the exact solution does not feature any contact discontinuities.

We equip the spatial domain $\Omega=(0,1)$ with reflecting wall boundaries (although other options are feasible).
In our first test, the dam is located at $x_0=0.5$ and the two values for the water height are set to $h_L=1$ and $h_R=0.1$.
As end time we choose $T=0.3$.
The exact solution to this problem can be found in \cite[Sec.~4.1.1]{delestre-arxiv}.

First, we perform a convergence study of LOW, MCL, and MCL-SDE schemes on a series of uniform meshes.
The convergence rates of the MCL and MCL-SDE approaches observed in \cref{tab:3-swe}, are optimal for examples such as this one in which the solution is discontinuous.

\begin{table}[ht!]
\centering
\begin{tabular}{c|cc|cc|cc}
$1/h$ & LOW & EOC & MCL & EOC & MCL-SDE & EOC \\
\hline
32& 7.93E-02 & & 3.28E-02 & & 3.66E-02 \\
64& 4.98E-02 & 0.67 & 1.67E-02 & 0.97 & 1.89E-02 & 0.95 \\
128 & 3.00E-02 & 0.73 & 8.47E-03 & 0.98 & 9.59E-03 & 0.98 \\
256 & 1.77E-02 & 0.76 & 4.28E-03 & 0.99 & 4.85E-03 & 0.98 \\
512 & 1.06E-02 & 0.75 & 1.94E-03 & 1.14 & 2.24E-03 & 1.11 \\
\end{tabular}
\caption{Convergence history of the wet dam break for the shallow water equations. The $\|\cdot\|_{\mathrm L^1(\Omega)}$~errors at $T = 0.3$ and the corresponding EOC.}
\label{tab:3-swe}
\end{table}

\begin{figure}[ht!]
\centering
\begin{subfigure}[b]{0.32\textwidth}
\caption{Water level}\label{fig:3-ssp2}
\includegraphics[width=\textwidth]{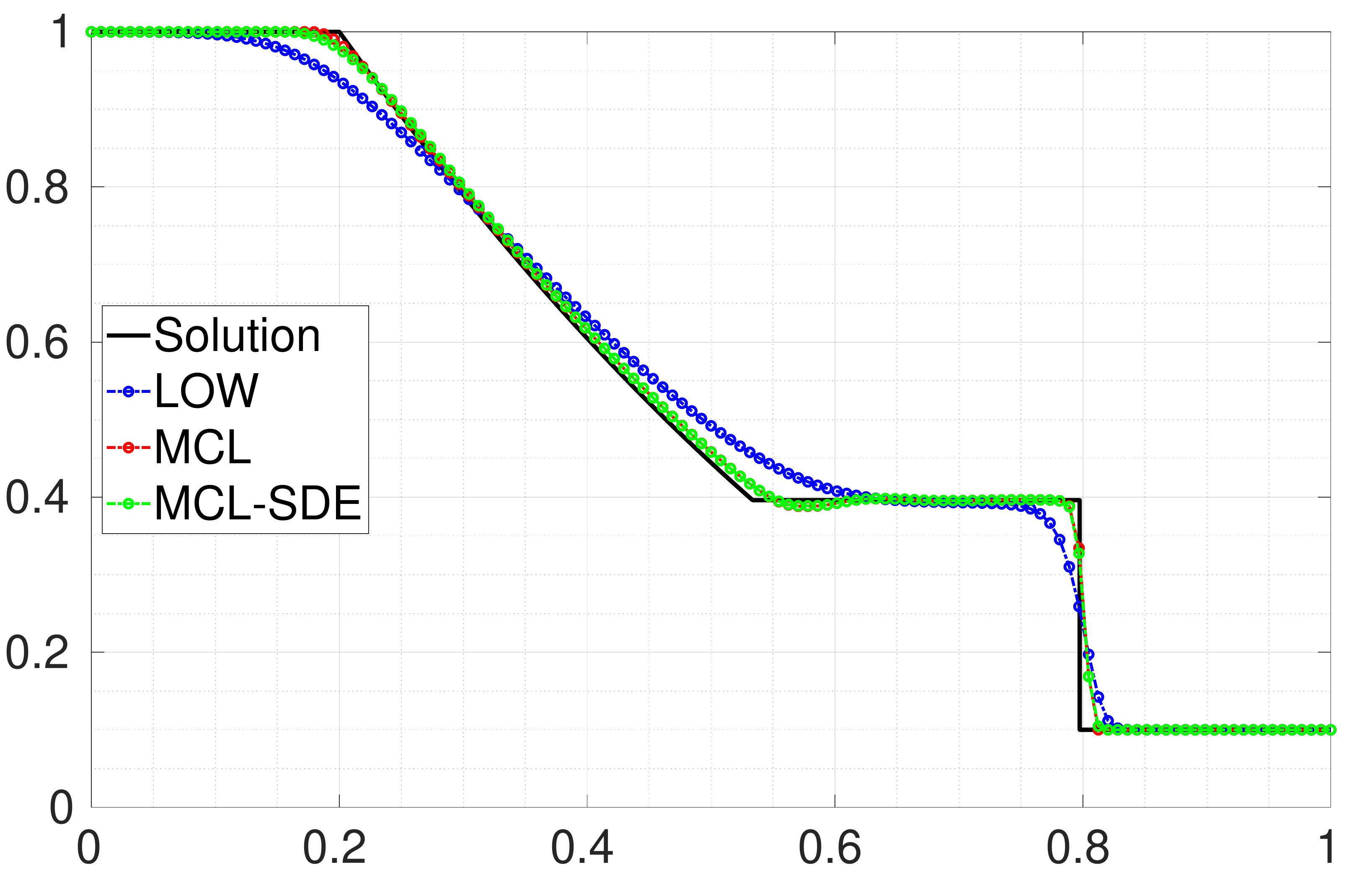}
\end{subfigure}
\begin{subfigure}[b]{0.32\textwidth}
\caption{Discharge}
\includegraphics[width=\textwidth]{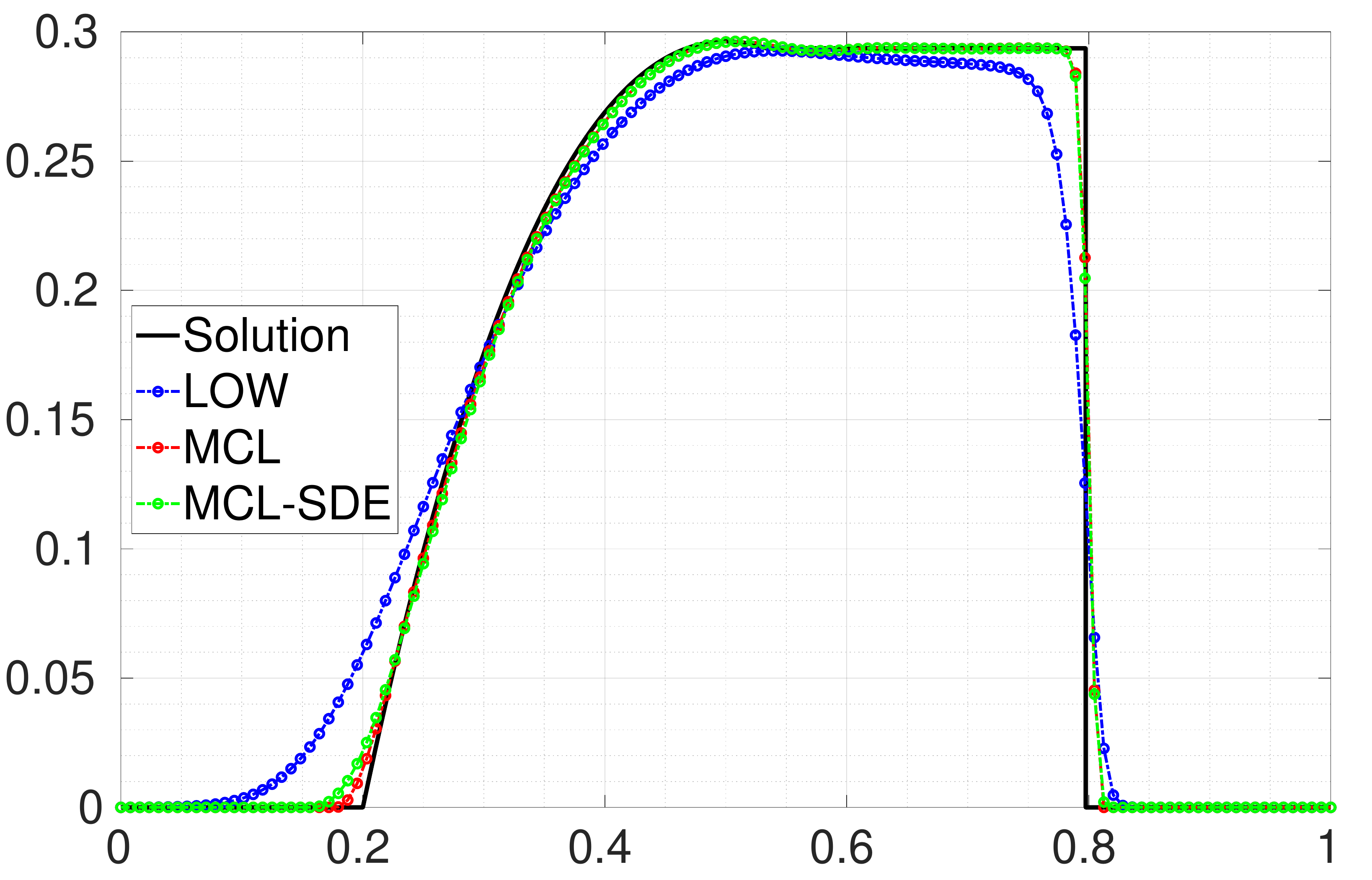}
\end{subfigure}
\begin{subfigure}[b]{0.32\textwidth}
\caption{Velocity}
\includegraphics[width=\textwidth]{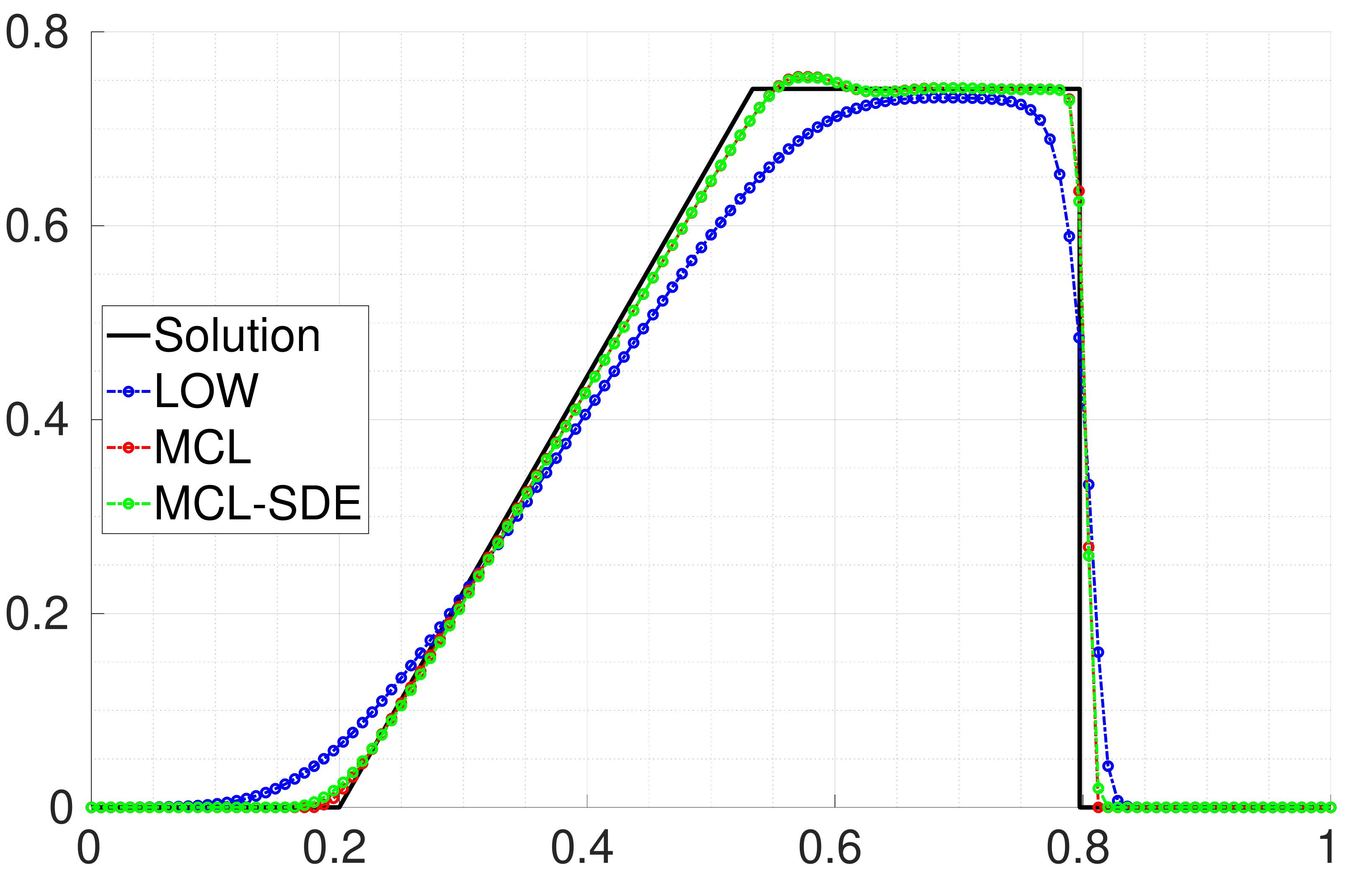}
\end{subfigure}
\caption{Wet dam break for the shallow water equations. Approximations at $T=0.3$ obtained with adaptive SSP2 RK time stepping and $\nu=0.5$ on a uniform mesh consisting of 128 elements.}\label{fig:3-wet}
\end{figure}

In \cref{fig:3-wet}, we display the approximations obtained on a uniform mesh consisting of 128 elements.
The low order profiles are significantly more diffusive than they tend to be in approximations to some scalar problems.
Similarly to shock tube examples for the Euler equations, there are some non-IDP-violating under- and overshoots in the flux-corrected solutions on the right of the rarefaction wave \cite{guermond2016}.

\begin{figure}[ht!]
\centering
\begin{subfigure}[b]{0.32\textwidth}
\caption{SSP3 RK with $\nu=0.67$}\label{3-ssp3}
\includegraphics[width=\textwidth]{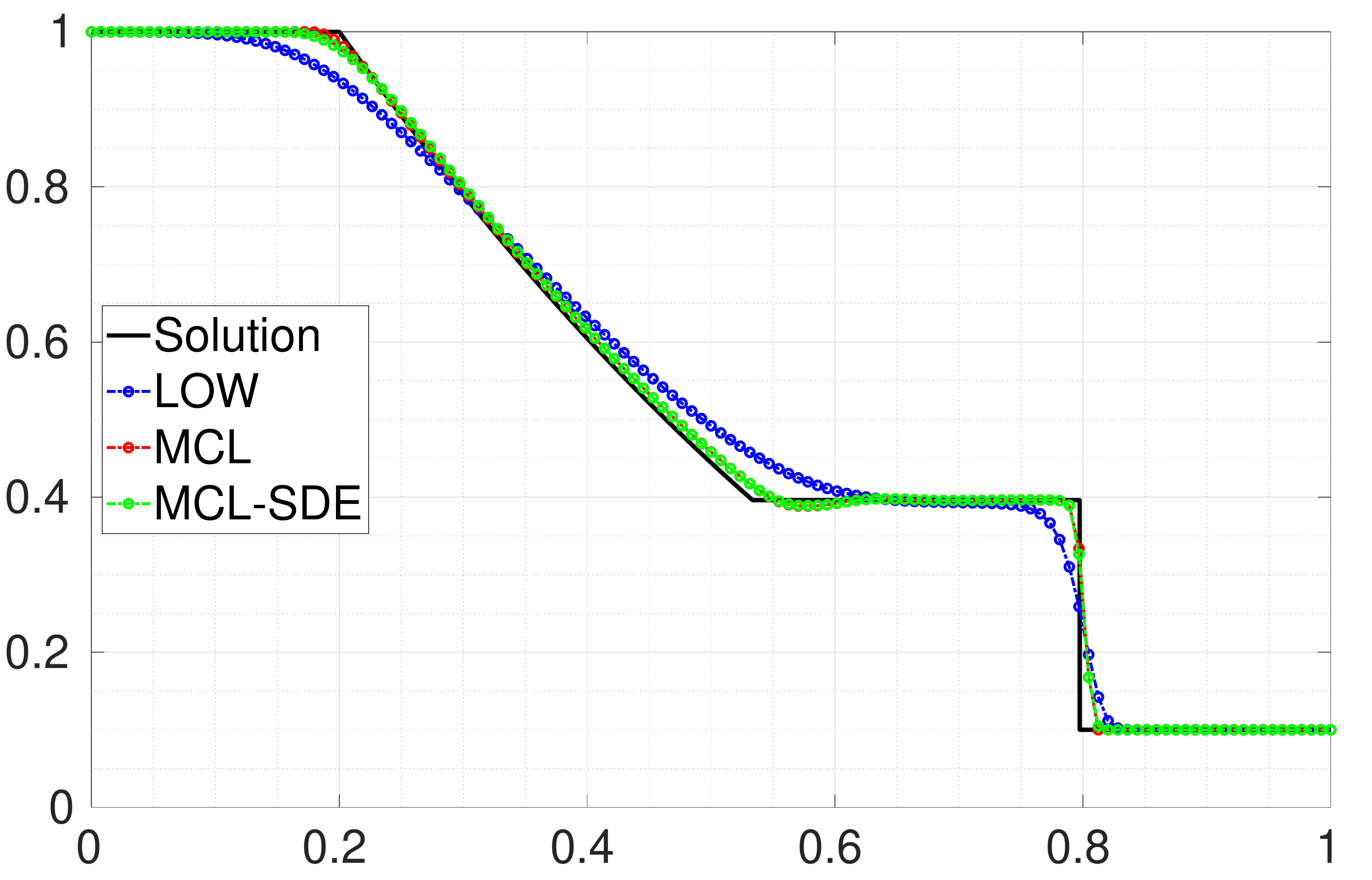}
\end{subfigure}
\begin{subfigure}[b]{0.32\textwidth}
\caption{SSP1 RK, $\nu=1$ with \eqref{eq:dij}}\label{3-ssp1}
\includegraphics[width=\textwidth]{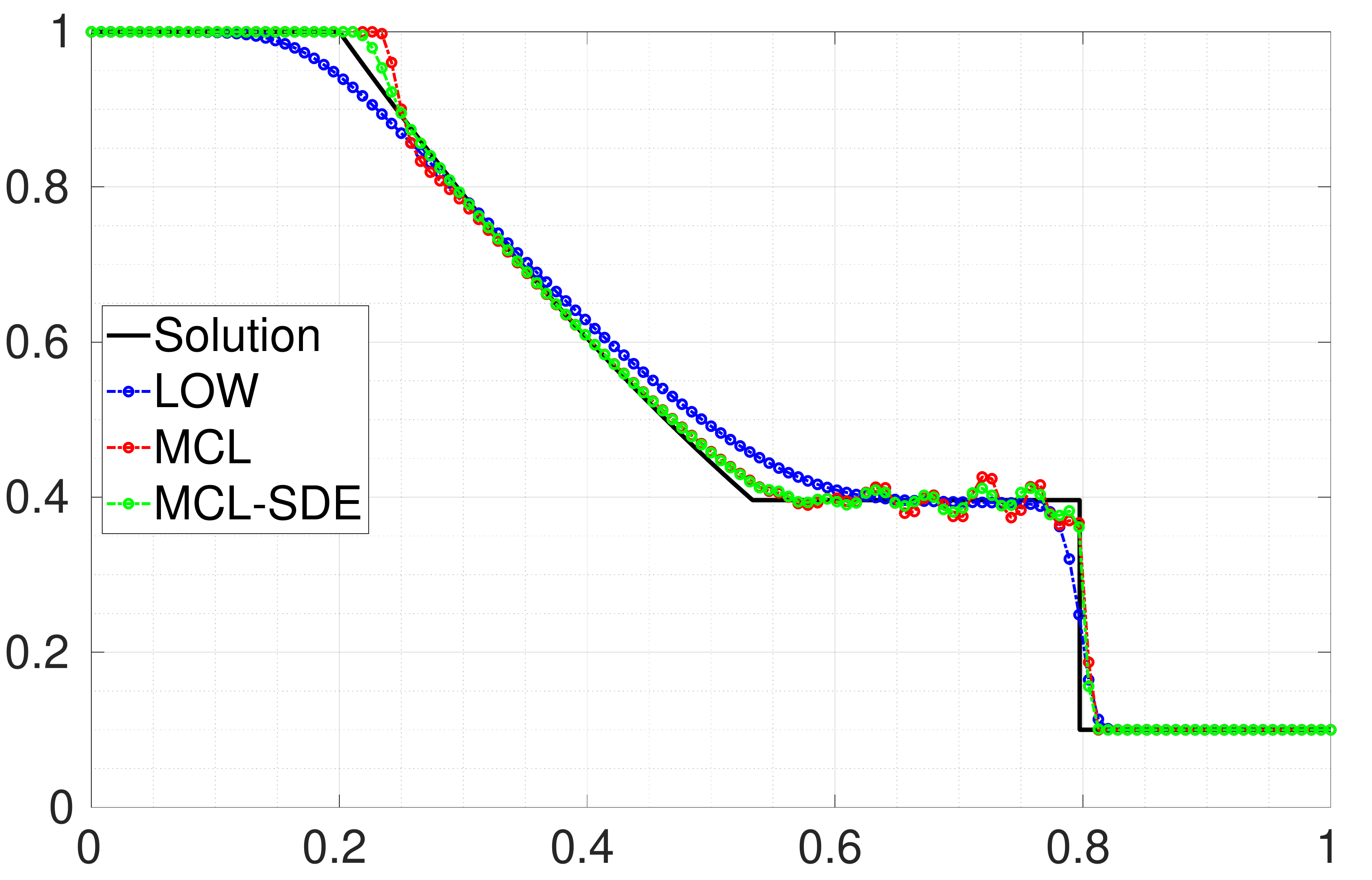}
\end{subfigure}
\begin{subfigure}[b]{0.32\textwidth}
\caption{SSP1 RK, $\nu=1$ with GMS}\label{3-ssp1-GMS}
\includegraphics[width=\textwidth]{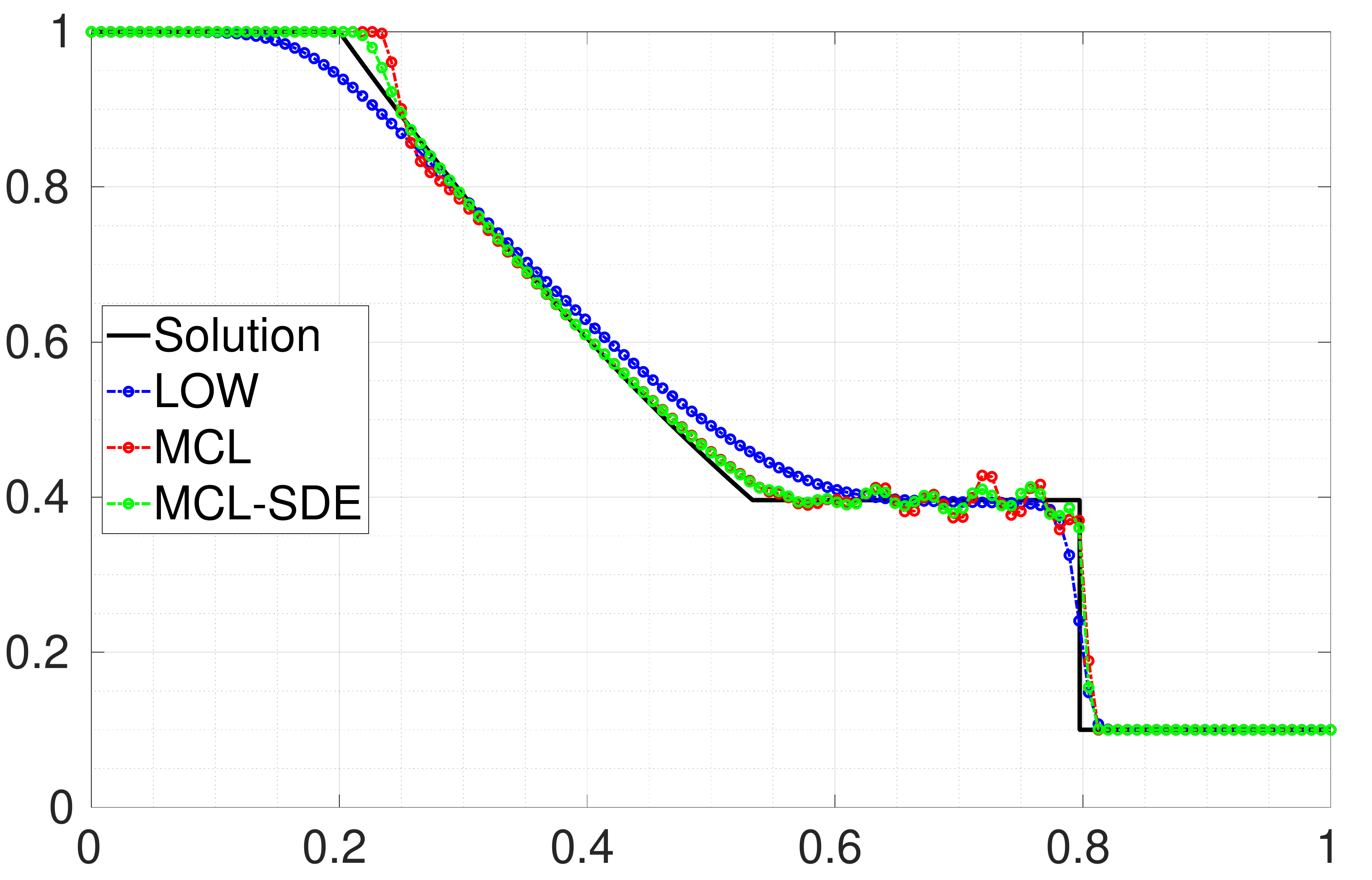}
\end{subfigure}
\caption{Wet dam break for the shallow water equations. Approximations at $T=0.3$ obtained with adaptive SSP RK time stepping on a uniform mesh consisting of 128 elements.}
\label{fig:3-wet-variations}
\end{figure}

Let us briefly study a few variations of time stepping schemes and artificial viscosities.
For completeness we run this problem with the SSP3 RK (Shu--Osher) method and leave the rest of the setup unchanged.
The results for the water heights in \cref{3-ssp3} are indistinguishable from the ones in \cref{fig:3-ssp2}.
In this example, setting $\nu=0.67$ was sufficient for the CFL condition \eqref{eq-cfl} to be satisfied at all times while values $\nu\ge 0.68$ required some repetitions of individual RK stages.

The theory presented in \cref{sec:conslaw} suggests that we may also employ forward Euler (SSP1 RK) time stepping and may even set the CFL parameter $\nu$ to one.
Indeed, the resulting approximations do not violate the IDP property, and one may assume the results to be reliable.
We investigate the validity of this assumption, first by using our nodal approximation \eqref{eq:dij} to the wave speeds, and, alternatively, the guaranteed maximum wave speed (GMS) proposed in \cite[Prop.~3.7]{azerad2017}, \cite[Sec.~4]{guermond2018a}.
In either case, oscillations and incorrect approximations in the left part of the rarefaction waves are visible in \cref{3-ssp1,3-ssp1-GMS}.
The fact that even the use of the GMS wave speed is not sufficient to prevent these nonphysical effects, implies that they are not a result of an incorrect wave speed approximation.
A reduction of the CFL parameter $\nu$ masks this issue in the sense that the amplitude of the oscillations becomes smaller.
The spurious approximations observed in this study originate from the combination of forward Euler time stepping with continuous finite element methods \cite[Sec.~4]{kuzmin2012a}.
Schemes based on discontinuous approximation spaces on the other hand can be safely employed in combination with forward Euler time stepping.
The fact that the low order method remains stable can be attributed to its equivalence to the vertex-centered finite volume scheme of local Lax--Friedrichs type \cite{selmin1993}.
For the shallow water equations we tested the GMS wave speed \cite[Prop.~3.7]{azerad2017}, \cite[Sec.~4]{guermond2018a} for multiple examples.
Although the low order method is derived based on assumptions that encourage the use of GMS instead of our approximation \eqref{eq:dij}, we encountered no example in which the use of GMS is actually necessary.
This observation was recently confirmed by Wu \etal \cite[Thm.~3.1]{wu-arxiv} who show that \eqref{eq:dij} preserves the IDP property for the SWE.

Using SSP$p$ RK time stepping with $p\in \{2,3\}$, we observe satisfactory agreement of MCL and MCL-SDE profiles with the exact solutions, not only for the conserved unknowns but also for the velocity.
The situation may be different if the test problem features dry or nearly dry states, which is why we consider such an example next.

\subsubsection{Dry dam break over flat topography}

Let us now set $h_R$ to zero, the end time to $T=0.15$, and leave the rest of the setup from the previous example unchanged.
Here the exact solution does not feature a shock wave, only a rarefaction wave is produced as a result of the dam break \cite[Sec.~4.1.2]{delestre-arxiv}.

\begin{figure}[ht!]
\centering
\begin{subfigure}[b]{0.32\textwidth}
\caption{Friction-based, water level}
\includegraphics[width=\textwidth]{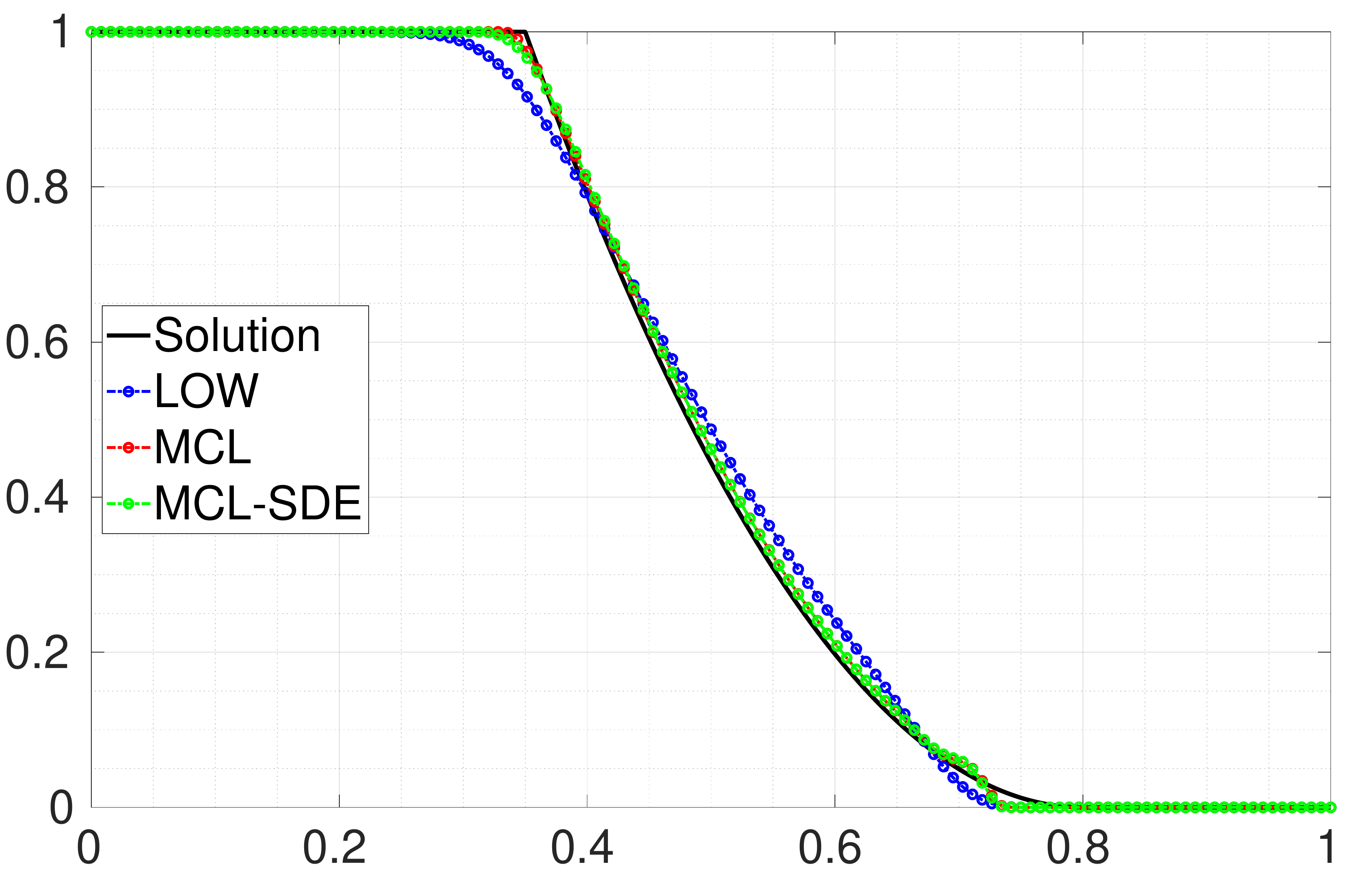}
\end{subfigure}
\begin{subfigure}[b]{0.32\textwidth}
\caption{Friction-based, discharge}
\includegraphics[width=\textwidth]{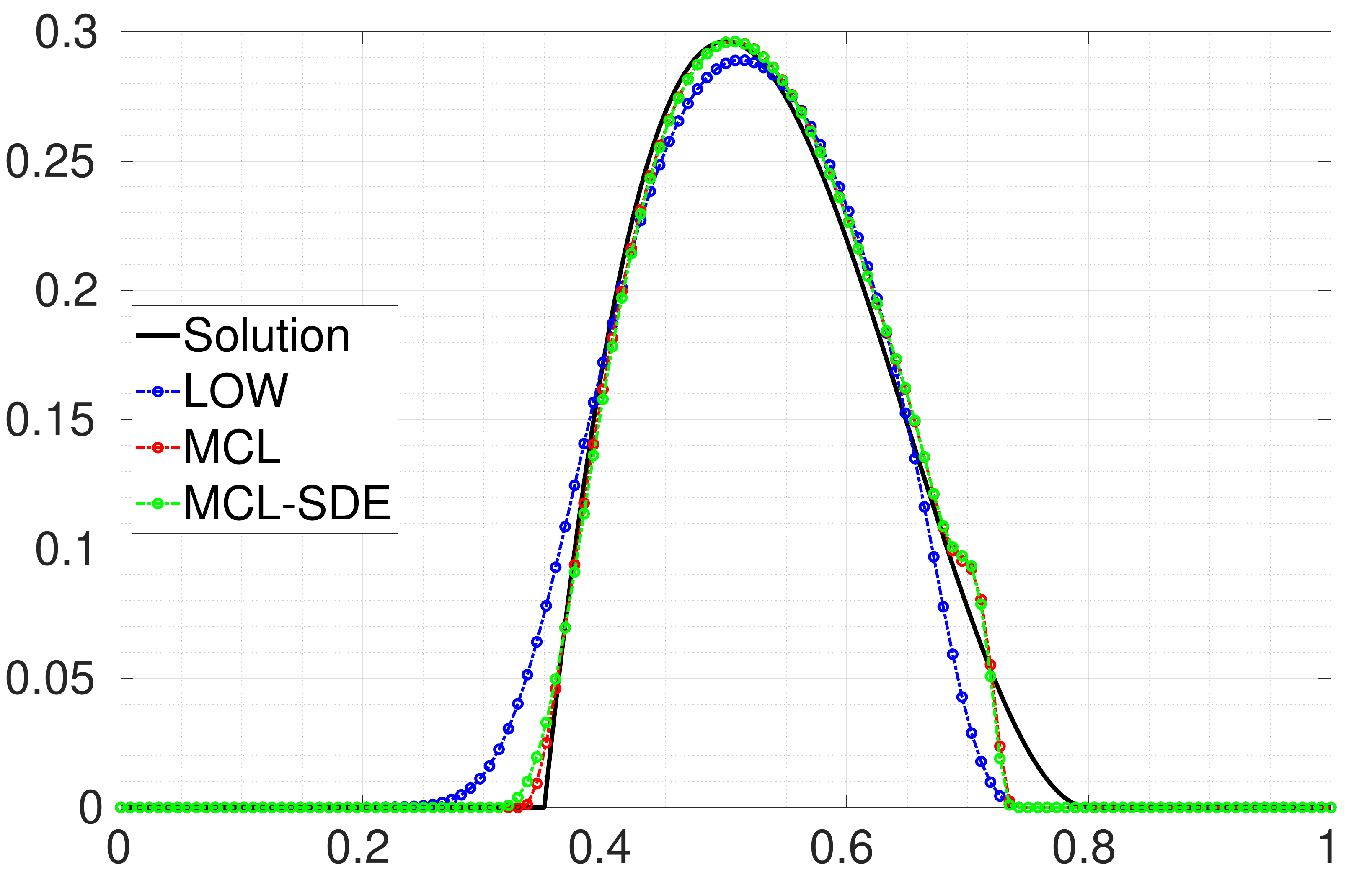}
\end{subfigure}
\begin{subfigure}[b]{0.32\textwidth}
\caption{Friction-based, velocity}
\includegraphics[width=\textwidth]{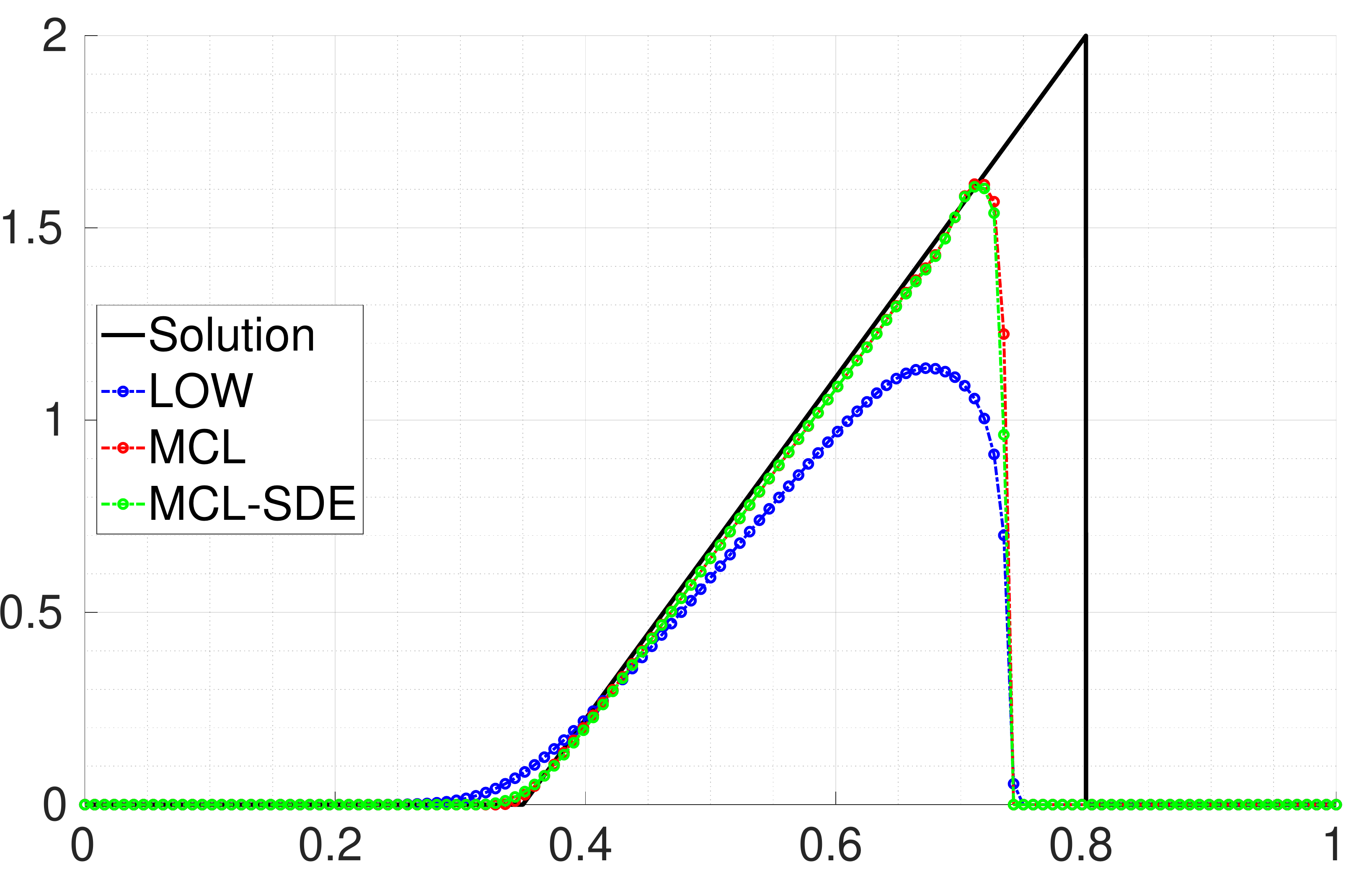}
\end{subfigure}
\begin{subfigure}[b]{0.32\textwidth}
\caption{Entropy-based, water level}
\includegraphics[width=\textwidth]{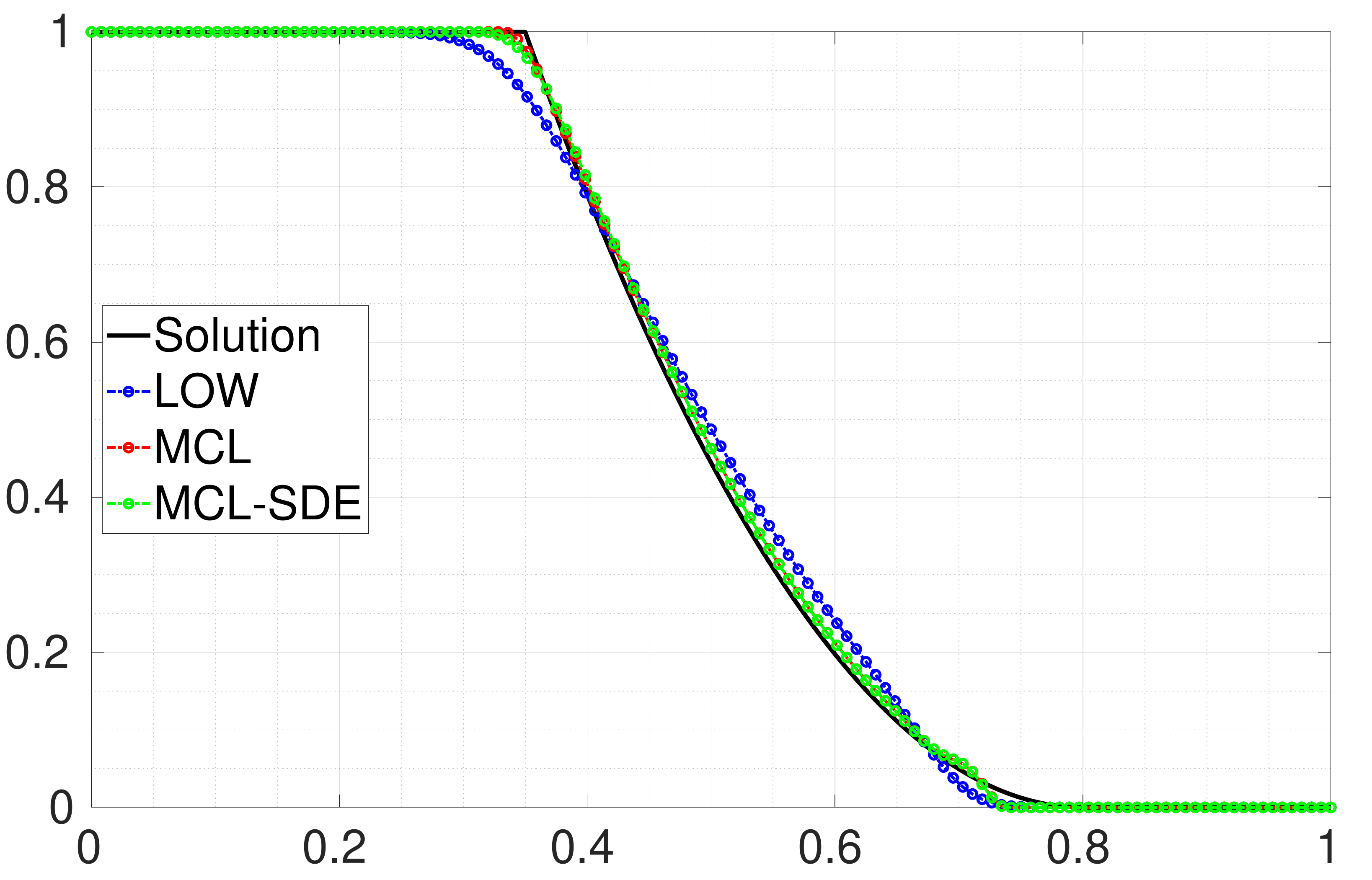}
\end{subfigure}
\begin{subfigure}[b]{0.32\textwidth}
\caption{Entropy-based, discharge}
\includegraphics[width=\textwidth]{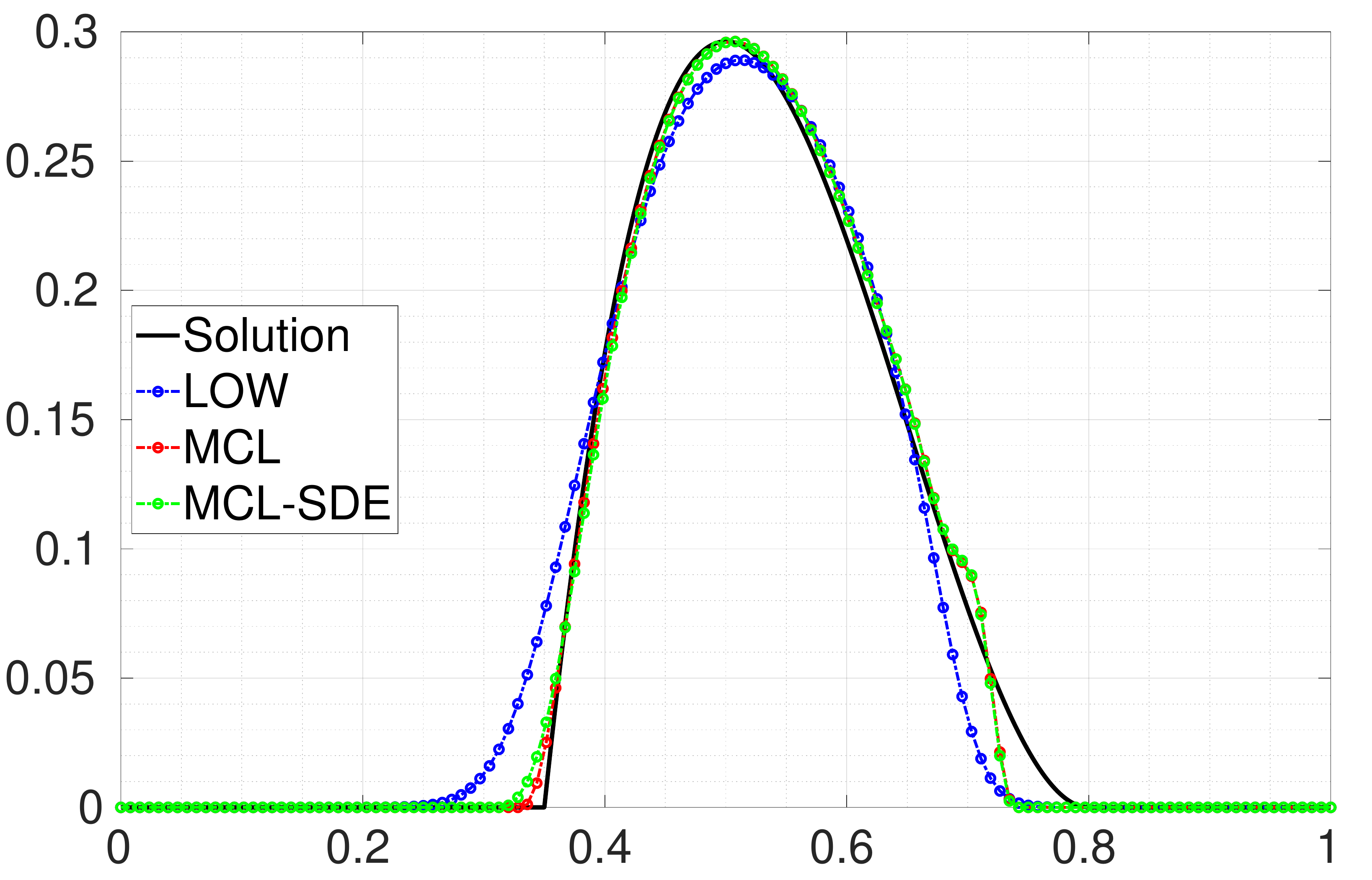}
\end{subfigure}
\begin{subfigure}[b]{0.32\textwidth}
\caption{Entropy-based, velocity}
\includegraphics[width=\textwidth]{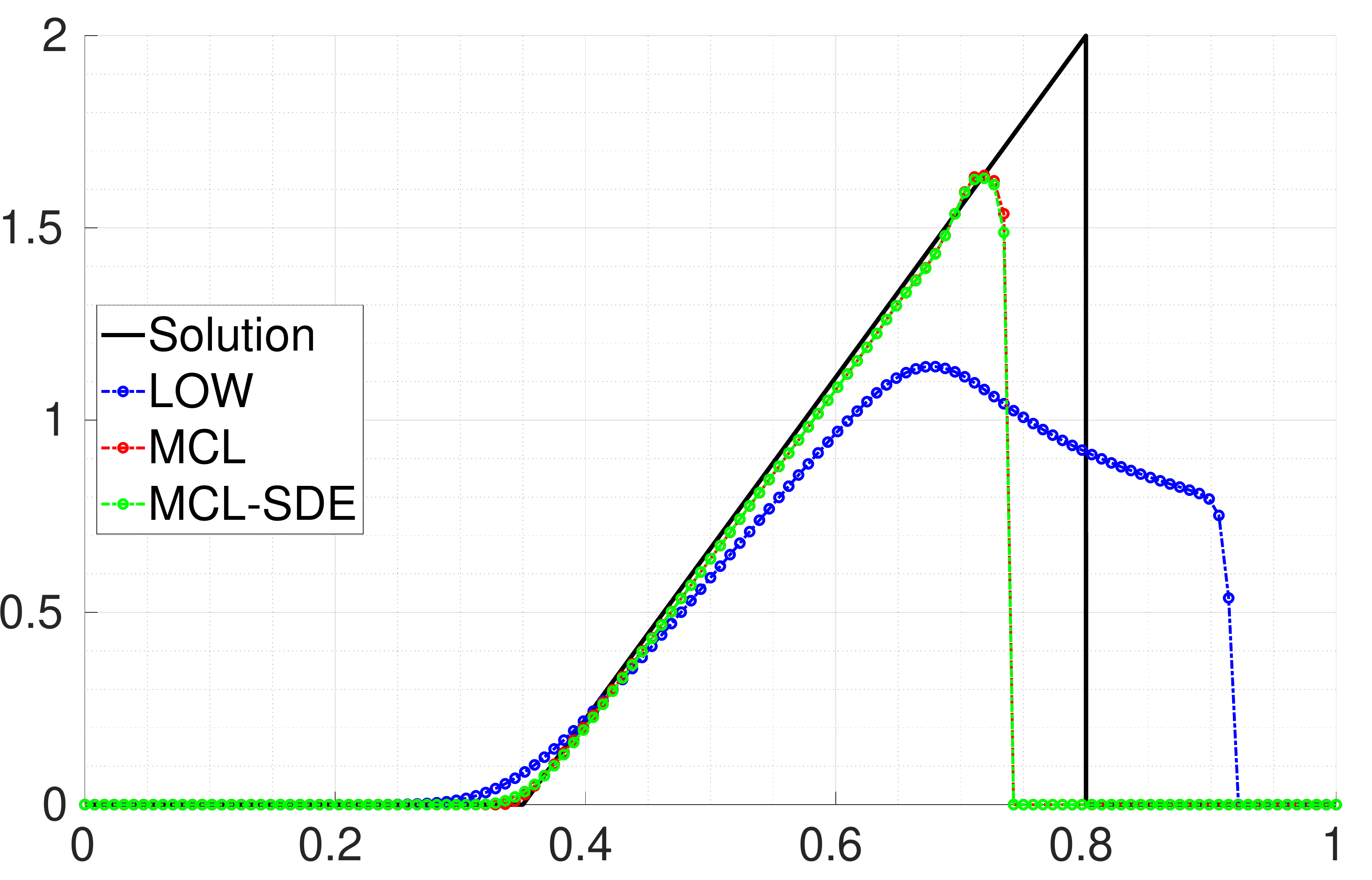}
\end{subfigure}
\caption{Dry dam break for the shallow water equations with new wetting and drying strategies.
Approximations at $T=0.15$ obtained with adaptive SSP2 RK time stepping and $\nu=0.5$ on a uniform mesh consisting of 128 elements.}
\label{fig:3-dam-dry}
\end{figure}

\begin{figure}[ht!]
\centering
\begin{subfigure}[b]{0.32\textwidth}
\caption{\cite{azerad2017}-based, water level}
\includegraphics[width=\textwidth]{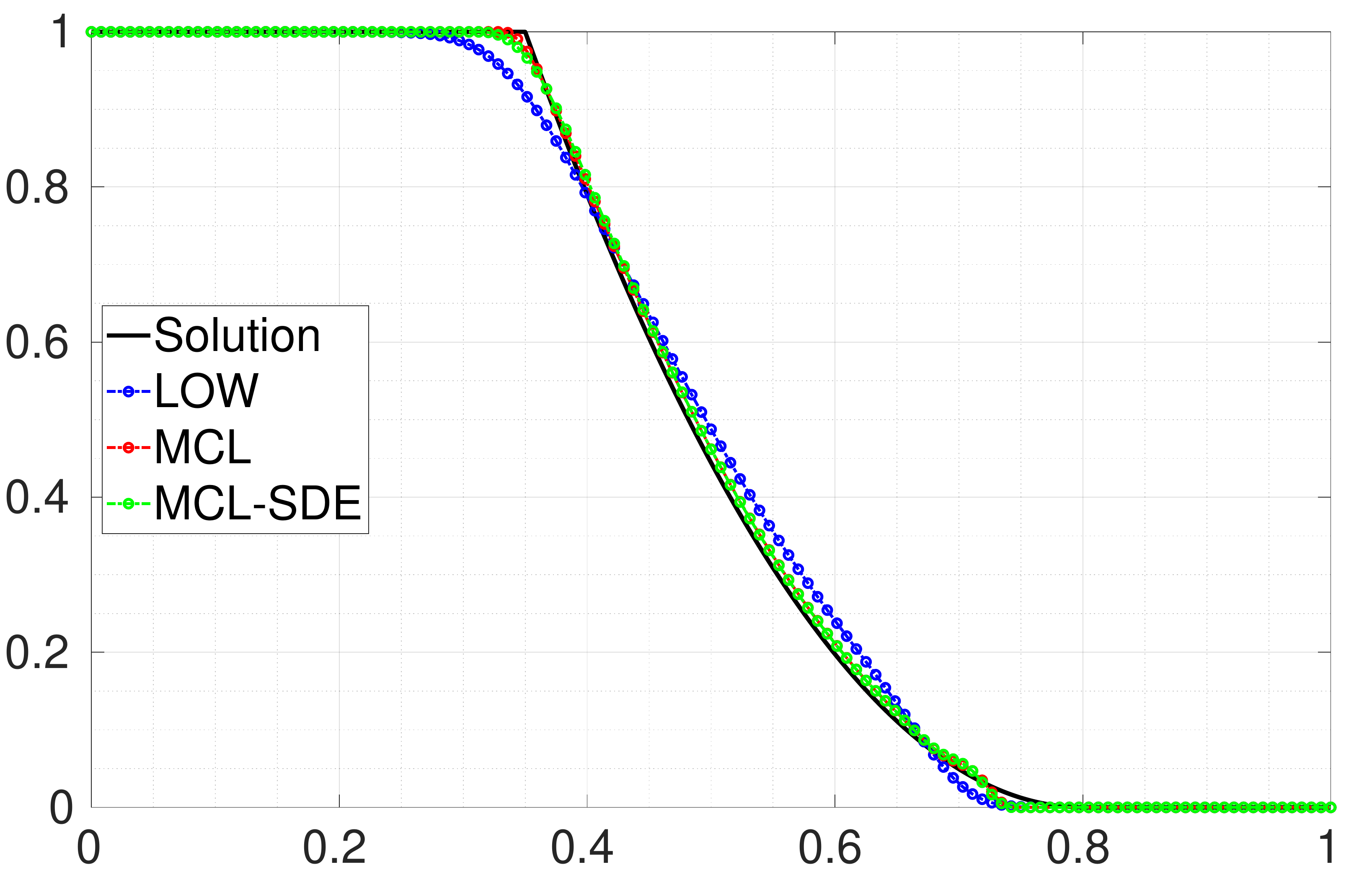}
\end{subfigure}
\begin{subfigure}[b]{0.32\textwidth}
\caption{\cite{azerad2017}-based, discharge}
\includegraphics[width=\textwidth]{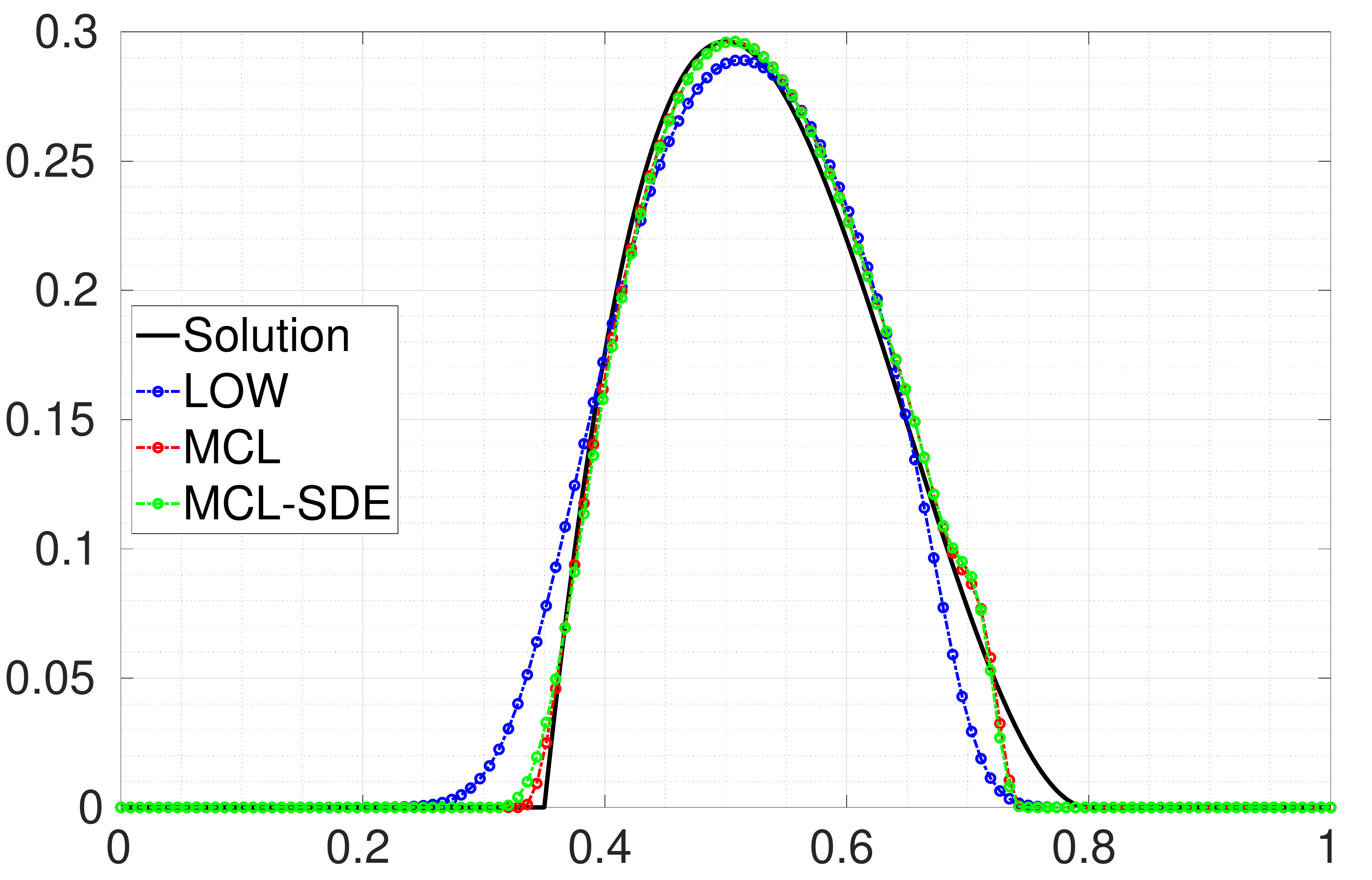}
\end{subfigure}
\begin{subfigure}[b]{0.32\textwidth}
\caption{\cite{azerad2017}-based, velocity}
\includegraphics[width=\textwidth]{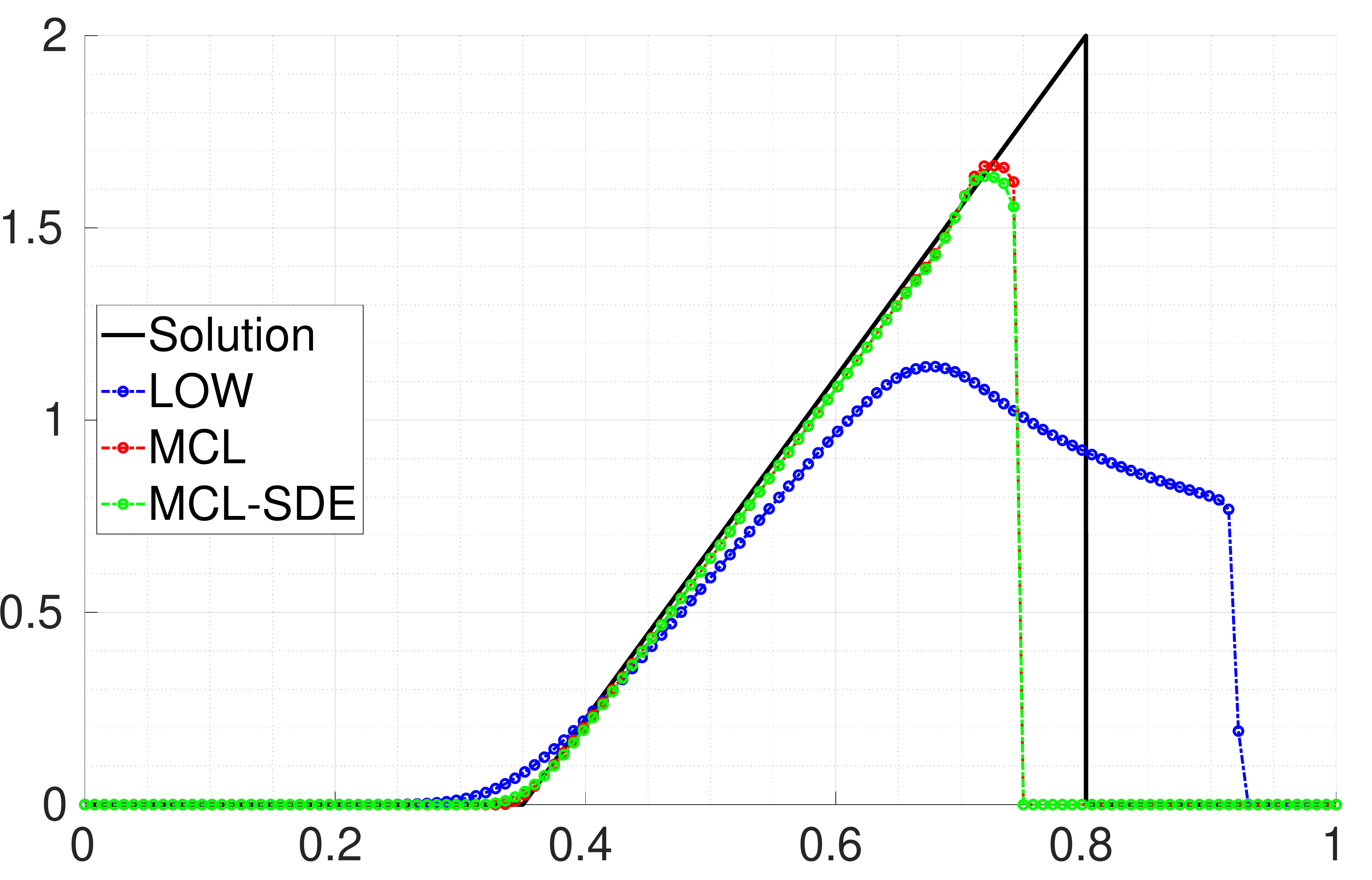}
\end{subfigure}
\begin{subfigure}[b]{0.32\textwidth}
\caption{\cite{kurganov2007a}-based, water level}
\includegraphics[width=\textwidth]{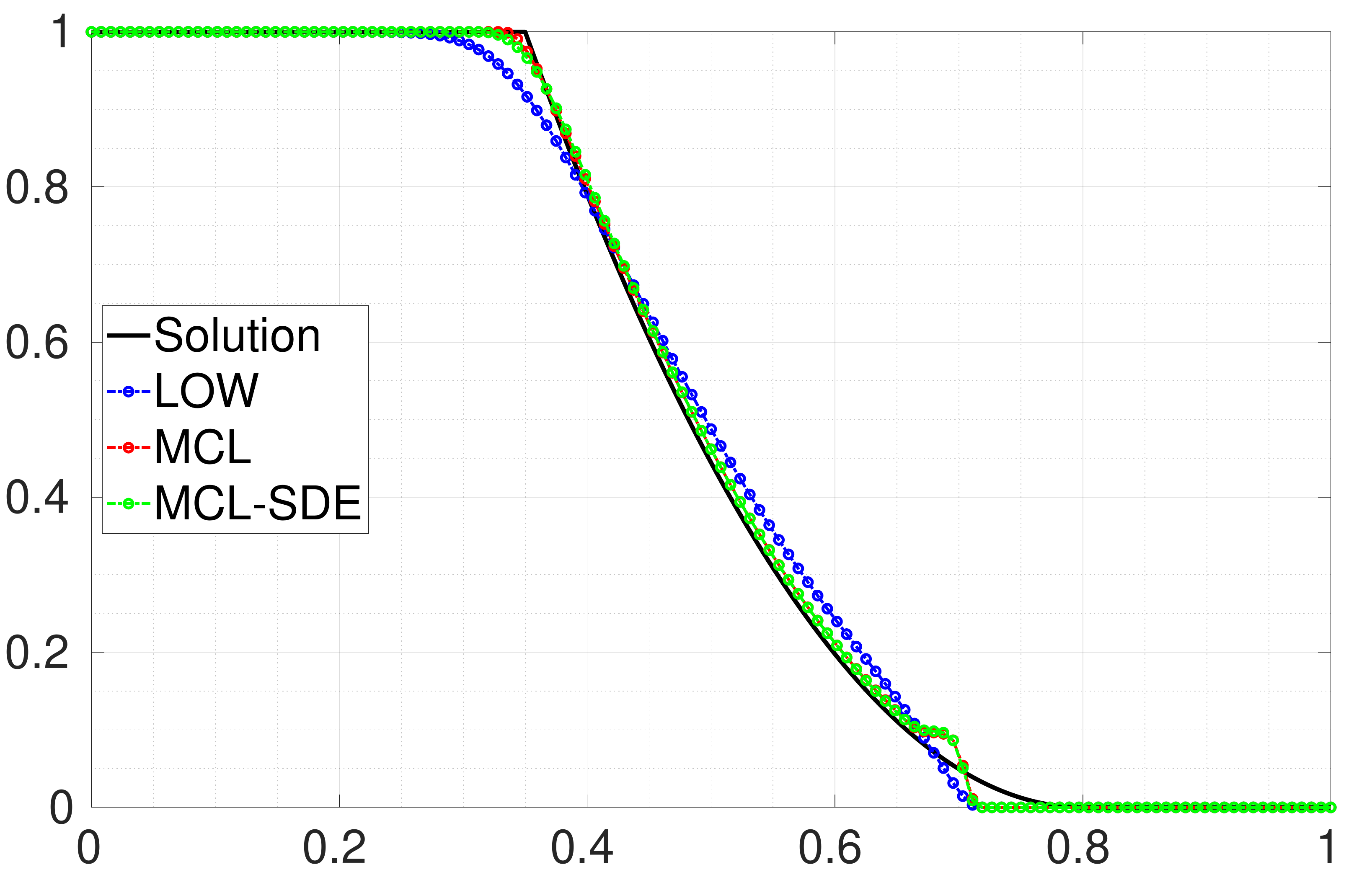}
\end{subfigure}
\begin{subfigure}[b]{0.32\textwidth}
\caption{\cite{kurganov2007a}-based, discharge}
\includegraphics[width=\textwidth]{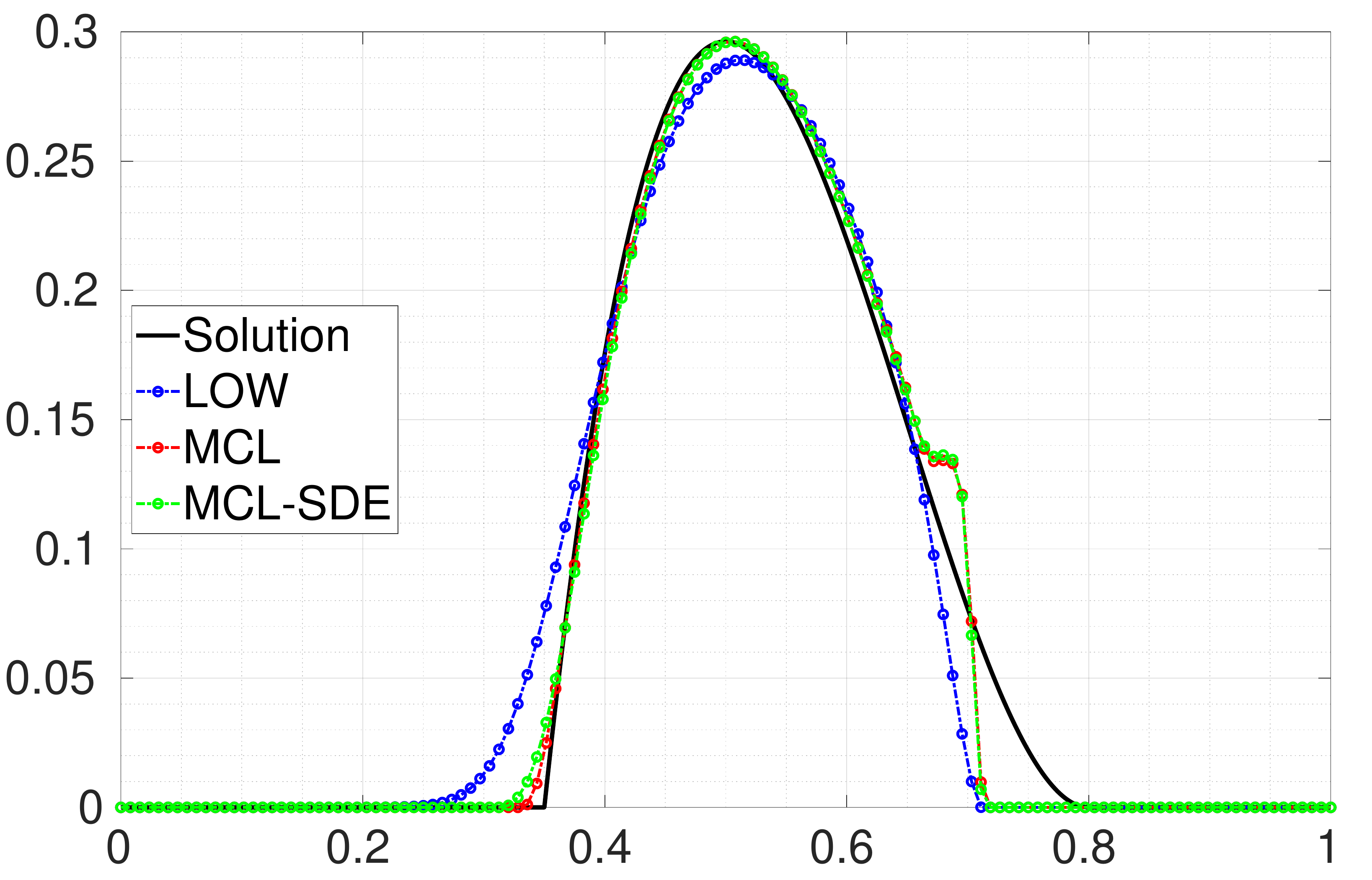}
\end{subfigure}
\begin{subfigure}[b]{0.32\textwidth}
\caption{\cite{kurganov2007a}-based, velocity}
\includegraphics[width=\textwidth]{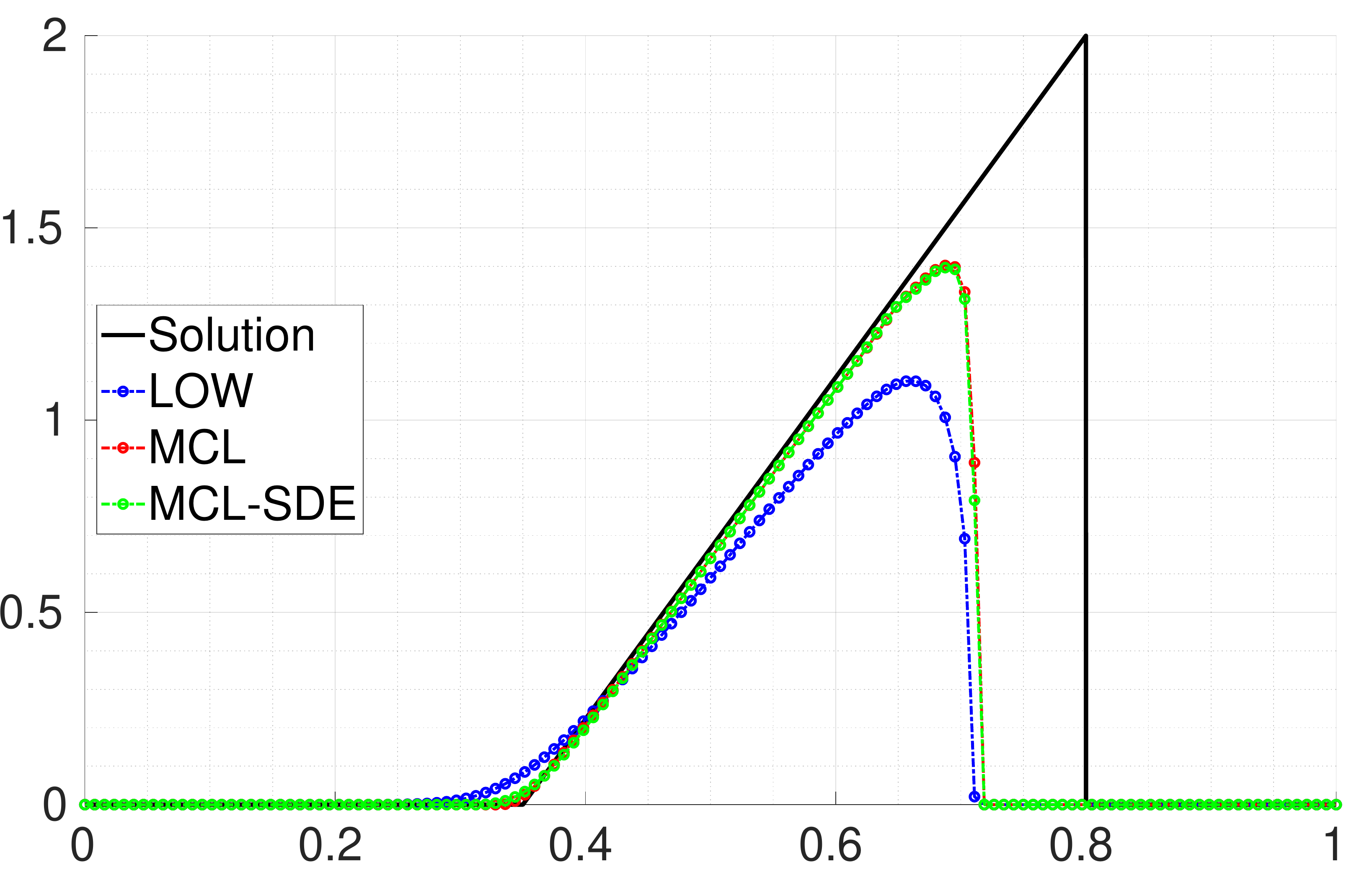}
\end{subfigure}
\begin{subfigure}[b]{0.32\textwidth}
\caption{\cite{ricchiuto2009}-based, water level}
\includegraphics[width=\textwidth]{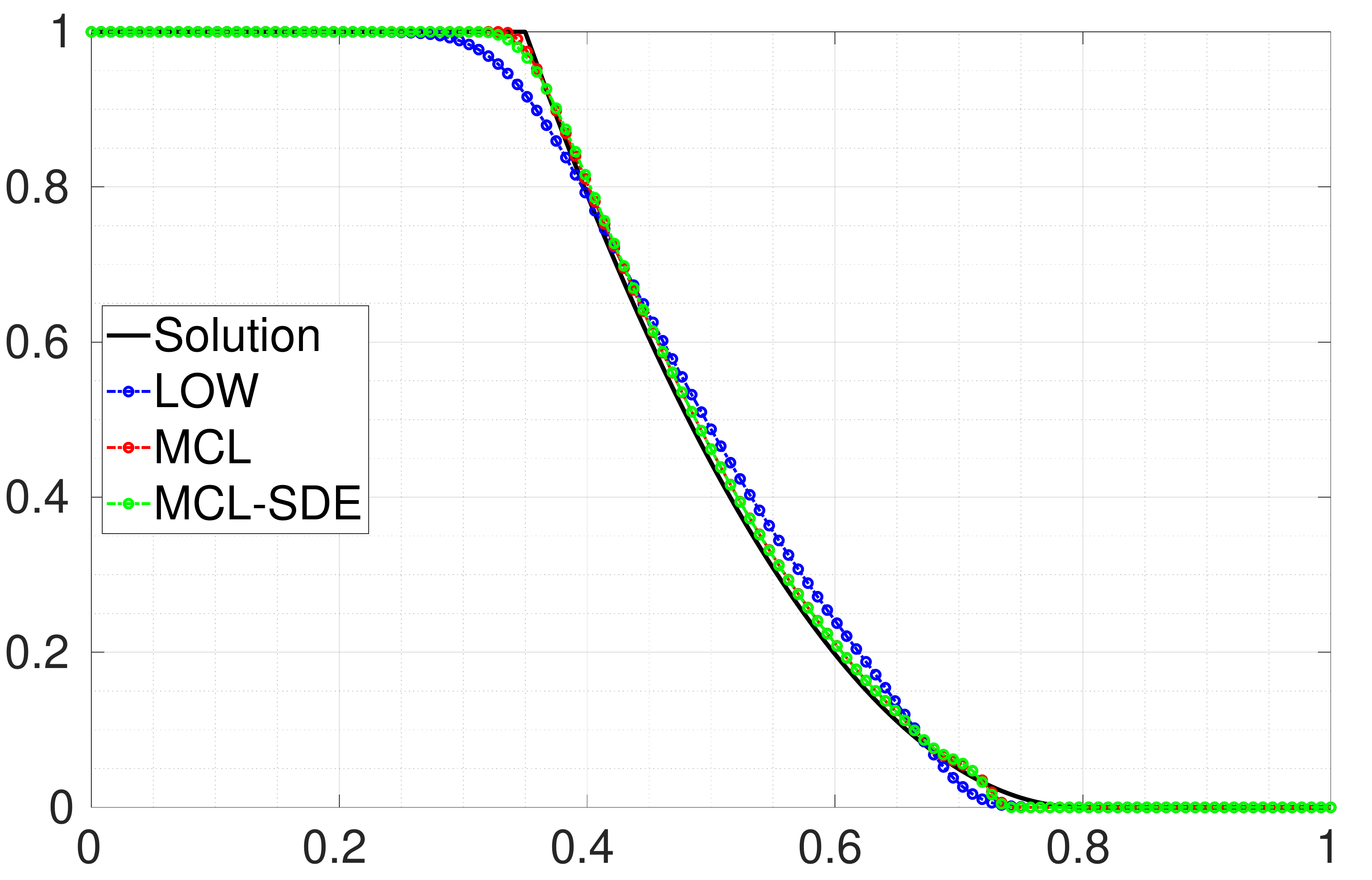}
\end{subfigure}
\begin{subfigure}[b]{0.32\textwidth}
\caption{\cite{ricchiuto2009}-based, discharge}
\includegraphics[width=\textwidth]{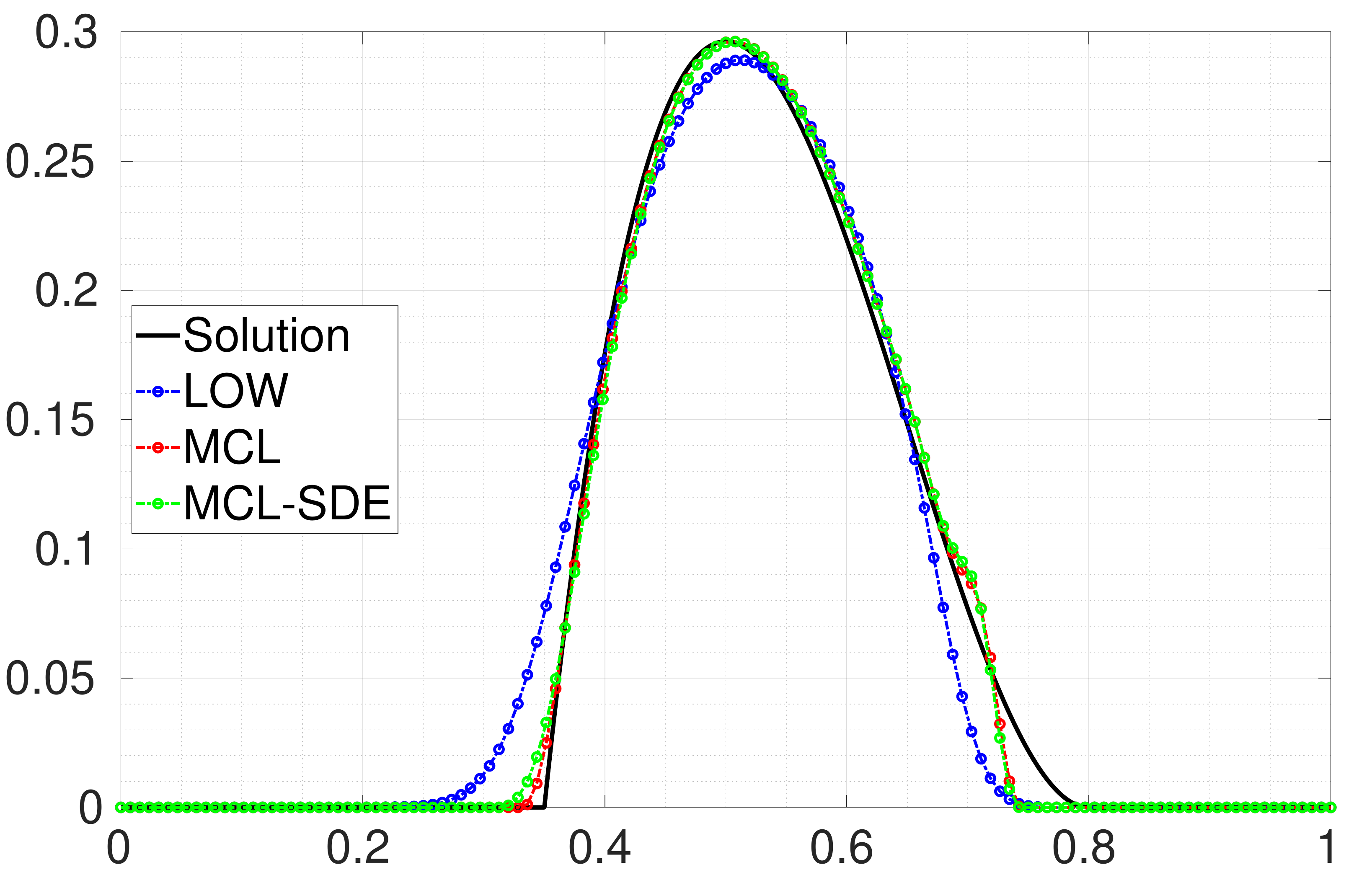}
\end{subfigure}
\begin{subfigure}[b]{0.32\textwidth}
\caption{\cite{ricchiuto2009}-based, velocity}
\includegraphics[width=\textwidth]{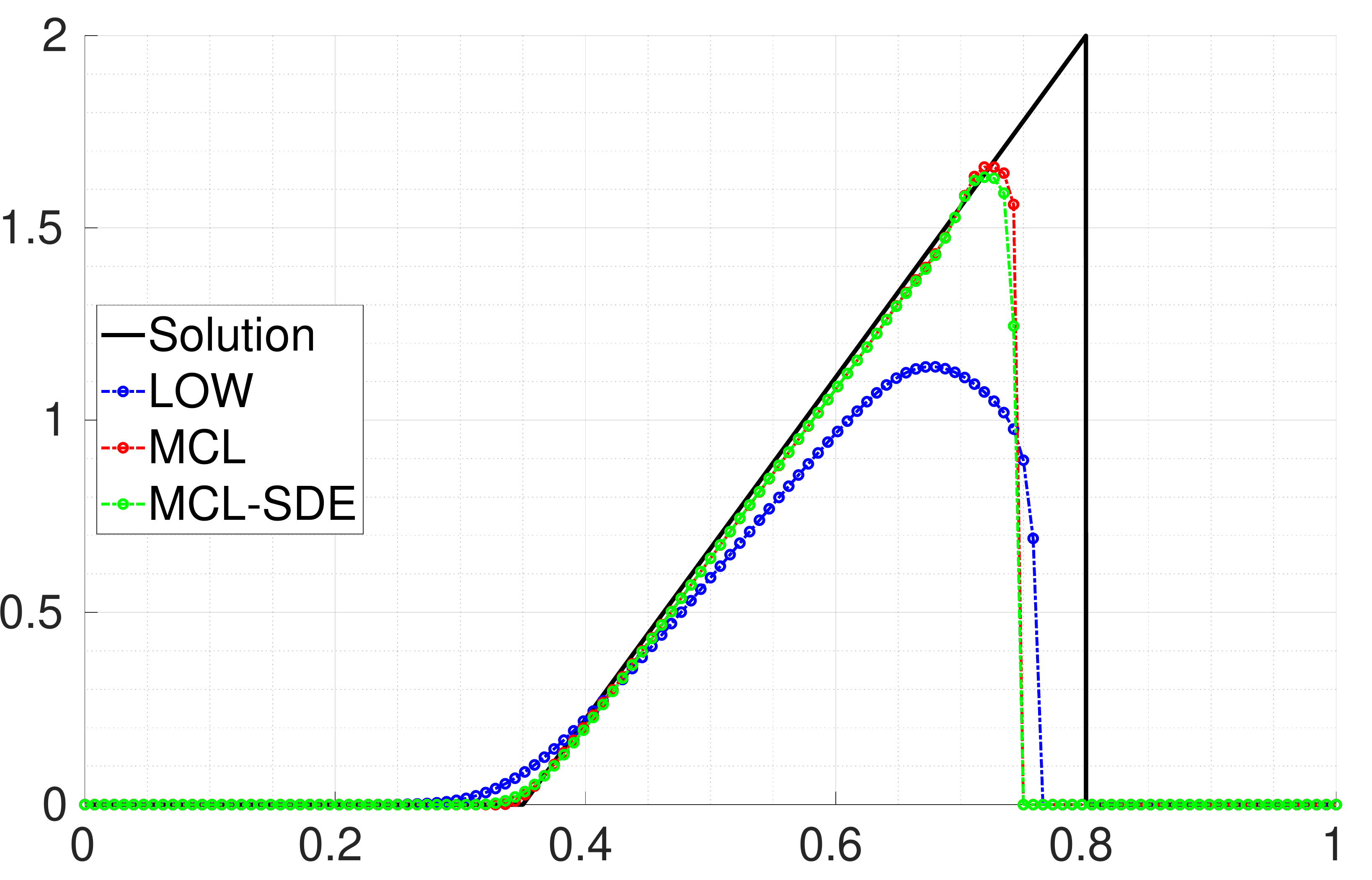}
\end{subfigure}
\caption{Dry dam break for the shallow water equations with wetting and drying strategies from the literature. Approximations at $T=0.15$ obtained with adaptive SSP2 RK time stepping and $\nu=0.5$ on a uniform mesh consisting of 128 elements.}\label{fig:3-dam-dry2}
\end{figure}

This example already requires some form of treatment to correctly capture the wet-dry transition.
If no such approach has been implemented, one can simply run this example by setting $h_R$ to a very small value, for instance $10^{-12}$.
Instead of following this approach, we compare the results obtained with our new friction- and entropy-based wetting and drying algorithms in \cref{fig:3-dam-dry}.
The approximations obtained with the existing schemes of Azerad \etal \cite{azerad2017}, Kurganov and Petrova \cite{kurganov2007a} as well as the one by Ricchiuto and Bollermann \cite{ricchiuto2009} are shown in \cref{fig:3-dam-dry2}.

All wetting and drying approaches produce acceptable numerical solutions for the water height, whereas approximations for the discharge close to the dry region are somewhat underresolved.
Improvements can only be obtained with refined meshes and time steps.
Among the five approaches for wetting and drying, the \cite{kurganov2007a}-based fix \eqref{eq-kp-fix} produces the most pronounced kink in the discharge and a similar artifact is visible in the corresponding water levels.
Moreover, this fix produces the smallest velocities among all considered approaches.
All other wetting and drying algorithms produce satisfactory results for this test problem.
The somewhat significant differences in the velocities, particularly for the low order solution are unsurprising to us because the calculation of $v=(hv)/h$ is quite sensitive to small water heights, which occur in almost dry areas.
Again, refinement is needed to obtain more accurately resolved velocity profiles.

\subsubsection{Wet dam break over a bump}

Next, we study a dam break problem proposed in \cite[Sec.~5.6]{winters2015}.
It involves a nonflat bottom topography.
The spatial domain $\Omega=(0,20)$ is again equipped with reflecting wall boundaries and the gravitational constant is $g=1$.
The bottom topography, and initial conditions read
\begin{align*}
b(x) = \begin{cases}\sin(0.25\pi x) & \text{if } |x-x_0| < 2, \\
0 & \text{otherwise},
\end{cases}\qquad h_0(x) = \begin{cases}
1.6-b(x) & \text{if } x < x_0, \\ 1.05-b(x) & \text{if } x > x_0,
\end{cases}
\end{align*}
where $x_0=10$ and $v_0 \equiv 0$.

\begin{figure}[ht!]
\centering
\begin{subfigure}[b]{0.32\textwidth}
\caption{Water level}
\includegraphics[width=\textwidth]{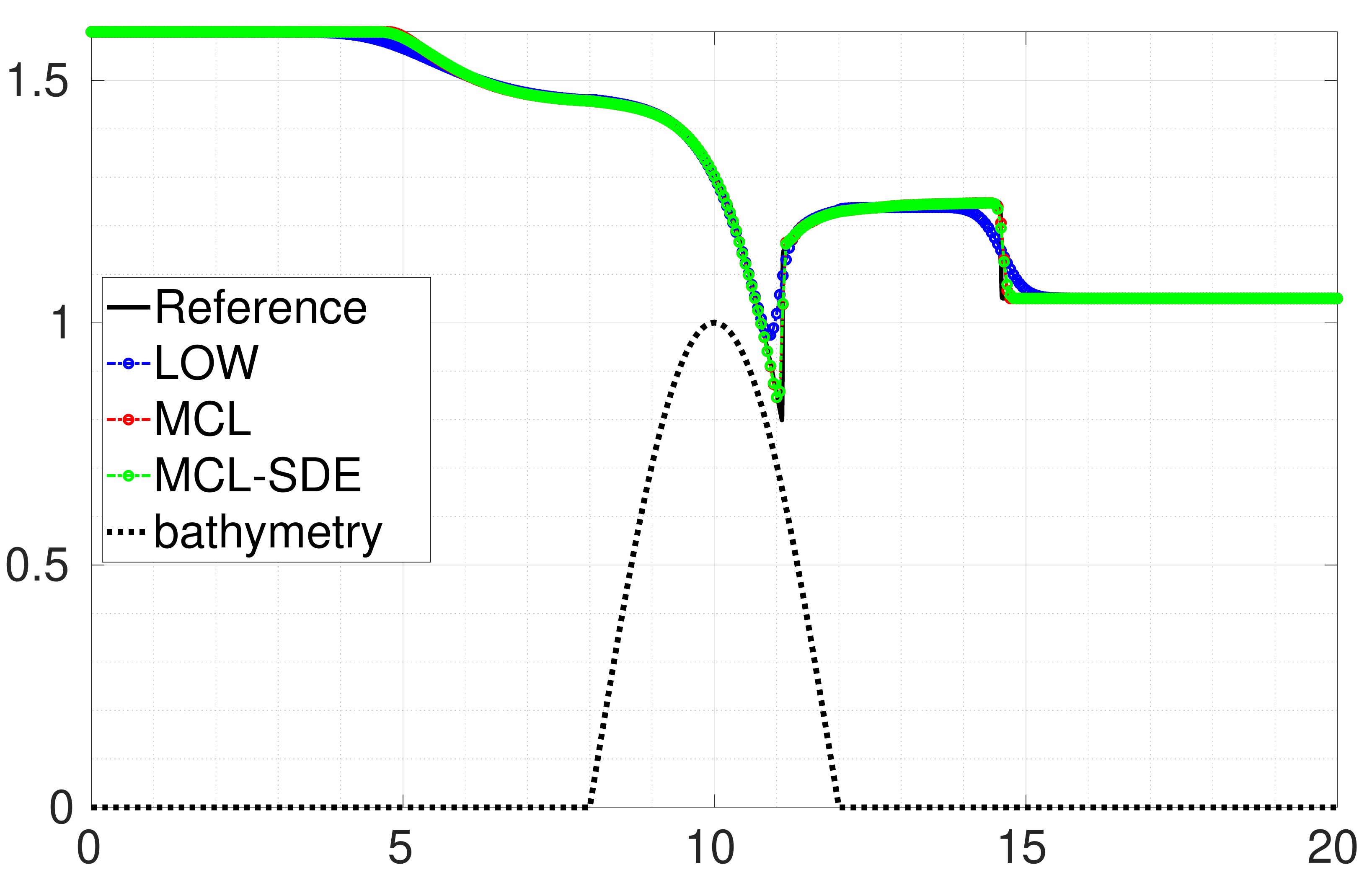}
\end{subfigure}
\begin{subfigure}[b]{0.32\textwidth}
\caption{Discharge}
\includegraphics[width=\textwidth]{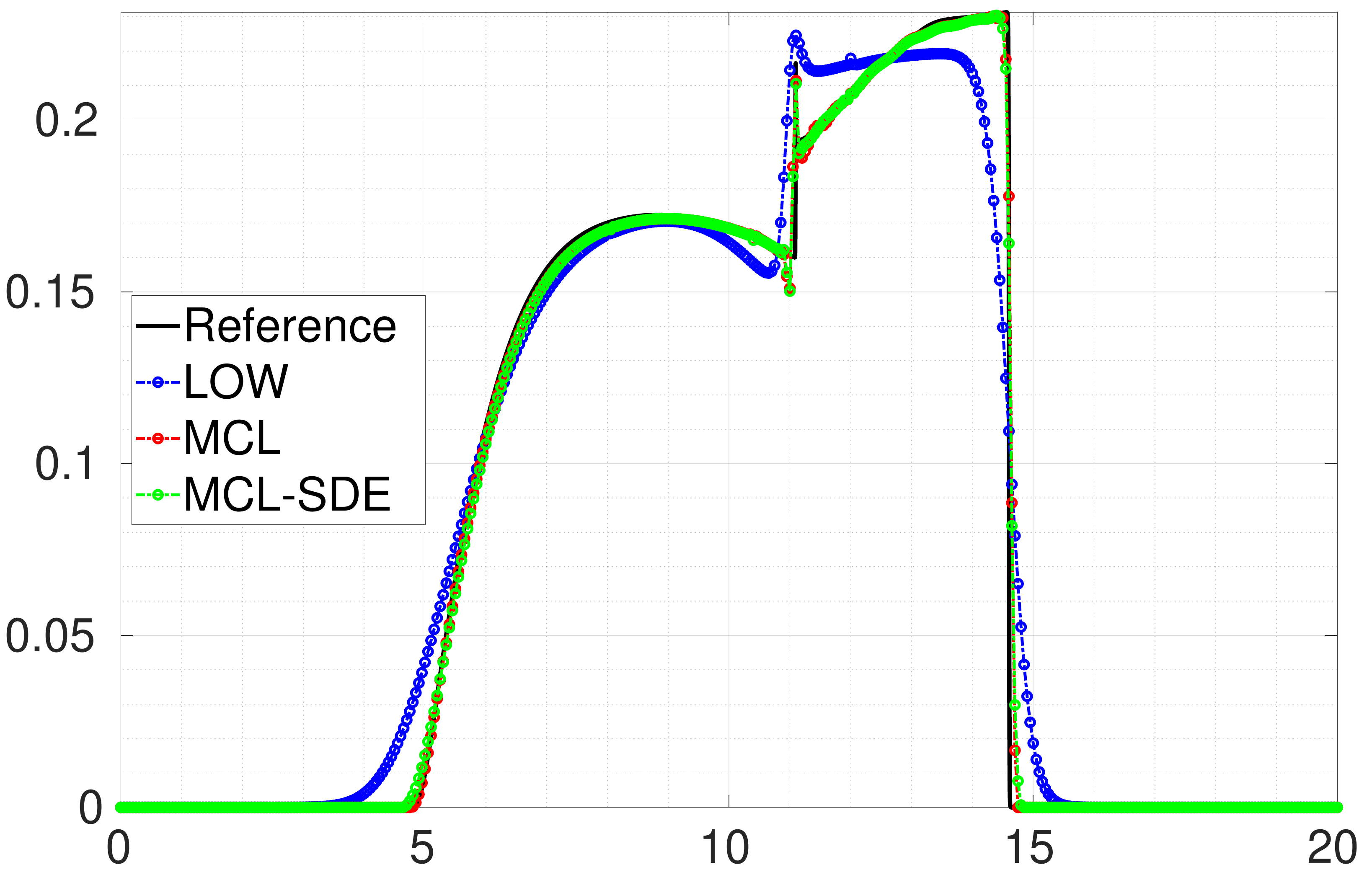}
\end{subfigure}
\begin{subfigure}[b]{0.32\textwidth}
\caption{Velocity}
\includegraphics[width=\textwidth]{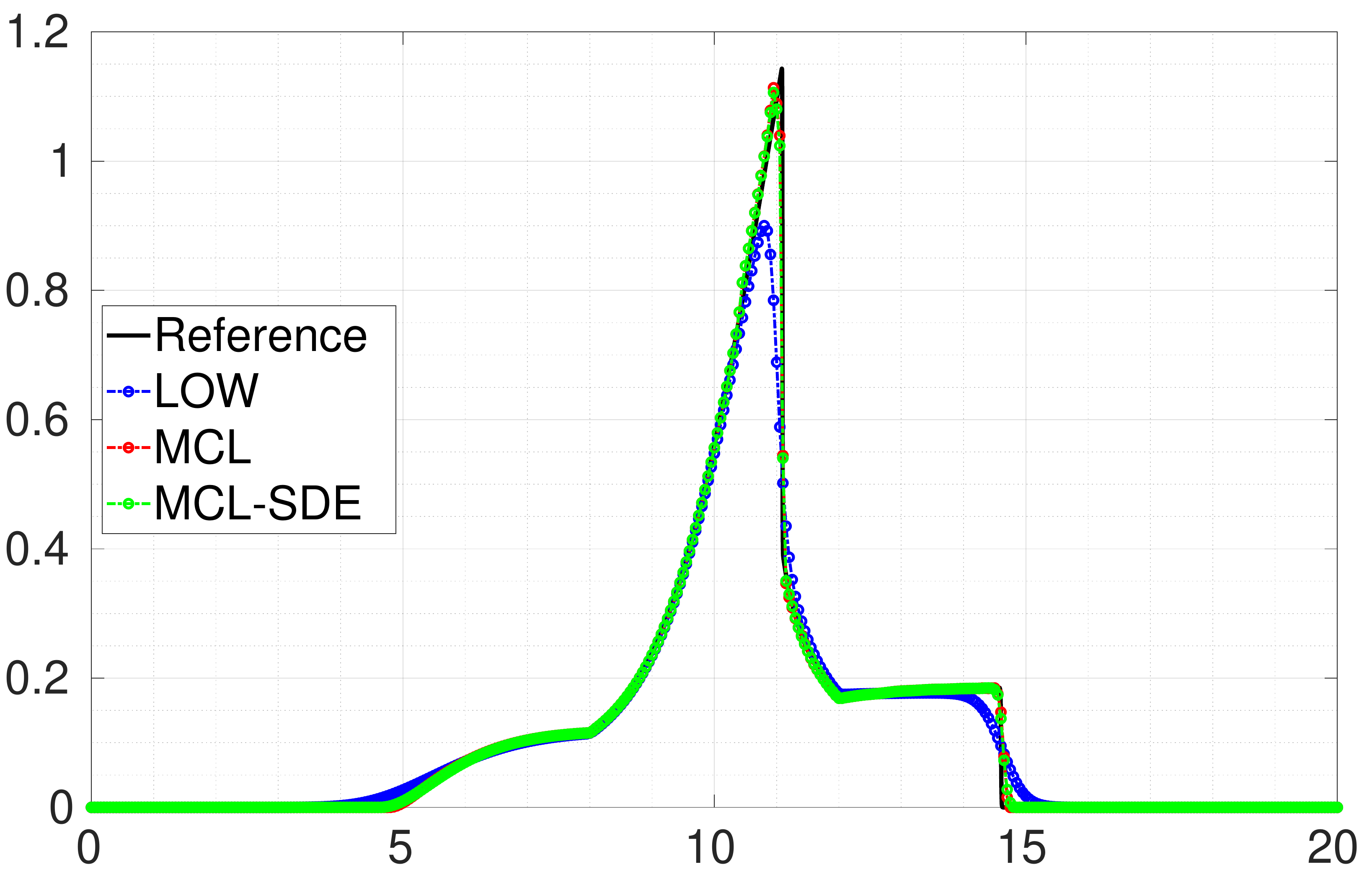}
\end{subfigure}
\caption{Dam break over a bump for the shallow water equations \cite{winters2015}.
Approximations at $T=4.5$ obtained with adaptive SSP2 RK time stepping and $\nu=0.5$ on a uniform mesh consisting of 400 elements.}\label{fig:3-topo}
\end{figure}

To facilitate a comparison of our results with the ones in \cite{winters2015}, we solve this problem up to time $T=4.5$ on a mesh consisting of 400 elements.
A reference solution is obtained with a finite volume method on a fine mesh consisting of $E=10^4$ elements.
Even though the initial water height on the right of the dam is quite small, our friction-based wetting and drying algorithm is never activated in this problem.
The results of this study are displayed in \cref{fig:3-topo}, where we observe excellent agreement with our reference solutions.
The obtained profiles also agree well with the ones in \cite[Sec.~5.6]{winters2015} with the exception that the peaks in the velocity profiles are slightly lower in our results.
This issue requires further investigations and comparisons with the methods in \cite{winters2015}.

\subsection{Oscillating surface in a parabolic lake}

In our final numerical example, we apply our schemes to one of Thacker's oscillatory lakes with a parabolic basin \cite{thacker1981}.
Such benchmarks are challenging tests for wetting and drying algorithms.
We use the same setup as in Vater \etal \cite[Sec.~4.4]{vater2015}, where
$\Omega=(-5000,5000)$, $g=9.81$, and $b(x) = h_0 (x/a)^2$ with $h_0=10$ and $a=3000$.
In the absence of friction, the exact solution is periodic and reads \cite{thacker1981,liang2009,vater2015}
\begin{align*}
x_\pm(t) ={}& -\frac{B}{\omega} \cos(\omega t)\pm a, \qquad B = 5, \qquad \omega = \frac{\sqrt{2gh_0}}{a},\\
H(x,t) ={}& \begin{cases}
h_0 - \frac{B^2}{4g} (1+\cos(2\omega t)) - \frac{Bx}{a} \sqrt{\frac{2h_0}{g}} \cos(\omega t)
& \text{if } x_-(t) \le x \le x_+(t),\\
b(x) & \text{otherwise,}
\end{cases} \\
v(x,t) ={}& \begin{cases}
\frac{Ba\omega}{\sqrt{2h_0g}}\sin(\omega t)
& \text{if } x_-(t) \le x \le x_+(t),\\
0 & \text{otherwise.}
\end{cases}
\end{align*}

\begin{figure}[ht!]
\centering
\begin{subfigure}[b]{0.32\textwidth}
\caption{$t=1000$}
\includegraphics[width=\textwidth]{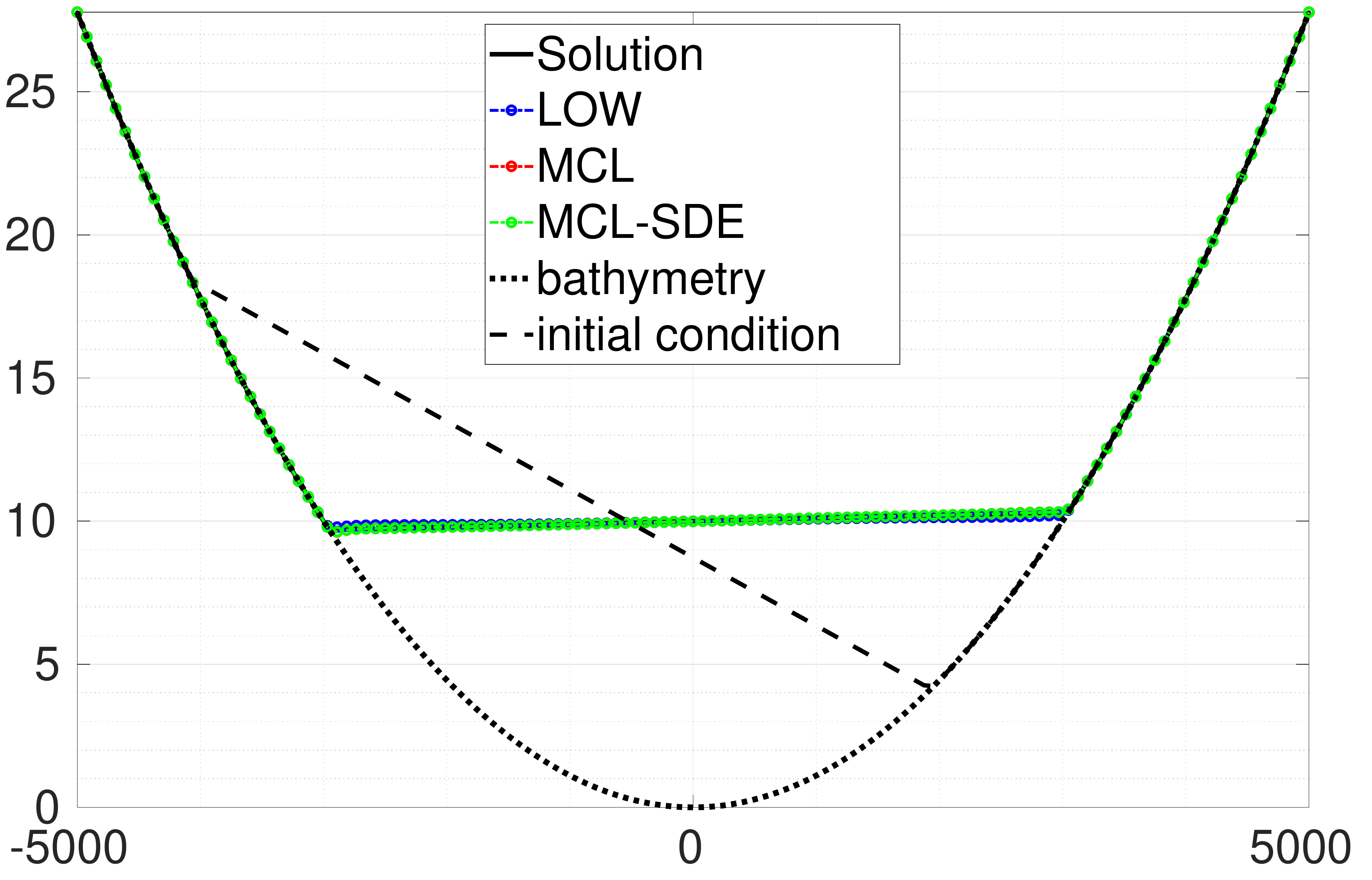}
\end{subfigure}
\begin{subfigure}[b]{0.32\textwidth}
\caption{$t=2000$}\label{fig:lake2}
\includegraphics[width=\textwidth]{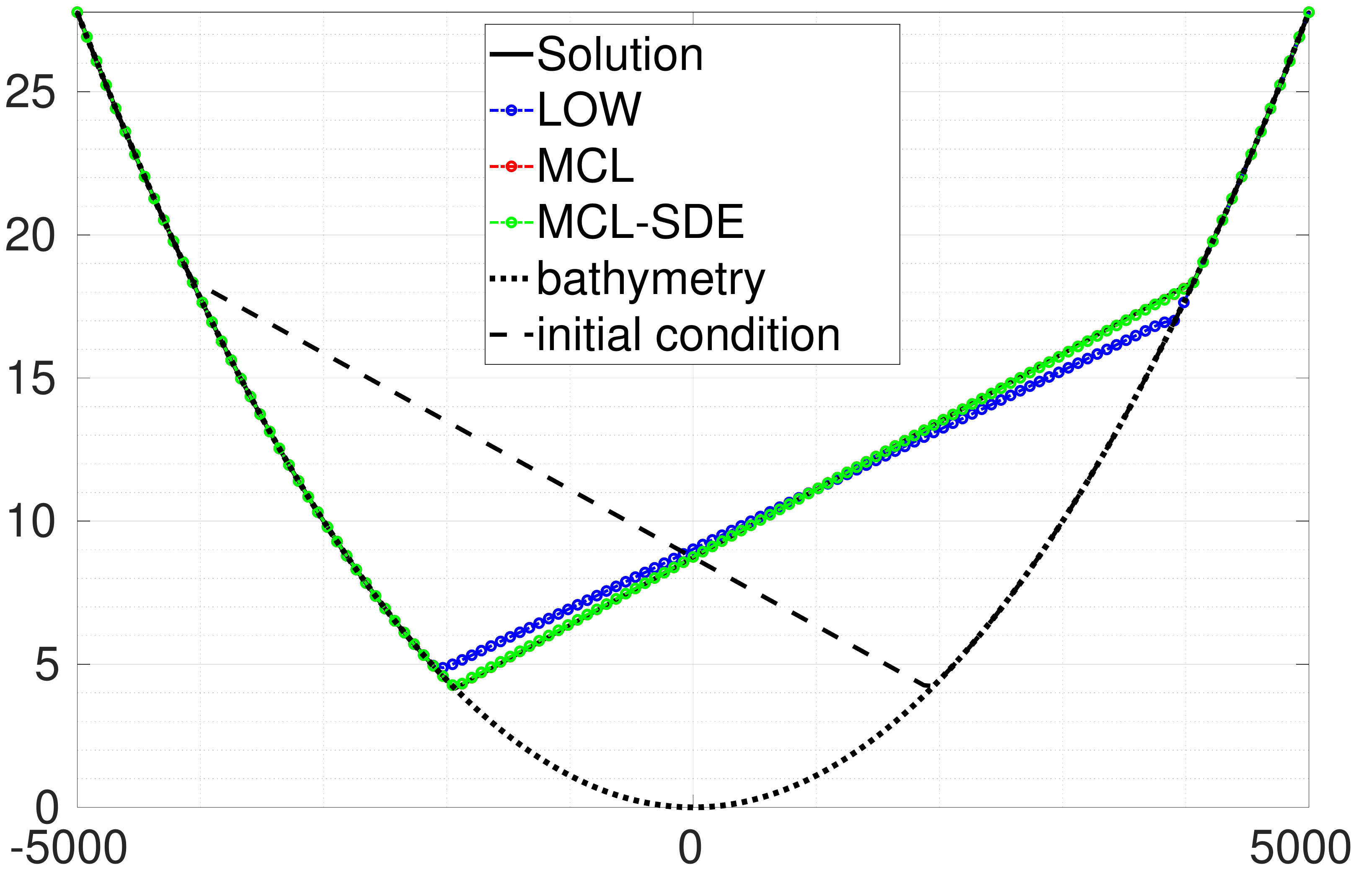}
\end{subfigure}
\begin{subfigure}[b]{0.32\textwidth}
\caption{$t=3000$}
\includegraphics[width=\textwidth]{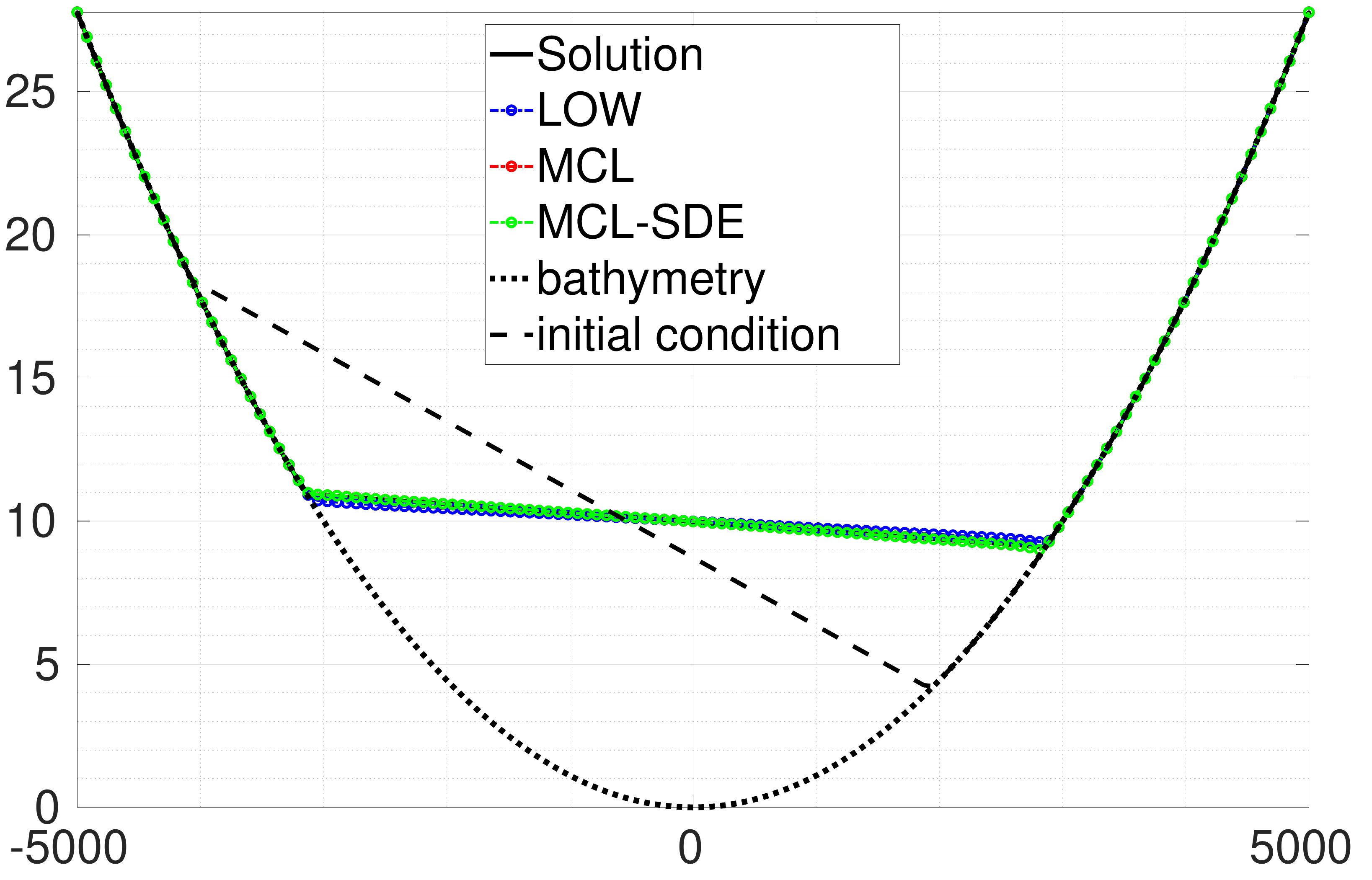}
\end{subfigure}
\caption{Oscillating surface in a parabolic lake for the shallow water equations \cite{thacker1981}.
Approximations to the free surface elevation at various times obtained with adaptive SSP2 RK time stepping and $\nu=0.05$ on a uniform mesh consisting of 128 elements.}\label{fig:lake}
\end{figure}

We employ a CFL parameter of $\nu=0.05$ in combination with our friction-based wetting and drying approach to solve this problem numerically up to end time $T=3000$.
Larger CFL parameters lead to either repetitions of single Runge--Kutta stages or increases of Rusanov diffusion coefficients for nodes around the wet-dry transitions.
For $\nu=0.05$, all schemes remain stable without the need for employing either of these adjustments, even in the case of adaptive time stepping.
\cref{fig:lake} displays LOW, MCL and MCL-SDE water levels at three different times along with the initial condition for illustrative purposes.
From \cref{fig:lake2} we can make out that the low order profile is trailing the exact solution and its flux-corrected counterparts.
Agreement of the flux-limited profiles with the exact water levels is again satisfactory.

We also tested whether we can employ other wetting and drying algorithms in this example.
With the fixes from \cite{ricchiuto2009} and \cite{azerad2017} our simulations crash.
The fix from \cite{kurganov2007a} produces profiles similar to the ones in \cref{fig:lake}.
It is actually possible to employ a larger CFL parameter $\nu$ with this wetting and drying approach.
This observation motivates further tests and adjustments of our friction-based strategy.
Specifically, nonlinear friction models should be considered and the parameters $\delta$ and $\sigma$ may need to be adjusted.
Since these studies should include multidimensional test cases, we have not yet conducted further research in this direction.

\section{Conclusions}

We presented an extension of the bound-preserving and entropy-stable monolithic convex limiting strategy to the inhomogeneous system of shallow water equations.
The proposed scheme is well-balanced \wrt lake at rest equilibria and represents a generalization of the corresponding methods for the SWE without topography.
In addition, we presented two new approaches to handle wet-dry transition regions numerically.
The results of 1D numerical experiments demonstrate the robustness, accuracy, and shock-capturing capabilities of our scheme.
An implementation of the proposed flux-correction schemes in the two-dimensional setting is currently under way to provide additional verification of our new friction-based wetting and drying model.
We are also planning to incorporate additional source terms, such as bottom friction and Coriolis forces, into the model.
Further interesting open problems include achieving well-balancedness \wrt more complicated steady states than the lake at rest, and using high-order baseline discretizations as target schemes.
These aspects require additional research and further generalizations of the proposed methodology.

\section*{Acknowledgment}
This work was supported by the German Research Association (DFG) under grant KU 1530/29-1.

\bibspacing=\dimen 100
\bibliographystyle{bibstyle-article}
\bibliography{bibliography}

\end{document}